\documentclass[12pt]{amsart}
\usepackage{amscd}
\usepackage[dvips]{graphicx}

\newtheorem{thm}{Theorem}[section]

\newtheorem{lma}[thm]{Lemma}

\theoremstyle{definition}
\newtheorem{dfn}[thm]{Definition}
\newtheorem{rmk}[thm]{Remark}

\newenvironment{pf}{\begin{proof}}{\end{proof}}

\numberwithin{equation}{section}

\newcommand{\R}{{\mathbb{R}}}
\newcommand{\C}{{\mathbb{C}}}
\newcommand{\Q}{{\mathbb{Q}}}
\newcommand{\Z}{{\mathbb{Z}}}
\newcommand{\chord}{{\mathcal{R}}}
\newcommand{\V}{{\mathbf{V}}}
\newcommand{\CC}{{\mathbf{C}}}
\newcommand{\A}{{\mathcal{A}}}

\newcommand{\Ordo}{{\mathcal{O}}}
\newcommand{\M}{{\mathcal{M}}}

\newcommand{\pa}{\partial}
\newcommand{\id}{\operatorname{id}}
\newcommand{\img}{\operatorname{im}}
\newcommand{\krn}{\operatorname{ker}}
\newcommand{\pos}{\operatorname{pos}}
\newcommand{\lt}{\textsc{lt}}
\newcommand{\I}{\vert\,}

\title[Rational Symplectic Field Theory\dots]
{Rational Symplectic Field Theory over $\Z_2$\\ for exact Lagrangian cobordisms}
\author{Tobias Ekholm}
\address{USC, Department of mathematics,
3620 S Vermont Ave, Los Angeles, CA 90089}
\email{tekholm{\@@}usc.edu}
\address{Department of mathematics, Uppsala University, Box 480, 751 06 Uppsala, Sweden}
\email{tobias{\@@}math.uu.se}

\begin{document}

\begin{abstract}
We construct a version of rational Symplectic Field Theory for pairs $(X,L)$, where $X$ is an exact symplectic manifold, where $L\subset X$ is an exact Lagrangian submanifold with components subdivided into $k$ subsets, and where both $X$ and $L$ have cylindrical ends. The theory associates to $(X,L)$ a $\Z$-graded chain complex of vector spaces over $\Z_2$, filtered with $k$ filtration levels. The corresponding $k$-level spectral sequence is invariant under deformations of $(X,L)$ and has the following property: if $(X,L)$ is obtained by joining a negative end of a pair $(X',L')$ to a positive end of a pair $(X'',L'')$, then there are natural morphisms from the spectral sequences of $(X',L')$ and of $(X'',L'')$ to the spectral sequence of $(X,L)$. As an application, we show that if $\Lambda\subset Y$ is a Legendrian submanifold of a contact manifold then the spectral sequences associated to $(Y\times\R,\Lambda_k^s\times\R)$, where $Y\times\R$ is the symplectization of $Y$ and where $\Lambda_k^s\subset Y$ is the Legendrian submanifold consisting of $s$ parallel copies of $\Lambda$ subdivided into $k$ subsets, give Legendrian isotopy invariants of $\Lambda$.
\end{abstract}

\thanks{The author acknowledges support from the Alfred P. Sloan Foundation, NSF-grant DMS-0505076, and the Knut and Alice Wallenberg Foundation.}
\maketitle

\section{Introduction}\label{S:intr}
Let $Y$ be a contact $(2n-1)$-manifold with contact $1$-form $\lambda$ (i.e., $\lambda\wedge (d\lambda)^{n-1}$ is a volume form on $Y$). The {\em Reeb vector field} $R_\lambda$ of $\lambda$ is the unique vector field which satisfies $\lambda(R_\lambda)=1$ and $\iota_{R_\lambda}d\lambda=0$, where $\iota$ denotes contraction. The {\em symplectization} of $Y$ is the symplectic manifold $Y\times\R$ with symplectic form $d(e^t\lambda)$ where $t$ is a coordinate in the $\R$-factor. A {\em symplectic manifold with cylindrical ends} is a symplectic $2n$-manifold $X$ which contains a compact subset $K$ such that $X-K$ is symplectomorphic to a disjoint union of two half-symplectizations $Y^+\times\R_+\cup Y^-\times\R_-$, for some contact $(2n-1)$-manifolds $Y^\pm$, where $\R_+=[0,\infty)$ and $\R_-=(-\infty,0]$. We will call (unions of) connected components of $Y^+\times\R_+$  and of $Y^-\times\R_-$ {\em positive- and negative ends} of $X$, respectively, and (unions of) connected components of $Y^+$ and of $Y^-$, {\em $(+\infty)$- and $(-\infty)$-boundaries} of $X$, respectively.

Symplectic Field Theory (SFT) was introduced in \cite{EGH}. It is a framework for extracting invariants in contact- and symplectic geometry by counting punctured $J$-holomorphic curves in symplectic manifolds with cylindrical ends, where $J$ is an almost complex structure adjusted to the symplectic form. The notion adjusted is borrowed from \cite{BEHWZ} where it is shown that if $J$ is adjusted to the symplectic form of a symplectic manifold $X$ with cylindrical ends, then $J$-holomorphic curves are asymptotic to Reeb orbit cylinders in either the positive- or the negative end of $X$ near their punctures and a version of Gromov compactness holds for the moduli space of $J$-holomorphic curves. The SFT of a symplectic manifold with cylindrical ends can be thought of as a theory with several levels. The first level is contact homology where one counts rational curves with one positive puncture, the second level is rational SFT where one counts all rational curves, also those with more than one positive puncture, and the third level is the full SFT where all curves of all genera are taken into account.

The relative counterpart of a symplectic manifold with cylindrical ends is a symplectic $2n$-manifold $X$ with a Lagrangian $n$-submanifold $L$ (i.e., the restriction of the symplectic form in $X$ to any tangent space of $L$ vanishes) such that outside a compact subset, $(X,L)$ is symplectomorphic to the disjoint union of $\bigl(Y^+\times\R_+,\Lambda^+\times\R_+\bigr)$ and $\bigl(Y^-\times\R_-,\Lambda^-\times\R_-\bigr)$, where $\Lambda^\pm\subset Y^{\pm}$ is a Legendrian $(n-1)$-submanifold (i.e., $\Lambda^\pm$ is everywhere tangent to the kernel of the contact form on $Y^\pm$). To define the counterpart of SFT in this setting one should count $J$-holomorphic curves in $X$ with boundary on $L$ asymptotic to Reeb orbit cylinders at interior punctures and to Reeb chord strips of $L$  at boundary punctures. (A Reeb chord is a flow segment of the Reeb vector field in $Y^\pm$ which begins and ends on $\Lambda^\pm$.) Turning such curve counts into algebra is however not straightforward because of a phenomenon often called boundary bubbling. For holomorphic curves without boundary in a symplectic manifold, generic bubbling off is described by the following local model: $\left\{(z,w)\in\C^2\colon z^2\pm w^2 = \epsilon\right\}$, where $\epsilon\in\C$, $\epsilon\to 0$. Hence it is a codimension two phenomenon and can often be disregarded when setting up homology theories. Boundary bubbling for holomorphic curves with boundary on a Lagrangian submanifold can be modeled by a restricted version of this local model as follows. The Lagrangian submanifold corresponds to $\R^2\subset\C^2$, the deformation parameter is constrained to be real, $\epsilon\in\R\subset\C$, and the curve is half of the curve in the model for curves without boundary. Thus, bubbling off at the boundary is a codimension one phenomenon which cannot be disregarded when setting up homology theories. In Lagrangian Floer homology, techniques for dealing with boundary bubbling have been developed in \cite{FOOO} and in \cite{CL}, see the discussion after Theorem \ref{t:main2}.

Let $(X,L)$ be a pair with properties as described above. If the symplectic manifold $X$ is exact (i.e., if the symplectic form $\omega$ on $X$ satisfies $\omega=d\beta$ for some $1$-form $\beta$) and if the Lagrangian submanifold $L$ is exact as well (i.e., if the restriction $\beta|L$ satisfies $\beta|L=df$ for some function $f$), then Stokes theorem implies that every $J$-holomorphic curve in $X$ with boundary on $L$ must have at least one positive puncture. For simpler notation, we call such a pair $(X,L)$ of exact manifolds an {\em exact cobordism}. For exact cobordisms, the extension of the first level of SFT to the relative case meets with no conceptual difficulties because boundary bubbling cannot appear: one of the holomorphic pieces after boundary bubbling would not have any positive puncture. The resulting theory gives a rich set of invariants, see e.g. \cite{Ch, El, EES1, EES3, Ng1, Ng2, K}. The aim of the present paper is to provide a generalization of the second level of SFT to the relative case of exact cobordisms with cylindrical ends.

In order to formulate our main results we introduce the following notation. Let $(X,L)$ be an exact cobordism. Let $(Y^\pm\times\R_{\pm},\Lambda^\pm\times\R_\pm)$ be the ends of $(X,L)$ and write $(\bar X,\bar L)$ for the compact part of $(X,L)$ obtained by cutting the infinite parts of the cylindrical ends off at $Y^{\pm}\times\{0\}$. We will sometimes think of Reeb chords of $\Lambda^\pm\subset Y^{\pm}$ in the $(\pm\infty)$-boundary of $(X,L)$ as lying in $\pa\bar X$ with endpoints on $\pa\bar L$. A formal disk in $(X,L)$ is a homotopy class of maps of the $2$-disk $D$, with $m$ marked disjoint closed subintervals in $\pa D$, into $\bar X$, where the $m$ marked intervals are required to map in an orientation preserving (reversing) manner to Reeb chords of $\pa \bar L$ in the $(+\infty)$-boundary (in the $(-\infty)$-boundary) and where remaining parts of the boundary $\pa D$ map to $\bar L$. The {\em action} of a Reeb chord is the integral of the contact form in the $(\pm\infty)$-boundary along it and the {\em $(+)$-action} of a formal disk is the sum of the actions of its Reeb chords in the $(+\infty)$-boundary. Assume that the set of connected components of $L$ has been subdivided into subsets (we call such subsets {\em pieces}). In Subsection \ref{s:admformdisk}, we define the notion of an {\em admissible} formal disk in $(X,L)$ in this situation. A $J$-holomorphic disk in $X$ with boundary on $L$ determines a formal disk and if this formal disk is admissible then boundary bubbling is impossible for topological reasons.

Let $\V(X,L)$ denote the $\Z$-graded vector space over $\Z_2$ consisting of all formal sums of admissible formal disks of $(X,L)$ which contain only a finite number of summands of $(+)$-action below any given number and with grading as follows: the degree of a formal disk is the formal dimension of the moduli space of $J$-holomorphic disks homotopic to the formal disk. The degree of a formal disk is thus computed in terms of the Maslov index of the boundary condition it determines and the number of Reeb chords it has, see Subsection \ref{s:grad}. We use the filtration $0\subset F^k\V(X,L)\subset\dots\subset F^2\V(X,L)\subset F^1\V(X,L)=\V(X,L)$, where $k$ is the number of pieces of $L$ and where the filtration level is determined by the number of Reeb chords of a formal disk which are in the $(+\infty)$-boundary of $(X,L)$, see Subsection \ref{s:filt}. In Subsection \ref{l:diffconc}, using the collections of all rigid admissible $J$-holomorphic disks in $(X,L)$ and of all $1$-parameter families of admissible $J$-holomorphic disks in the symplectizations of the contact manifolds in the $(\pm\infty)$-boundary of $(X,L)$, we define a differential $d^f\colon \V(X,L)\to\V(X,L)$. The differential increases grading by $1$, respects the filtration, and is $(+)$-action non-decreasing. Let $\alpha>0$. The finite $(+)$-action subspace $\V_{[\alpha]}(X,L)\subset\V(X,L)$ is the subspace of formal sums of admissible disks of $(+)$-action at most $\alpha$. It inherits a filtration from the filtration on $\V(X,L)$. Since $d^f$ does not decrease $(+)$-action, it induces filtration preserving differentials $d^f_\alpha$ on $\V_{[\alpha]}(X,L)$ and the natural projection maps $\pi^\alpha_\beta\colon\V_{[\alpha]}(X,L)\to\V_{[\beta]}(X,L)$, $\alpha>\beta$, are filtration preserving chain maps. Define the {\em rational admissible SFT spectral sequence} $\bigl\{E^{p,q}_{r\,;[\alpha]}(X,L)\bigr\}_{r=1}^k$ as the spectral sequence induced by the filtration preserving differentials $d^f_\alpha\colon\V_{[\alpha]}(X,L)\to\V_{[\alpha]}(X,L)$, see Subsection \ref{s:spsq}.

Throughout the paper we will assume that all exact cobordisms have {\em good ends}. The details of this technical condition are specified in Subsection \ref{s:goodends} of Appendix \ref{A:B}. It restricts the contact forms on the contact manifolds at the $(\pm\infty)$-boundaries of exact cobordisms, but allows for certain non-compact manifolds as ends, in particular $1$-jet spaces.
\begin{thm}\label{t:main1}
Let $(X,L)$ be an exact cobordism with good ends and with a subdivision $L=L_1\cup\dots\cup L_k$ of $L$ into pieces and let $\alpha>0$. Then $\bigl\{E^{p,q}_{r\,;[\alpha]}(X,L)\bigr\}$ does not depend on the choice of adjusted almost complex structure $J$, and is invariant under compactly supported exact deformations of $(X,L)$.
\end{thm}

Theorem \ref{t:main1} is proved in Section \ref{S:proofs}. The spectral sequence in Theorem \ref{t:main1} share many of the familiar properties of rational SFT in the non-relative case, see \cite{EGH}. For example, if a pair $(Y,\Lambda)$ of a contact manifold with a Legendrian submanifold is a $(+\infty)$-boundary of an exact cobordism $(X^b,L^b)$ and a $(-\infty)$-boundary of another exact cobordism $(X^a,L^a)$ then the two cobordisms can be joined along $(Y,\Lambda)$ to form a new exact cobordism $(X^{ba},L^{ba})$, and there are filtration preserving chain maps $\V(X^b,L^b)\to\V(X^{ba},L^{ba})$ and $\V(X^a,L^a)\to\V(X^{ba},L^{ba})$ which induce morphisms of spectral sequences, see Lemma \ref{l:chmap}.

If $\Lambda\subset Y$ is a Legendrian submanifold of a contact manifold $Y$ and if $\Lambda$ is subdivided into pieces $\Lambda=\Lambda_1\cup \Lambda_2\cup\dots\cup \Lambda_k\subset Y$ then $(Y\times\R,\Lambda\times\R)$ is an exact cobordism with $\Lambda\times\R$ subdivided into pieces. In Subsection \ref{d:spsqLeg} we define the {\em rational admissible SFT-invariant of $(Y,\Lambda)$}. In general, this invariant is the graded vector space $\bigl\{E^{p,q}_r(Y,\Lambda)\bigr\}_{r=1}^k$ which is the inverse limit of the inverse system
$$
\begin{CD}
\bigl\{E^{p,q}_{r;[\beta]}(Y\times\R,\Lambda\times\R)\bigr\}
@<{{\pi^\alpha_\beta}_\ast}<<
\bigl\{E^{p,q}_{r;[\alpha]}(Y\times\R,\Lambda\times\R)\bigr\}
\end{CD},
\quad 0<\beta<\alpha.
$$
In the special case that $\Lambda\subset Y$ has only finitely many Reeb chords (e.g., if $Y=J^1(M)$ is the $1$-jet space of some smooth manifold $M$) then the limit $\bigl\{E^{p,q}_r(Y,\Lambda)\bigr\}$ has a naturally induced structure as a spectral sequence which is also part of the invariant, see Subsection \ref{s:spsq}.

\begin{thm}\label{t:main2}
If $\Lambda\subset Y$ is a Legendrian submanifold of a contact manifold $Y$, with a good contact form, and if $\Lambda$ is subdivided into pieces then $\bigl\{E^{p,q}_r(Y,\Lambda)\bigr\}$ is independent of good contact form and invariant under contact isotopies of the pair $(Y,\Lambda)$ (in particular under Legendrian isotopies of $\Lambda$).
\end{thm}

Theorem \ref{t:main2} is proved in Section \ref{S:proofs}. We give a brief description of how it can be applied. This also demonstrates that although the strategy for dealing with boundary bubbling in this paper is different from that of \cite{FOOO} and that of \cite{CL}, the constructions are related. If $\Lambda\subset Y$ is a Legendrian submanifold (possibly connected) of a contact manifold, then we consider the many component Legendrian submanifold $\tilde\Lambda$ consisting of finitely many nearby parallel (i.e., parallel along the Reeb flow) copies of $\Lambda$. Partitioning the collection of parallel copies into pieces we can apply the rational admissible SFT invariant. In fact, much like in \cite{EENS}, the differential on $\V(Y\times\R,\tilde\Lambda\times\R)$ can in this case be computed in terms of (all) moduli spaces of holomorphic disks in $Y\times\R$ with boundary on $\Lambda\times\R$ in combination with the spaces of Morse flow trees in $\Lambda$, see \cite{FO,E1}, determined by a finite collection of Morse functions on $\Lambda$. In this sense, the method for dealing with boundary bubbling in the present paper is related to \cite{FOOO} and to \cite{CL}, though in \cite{CL}, only one Morse function and flow lines appears, rather than as here and in \cite{FOOO}, many Morse functions and flow trees. We will address the exact expression of the rational admissible SFT invariant in terms of Morse flow trees elsewhere.

In \cite{E2}, we give an application of Theorem \ref{t:main2}. It is used to show that the members of a basic family of $2n$-dimensional Legendrian spheres in $J^1(\R^{2n})$, all with the same classical invariants, and all with vanishing contact homology, are pairwise non-Legendrian isotopic. We plan to describe a more complete version of admissible rational SFT for exact cobordisms, removing the condition of ends being good, including  oriented moduli spaces (which requires relatively spin Legendrian submanifolds) and working with vector spaces over $\Q$, in a future paper.

The paper is organized as follows. In Section \ref{S:vsp+op}, the vector space underlying the chain complex associated to an exact cobordism is defined and some operations on this vector space are discussed. In Section \ref{S:grafildiff}, the grading, the filtration, and the differential of the chain complex is introduced, the corresponding spectral sequences are defined, and induced chain maps are studied. In Section \ref{S:inv}, lemmas necessary to demonstrate invariance properties of the spectral sequences are proved and in Section \ref{S:proofs}, these results are collected into proofs of the main results. In Appendix \ref{A:A}, a basic construction in symplectic geometry is discussed and in Appendix \ref{A:B}, properties of holomorphic disks in exact cobordisms are described and a perturbation scheme for transversality is constructed using the polyfold language developed in \cite{H, HWZ}.

\subsection*{Acknowledgements} The author thanks Y. Eliashberg, K. Honda, and L. Ng for useful discussions.
\section{The vector space of an exact cobordism}\label{S:vsp+op}
In this section we introduce admissible formal disks. Using them, we associate a graded filtered vector space to an exact cobordism. The section is organized as follows. In Subsection \ref{s:admformdisk}, admissible formal disks are defined and in Subsection \ref{s:compdisk}, gluings of such disks in two exact cobordisms joined at a common end are studied. In Subsection \ref{s:vsp+op}, the vector space of an exact cobordism is defined and gluings of formal disks are interpreted as operations on such vector spaces.

\subsection{Admissible formal disks}\label{s:admformdisk}
Let $(X,L)$ be an exact cobordism with $(\pm\infty)$-boundary $(Y^\pm,\Lambda^\pm)$. Recall the subdivision $(X,L)=(\bar X,\bar L)\cup (Y^+\times\R_+,\Lambda^+\times\R_+)\cup(Y^-\times\R_-,\Lambda^-\times\R_-)$ into a compact part and ends. Make the identification $(\pa\bar X,\pa\bar L)=(Y^+,\Lambda^+)\cup(Y^-,\Lambda^-)$. Assume that $L$ comes equipped with a subdivision $L=L_1\cup\dots\cup L_k$ into pieces where each piece $L_j$ is a union of connected components of $L$. This subdivision induces a subdivision of the ends, $\Lambda^\pm=\Lambda^\pm_1\cup\dots\cup\Lambda^{\pm}_{k}$. Let $\chord^+$ and $\chord^-$ denote the sets of Reeb chords of $\Lambda^+$ and $\Lambda^-$, respectively. We write $\chord^\pm=\chord^\pm_{\rm pu}\cup\chord^{\pm}_{\rm mi}$. Here $\chord^{\pm}_{\rm pu}$ contains all {\em pure} Reeb chords with both endpoints on the same piece in the subdivision of $\Lambda^\pm$, and $\chord^{\pm}_{\rm mi}$ contains all {\em mixed} Reeb chords with endpoints on distinct pieces. Note that a Reeb chord is oriented (by the Reeb flow).

A {\em formal disk map} is a map from a $2$-disk $D$ with $m$ marked disjoint boundary segments into $\bar X$ with the following properties. Each marked boundary segment either maps in an orientation preserving way to a Reeb chord of $\Lambda^+\subset\pa\bar L$ in $Y^+\subset\pa\bar X$, or maps in an orientation reversing way to a Reeb chord of $\Lambda^-\subset\pa\bar L$ in $Y^-\subset\pa\bar X$. Each unmarked boundary segment maps to $\bar L$. We say that two formal disk maps are {\em homotopic} if they are homotopic through formal disk maps. In particular, two formal disk maps are homotopic only if they have the same Reeb chords and the respective induced cyclic orderings on these Reeb chords agree. Furthermore, if two formal disk maps have the same Reeb chords in the same cyclic order then their unmarked boundary arcs which connect corresponding Reeb chord endpoints form difference-loops (i.e., the path of one disk followed by the inverse path of the other) in $\bar L$ and the two formal disk maps are homotopic only if all these difference-loops are contractible. Finally, in case all the difference-loops are contractible, choosing homotopies between the boundary loops, the two formal disk maps together with these homotopies give a difference-map of a $2$-sphere into $\bar X$ and the formal disk maps are homotopic if and only if the homotopy class of this difference-map lies in the image of $\pi_2(\bar L')\to\pi_2(\bar X)$, where $\bar L'=\bar L_1'\cup\dots\cup\bar L_m'\subset\bar L$ is the union of the connected components of $\bar L$ which contain some boundary component of the formal disk.

In conclusion, if for a fixed cyclic word of Reeb chords there is a formal disk map which realizes this word then the homotopy classes of formal disk maps is a principal homogeneous space over the product of the kernel of a map $\pi_1(\bar L_1')\times\dots\times\pi_1(\bar L_m')\to\pi_1(\bar X)$ determined by the Reeb chord endpoints and the quotient $\pi_2(\bar X)/\img(\pi_2(\bar L')\to\pi_2(\bar X))$. A {\em formal disk} is a homotopy class of formal disk maps.

We call the Reeb chords of the formal disk its {\em punctures}. When speaking of formal disks we will contract the marked intervals which map to the Reeb chords to points and call them punctures as well. A component of the complement of the punctures in the boundary of a formal disk will be called a {\em boundary component}. We say that a puncture of a formal disk is {\em positive} if it maps to a chord in $\chord^+$ and that it is {\em negative} if it maps to a chord in $\chord^-$. We say that a puncture of a formal disk is {\em mixed} if it maps to a chord in $\chord^\pm_{\rm mi}$ and that it is {\em pure} if it maps to a chord in $\chord^\pm_{\rm pu}$.

Let $D$ be the source of a formal disk. A neat arc in $D$ which connects two of its boundary components (and which is transverse to the boundary at the boundary) will be called a {\em collapsing arc}. Note that a collapsing arc $a$ subdivides $D$ into two components $D_1(a)$ and $D_2(a)$.

\begin{dfn}\label{d:adm}
A formal disk parameterized by $D$ is {\em admissible} if it meets the following two conditions.
\begin{itemize}
\item[$({\bf a1})$] $D$ has at least one puncture which is positive.
\item[$({\bf a2})$] For every collapsing arc $a$ in $D$ with endpoints on boundary components mapping to the same piece of $\bar L$, one component $D_1(a)$ or $D_2(a)$ has either no punctures or only pure negative punctures.
\end{itemize}
\end{dfn}
\begin{rmk}
Definition \ref{d:adm} is not the only possible choice. All technical results proved below hold true if one uses the following alternative definition instead:
keep $({\bf a2})$ and change $({\bf a1})$ to
\begin{itemize}
\item[$({\bf a1'})$] $D$ has at least one puncture which is positive or mixed.
\end{itemize}
This gives rise to a somewhat different theory discussed in \cite{E2}.
\end{rmk}

\begin{rmk}
Another possibility for an alternative definition of admissibility is as follows.
Keep $({\bf a1})$ and change $({\bf a2})$ to
\begin{itemize}
\item[$({\bf a2\,'})$] For every collapsing arc $a$ in $D$ with endpoints on boundary components mapping to the same piece of $\Lambda$, one component $D_1(a)$ or $D_2(a)$ has either no punctures or only negative punctures.
\end{itemize}
However, to prove the invariance of the theory which arises from this definition it seems, because of problems with ordering of punctures, one must introduce a distinguished puncture in each disk and average, and hence work over $\Q$ (instead of over $\Z_2$) and consider orientations. It would be interesting to understand the relation between the Legendrian isotopy invariants which arise from $({\bf a2\,'})$ and the collection of such invariants which arises from Definition \ref{d:adm}, when the connected components of a given Legendrian submanifold are reorganized into pieces in all possible ways.
\end{rmk}

\begin{lma}\label{l:admpurepos}
An admissible disk with a pure positive puncture can neither have mixed punctures nor have other positive punctures.
\end{lma}
\begin{pf}
If it did, a small collapsing arc cutting the pure positive puncture off from the rest of the disk contradicts $({\bf a2})$.
\end{pf}

\subsection{Joining exact cobordisms and gluing formal disks}\label{s:compdisk}
Let $(X^b,L^b)$ be an exact cobordism with a $(+\infty)$-boundary equal to $(Y,\Lambda)$ and let $(X^a,L^a)$ be an exact cobordism with a $(-\infty)$-boundary equal to $(Y,\Lambda)$. Then we can join the exact cobordisms over $(Y,\Lambda)$ to an exact cobordism $(X^{ba},L^{ba})$. Since $\Lambda\subset\pa {\bar L}^b$ and $\Lambda\subset\pa {\bar L}^a$, we may view $\Lambda$ as a submanifold of $L^{ba}\subset X^{ba}$. A {\em partial formal disk} of $(X^{ba},L^{ba})$ is defined as a formal disk except that it is allowed to have punctures also at Reeb chords of $\Lambda\subset\bar L^{ba}\subset \bar X^{ba}$.

Let $\chord$ denote the set of Reeb chords of $\Lambda$. Let $g^b$ be a formal disk in $(X^b,L^b)$ with one of its positive punctures at a Reeb chord $c\in \chord$. Let $g^a$ be a formal disk in $(X^a,L^a)$ with one of its negative punctures at $c$. Then we can attach $g^b$ to $g^a$ at $c$ in an obvious way, to form a partial formal disk $g^{ba}_1$ in $(X^{ba},L^{ba})$. This construction can be repeated: at a positive- or negative puncture of $g^{ba}_1$ which lies in $\chord$, a formal disk of $(X^a,L^a)$ or $(X^b,L^b)$, respectively, can be attached to $g^{ba}_1$ to form a new partial formal disk $g^{ba}_2$, etc. We say that the formal disks in $(X^a,L^a)$ and $(X^b,L^b)$ used to build a partial formal disk in $(X^{ba},L^{ba})$ in this way are its {\em factors}.

\begin{lma}\label{l:nonadmsub}
Let $w$ be a formal disk in $(X^{ba},L^{ba})$ with factors which are formal disks in $(X^a,L^a)$ and $(X^b,L^b)$, all with at least one mixed- or positive puncture. If some partial formal sub-disk of $w$ is non-admissible then so is $w$.
\end{lma}

\begin{pf}
No formal disk has only one mixed puncture. Consider a collapsing arc of the sub-disk with mixed- or positive punctures on both sides. Attaching at a mixed puncture leaves a mixed puncture and attaching at a positive puncture leaves a positive puncture. Thus, the collapsing arc of the sub-disk shows also that the final formal disk is non-admissible.
\end{pf}

\subsection{The vector space of formal disks and operations}\label{s:vsp+op}
If $\gamma$ is curve in a contact manifold $Y$ with contact form $\lambda$, then the {\em action} of $\gamma$ is
$$
\A(\gamma)=\int_\gamma\lambda.
$$
Let $\Lambda\subset Y$ be a Legendrian submanifold. Reeb chords of $\Lambda$ are critical points for the action functional acting on curves with endpoints on $\Lambda$. If the contact form $\lambda$ is chosen sufficiently generic then the Reeb chords of $\Lambda$ are isolated and, if $\chord$ is the set of Reeb chords of $\Lambda$, then the set $\{\A(c)\colon c\in\chord\}$ is a discrete subset of $\R$ which is bounded below by some number $\A_0(Y,\Lambda)>0$. (The existence of a lower bound is true in greater generality and holds since the Reeb vector field is never tangent to a Legendrian submanifold).

Let $(X,L)$ be an exact cobordism with $(\pm\infty)$-boundary $(Y^\pm,\Lambda^\pm)$. If $g$ is a formal disk in $(X,L)$ with positive punctures at Reeb chords $a_1,\dots,a_p$ and negative punctures at Reeb chords $b_1,\dots, b_q$ then define the {\em $(+)$-action} of $g$ as
\begin{equation}\label{e:posaction}
\A^+(g)=\sum_{j=1}^p\A(a_j),
\end{equation}
{\em the $(-)$-action} of $g$ as
\begin{equation}\label{e:negaction}
\A^-(g)=\sum_{j=1}^q\A(b_j),
\end{equation}
and the {\em action} of $g$ as
\begin{equation}\label{e:action}
\A(g)=\A^+(g)-\A^-(g).
\end{equation}

Let $\tau=(\tau_1,\dots,\tau_n)$ be a vector of formal variables. Define the module $\V(X,L;\tau)$ over the polynomial ring $\Z_2[\tau]$ as the vector space over $\Z_2$ of all formal $\Z_2[\tau]$-sums of admissible formal disks of $(X,L)$ which satisfy the following finiteness condition: for any $\alpha>0$, the number of formal disks in a formal sum of $(+)$-action at most $\alpha$ is finite.

Using formulas, $\V(X,L;\tau)$ can be described as follows. Let $G$ denote the set of formal disks in $(X,L)$. An element $v\in\V(X,L;\tau)$ can be written as
\begin{equation}\label{e:sumrep}
v=\sum_{g\in G} p_g(\tau)\,g,
\end{equation}
where $p_g(\tau)\in\Z_2[\tau]$ and where for any $\alpha>0$ the set
$$
\bigl\{g\in G\colon \A^+(g)\le \alpha,\, p_g(\tau)\ne 0\bigr\}
$$
is finite. As indicated in \eqref{e:sumrep}, we think of a vector $v\in\V(X,L;\tau)$ as a collection of basis elements (denoted $g$) with non-zero weights in $\Z_2[\tau]$ (denoted $p_g(\tau)$). Addition in $\V(X,L;\tau)$ is defined by componentwise addition: if
$$
v=\sum_{g\in G} p_g(\tau)\,g\quad\text{and}\quad
w=\sum_{g\in G} q_g(\tau)\,g
$$
then
$$
v+w=\sum_{g\in G} \bigl(p_g(\tau)+q_g(\tau)\bigr)\,g.
$$
Multiplication with a polynomial is also defined componentwise: if $q(\tau)\in\Z_2[\tau]$ then
$$
q(\tau)\, v =
q(\tau)\left(\sum_{g\in G} p_g(\tau)\,g\right)=
\sum_{g\in G} q(\tau)p_g(\tau)\,g.
$$
Define $\V^+(X,L;\tau)\subset\V(X,L;\tau)$ as the submodule of formal sums in which all formal disks have positive action.

Let $\alpha>0$. Define the finite $(+)$-action submodules $\V_{[\alpha]}(X,L)\subset\V(X,L)$ as the $\Z_2[\tau]$-module generated by all formal disks $g$ with $\A^+(g)\le\alpha$.
\begin{rmk}\label{r:invlim}
If $\alpha>\beta>0$ then there are natural projection maps
$$
\begin{CD}
\V_{[\beta]}(X,L;\tau) @<\pi^\alpha_\beta<< \V_{[\alpha]}(X,L;\tau)
\end{CD}
$$
where $\pi^\alpha_\beta\colon\V_{[\alpha]}(X,L)\to\V_{[\alpha]}(X,L)/K_\beta^\alpha\approx\V_{[\beta]}(X,L)$, with $K^\alpha_\beta$ the submodule of $\V_{[\alpha]}(X,L)$ generated by all formal disks $g$ with $\A^+(g)>\beta$, and  $\V(X,L;\tau)$ is expressible as the inverse limit
$$
\V(X,L;\tau)=\lim_{\longleftarrow\, \alpha} \V_{[\alpha]}(X,L;\tau).
$$
\end{rmk}

\subsection{Gluing pairing}\label{s:glupair}
As in Subsection \ref{s:compdisk}, let $(X^b,L^b)$ and $(X^a,L^a)$ be exact cobordisms such that $(Y,\Lambda)$ is a $(+\infty)$-boundary of $(X^b,L^b)$ and a $(-\infty)$-boundary of $(X^a,L^a)$, and let $(X^{ba},L^{ba})$ be the joined cobordism. We define gluing pairings
$$
\V(X^b,L^b;\tau^b)\times\V(X^a,L^a;\tau^a)
\to\V(X^{ba},L^{ba};\tau^a,\tau^b),
$$
as follows. If $v^b\in\V(X^b,L^b;\tau^b)$ and $v^a\in\V(X^a,L^a;\tau^a)$ then the {\em gluing pairing} of $v^b$ and $v^a$
\begin{equation}\label{e:glupair}
\left(v^a\I v^b\right)\in\V(X^{ab},L^{ba};\tau^a,\tau^b)
\end{equation}
is the vector of all admissible formal disks, with factors from $v^a$ and $v^b$, weighted by the product of the weights of its factors. In order to make this definition precise we discuss properties of the gluing operation and explain how to count glued disks.

\begin{lma}\label{l:factor1}
If $u$ and $v$ are partial formal disks in $(X^{ab},L^{ab})$ and if $w$ is the partial formal disk obtained by joining $u$ and $v$ at a Reeb chord of $(Y,\Lambda)$ then
$$
\A^+(w)\ge \min\{\A^+(u),\A^+(v)\}.
$$
In particular, if $g^{ab}$ is a formal disk in $(X^{ba},L^{ba})$ with factors $\{g^a_j\}_{j=1}^s$ and $\{g^b_k\}_{k=1}^t$, formal disks in $(X^a,L^a)$ and in $(X^b,L^b)$, respectively, then
$$
\A^+(g^{ba})\ge\min\left\{\A^+(g^a_j), \A^+(g^b_k)\right\}_{1\le j\le s,\, 1\le k\le t}.
$$
\end{lma}
\begin{pf}
Either all Reeb chords at the positive punctures of $u$ or all Reeb chords at the positive punctures of $v$ are Reeb chords at the positive punctures of $w$ as well.
\end{pf}

\begin{lma}\label{l:no2mixed}
Let $u$ be a mixed formal disk in $(X^b,L^b)$ or in $(X^a,L^a)$ with at least one positive puncture. Consider a formal disk  $w$ in $(X^{ba},L^{ba})$ with factors which are formal disks in  $(X^b,L^b)$ or in $(X^a,L^a)$. Assume that all factors of $w$ have at least one positive puncture and that $w$ has at least two factors which both equal $u$. Then $w$ is non-admissible.
\end{lma}

\begin{figure}[htbp]
\begin{center}
\includegraphics[angle=0, width=10cm]{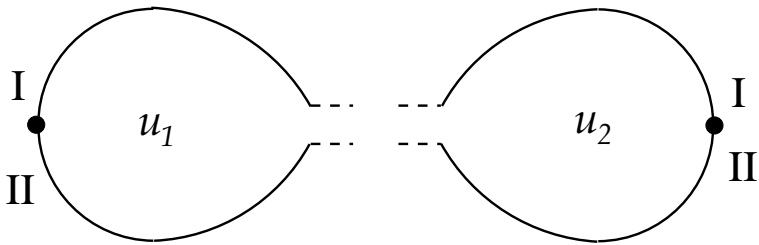}
\end{center}
\caption{${\rm I}$ and ${\rm II}$ label pieces of $L^{ba}$. If the collapsing arc connecting ${\rm I}$ to ${\rm I}$ does not have mixed punctures on both sides then the collapsing arc connecting ${\rm II}$ to ${\rm II}$ does.}
\label{fig:collarc}
\end{figure}

\begin{pf}
The formal disk $w$ has a tree structure where its factors are the nodes and glued punctures are edges. Pick the shortest path in the tree of $w$ which connects two distinct $u$-factors and consider the corresponding sub-disk $w'$ of $w$. Lemma \ref{l:nonadmsub} implies that it is enough to show that $w'$ is non-admissible. We write $w'=w''\sharp\, u_1\sharp\, u_2$, where $w''$ is the sub-disk of $w'$ obtained by removing the two $u$-factors $u_1$ and $u_2$.

If there exists some mixed puncture $p$ of $u$ such that the corresponding punctures $p_1$ in $u_1$ and $p_2$ in $u_2$ are both punctures of $w'$, then an arc in $w'$ starting on the incoming boundary component near $p_1$ and ending at the incoming boundary component near $p_2$ gives a collapsing arc separating $p_1$ from $p_2$ and hence $w'$ is non-admissible. If $u$ is a formal disk in $(X^a,L^a)$ then such a puncture always exists: take $p$ equal to any of the positive punctures of $u$. If $u$ is a formal disk in $(X^b,L^b)$ and if there is no such puncture $p$, then $u$ must have two mixed positive punctures and no other mixed punctures. Since $w''$ must contain some $(X^a,L^a)$-disk, $w'$ has at least one more mixed positive puncture.  (For example, the positive puncture of an $(X^a,L^a)$-factor, which must be mixed since it connects to mixed disks, see Lemmas \ref{l:admpurepos} and \ref{l:nonadmsub}.) At least one of the arcs in $w'$ connecting matching boundary components of $u_1$ and $u_2$ is then a collapsing arc with mixed punctures on both sides showing that $w'$ is non-admissible, see Figure \ref{fig:collarc}.
\end{pf}

\begin{lma}\label{l:factor2}
Let $v^a\in\V(X^a,L^a)$ and $v^b\in\V(X^b,L^b)$.
For any $\alpha>0$, there are only finitely many admissible formal disks $g^{ba}$ in $(X^{ba},L^{ba})$ with $\A^+(g^{ba})\le \alpha$ and with factors from $v^a$ and $v^b$. Furthermore each such admissible formal disk $g^{ba}$ has only finitely many factorizations into factors from $v^a$ and $v^b$.
\end{lma}

\begin{pf}
Lemma \ref{l:factor1} implies that for any factor $g$ of $g^{ba}$ as above $\A^+(g)\le \alpha$. Thus there are only finitely many formal disks from $v^a$ and $v^b$ which can arise as factors of $g^{ba}$. Lemma \ref{l:no2mixed} implies that any mixed formal disk in $v^a$ or $v^b$ can only appear once as a factor in $g^{ba}$. If a factor from $v^a$ is pure then all factors of $g^{ba}$ are pure by admissibility. Hence the total number of negative punctures of disks in $v^a$ where pure factors from $v^b$ can be attached is finite. The lemma follows.
\end{pf}

We use Lemmas \ref{l:factor1} and \ref{l:factor2} to define the gluing pairing \eqref{e:glupair} as follows. Let $v^b\in\V(X^b,L^b;\tau^b)$ and $v^a\in\V(X^a,L^a;\tau^a)$. In order to find the coefficient of one of the finitely many formal disks $g^{ba}$ of $(+)$-action below $\alpha>0$ which may contribute to $\left(v^b\I v^a\right)$ we need only check the finite parts $\bar v^a$ of $v^a$ and $\bar v^b$ of $v^b$ which are sums of formal disks of $(+)$-action below $\alpha$. We count disks as follows. Each formal admissible disk in $(X^{ba},L^{ba})$ with factors in $\bar v^a$ and $\bar v^b$ has a tree structure: vertices correspond to disk factors and edges to glued Reeb chords. Label the vertices in the tree of a formal disk in $(X^{ba},L^{ba})$ by the corresponding formal disk from $\bar v^a$ or $\bar v^b$ and weight the formal admissible disk in $(X^{ba},L^{ba})$ by the product over the vertices in the tree of the polynomials which are the weights at its vertices. Summing the weighted disks corresponding to all such trees with vertices from $\bar v^a$ and $\bar v^b$ (there are only finitely many) in $\V(X^{ba},L^{ba};\tau^b,\tau^a)$ gives an element $\xi\in\V(X^{ba},L^{ba};\tau^b,\tau^a)$. Let $\bar p_{g^{ba}}(\tau^b,\tau^a)$ denote the polynomial coefficient of the formal disk $g^{ba}$ in the (finite) expansion of $\xi$. We define the polynomial coefficient $p_{g^{ba}}(\tau^b,\tau^a)$ of the formal disk $g^{ba}$ with $\A^+(g^{ba})<\alpha$ in the (possibly infinite) expansion of $\left(v^b\I v^a\right)$ as $p_{g^{ba}}(\tau^b,\tau^a)=\bar p_{g^{ba}}(\tau^b,\tau^a)$.

We next consider compositions of gluing pairings. Let $(X^c, L^c)$ be an exact cobordism with a $(+\infty)$-boundary $(Y^0,\Lambda^0)$, let $(X^b, L^b)$ an exact cobordism with a $(-\infty)$-boundary $(Y^0,\Lambda^0)$ and a $(+\infty)$-boundary $(Y^1,\Lambda^1)$, and $(X^a,L^a)$ an exact cobordism with a $(-\infty)$-boundary $(Y^1,\Lambda^1)$. Let $v^\ast\in\V(X^\ast,L^\ast;\tau^\ast)$, $\ast\in\{c,b,a\}$. Let $(X^{cba},L^{cba})$ denote the exact cobordism obtained by joining these three cobordisms.
\begin{lma}\label{l:assoc}
The following equation holds
\begin{equation*}
\Bigl(v^c\I\left(v^b\I v^a\right)\Bigr)
=\Bigl(\left(v^c\I v^b\right)\I v^a\Bigr)\in\V(X^{cba},L^{cba};\tau),
\end{equation*}
where $\tau=(\tau^a,\tau^b,\tau^c)$.
\end{lma}
\begin{pf}
Consider any formal disk in $(X^{cba},L^{cba})$ which arises after gluing according to the prescription in the left hand side. By subdividing it differently we see that it also arises after gluing as prescribed in the right hand side. Moreover, the weights are in both cases the product of the weights of all the factors. The same argument shows that any formal disk contributing to the right hand side also contributes to the left hand side. The lemma follows.
\end{pf}

We will often use linearized versions of the gluing pairing \eqref{e:glupair} defined as follows. Let $(X^b,L^b)$ and $(X^a,L^a)$ be as above. Fix $f^b\in\V(X^b,L^b)$ and consider for $v^a\in\V(X^a,L^a)$, $v^a\mapsto \left(f^b\I v^a\right)$ as a function $\V(X^a,L^a)\to\V(X^{ba},L^{ba})$. We define the {\em linearization} of this function at $f^a\in\V(X^a,L^a)$
\begin{equation}\label{e:linglupair1}
\left[\pa_{f^a}\!\left(f^b\I\right.\right]\colon\V(X^a,L^a)\to\V(X^{ba},L^{ba}),
\end{equation}
through the following equation
$$
\left(f^b\I f^a+\epsilon w^a\right)=\left(f^a\I f^b\right)
+\epsilon\left[\pa_{f^a}\!\left(f^b\I\right.\right](w^a)+\Ordo(\epsilon^2)
\in\V(X^{ba},L^{ba};\epsilon),
$$
for $w^a\in\V(X^a,L^a)$.

In a similar way, fixing $f^a\in\V(X^a,L^a)$ we consider $\left(v^b\I f^a\right)$, as a function $\V(X^b,L^b)\to\V(X^{ba},L^{ba})$. Its {\em linearization} at $f^b\in\V(X^b,L^b)$
\begin{equation}\label{e:linglupair2}
\left[\pa_{f^b}\!\left.\I f^a\right)\right]\colon\V(X^b,L^b)\to\V(X^{ba},L^{ba}),
\end{equation}
is defined through
$$
\left(f^b+\epsilon w^b\I f^a\right)=\left(f^a\I f^b\right)
+\epsilon\left[\pa_{f^b}\!\left.\I f^a\right)\right](w^b)+\Ordo(\epsilon^2)
\in\V(X^{ba},L^{ba};\epsilon),
$$
for $w^b\in\V(X^b,L^b)$.

\begin{lma}\label{l:compcoblin}
Let $(X^\ast,L^\ast)$, $\ast\in\{c,b,a\}$, be exact cobordisms as in Lemma \ref{l:assoc} and let $f^\ast\in\V(X^\ast,\Lambda^\ast)$. Then
\begin{align*}
&\left[\pa_{f^c}\left\I\left(f^b\I f^a\right)\Bigr)\right.\right](w^c)
=\left[\pa_{\left(f^c\I f^b\right)}\!\left.\I f^a\right)\right]
\Bigl(\left[\pa_{f^c}\!\!\left.\I f^b\right)\right](w^c)\Bigr),\quad w^c\in\V(X^c,\Lambda^c), \text{ and}\\
&\left[\pa_{f^a}\!\left.\Bigl(\left(f^c\I f^b\right)\right\I\right](w^a)
=\left[\pa_{\left(f^b\I f^a\right)}\!\left(f^c\I\right.\right]
\Bigl(\left[\pa_{f^a}\!\left(f^b\I\right.\right](w^a)\Bigr),\quad w^a\in\V(X^a,\Lambda^a)
\end{align*}
\end{lma}

\begin{pf}
Let $w^c\in\V(X^c,L^c)$. Lemma \ref{l:assoc} implies that
$$
\Bigl(\left(f^c+\epsilon w^c\I f^b\right)\left\I f^a\Bigr)\right.=
\left.\Bigl(f^c+\epsilon w^c\right\I \left(f^b\I f^a\right)\Bigr).
$$
The left hand side of this equation can be rewritten as
\begin{align*}
&\Bigl(\left(f^c\I f^b\right)
+\epsilon\left[\pa_{f^c}\left.\I f^b\right)\right](w^c)+\Ordo(\epsilon^2)
\left\I f^a\Bigr)\right.=\\
&\Bigl(\left(f^c\I f^b\right)\left\I f^a\Bigr)\right. +
\epsilon\left[\pa_{\left(f^c|f^b\right)}\!\left\I f^a\Bigr)\right.\right]
\Bigl(\left[\pa_{f^c}\!\!\left.\I f^b\right)\right](w^c)\Bigr)+\Ordo(\epsilon^2)
\end{align*}
and the right hand side as
$$
\Bigl(f^c\left\I\left(f^b\I f^a\right)\Bigr)\right.
+\epsilon\left[\pa_{f^c}\!\left\I \left( f^b\I f^a\right)\Bigr)\right.\right](w^c) +\Ordo(\epsilon).
$$
The first equation follows. The second equation is proved in a similar way.
\end{pf}
\section{Grading, filtration, and differential}\label{S:grafildiff}
In this section we introduce grading and filtration on the vector space of an exact cobordism. We also define a grading preserving differential on this space using the gluing pairing in combination with properties of moduli spaces of $J$-holomorphic disks. The grading preserving differential determines a spectral sequence and we show that if two exact cobordisms are joined at a common end then there are induced filtration preserving chain maps and hence morphisms of spectral sequences. We refer to Appendix \ref{A:B} for a more detailed discussion of moduli spaces of $J$-holomorphic disks in exact cobordisms. We assume throughout this section, that every exact cobordism comes equipped with a sufficiently generic almost complex structure $J$ which is adjusted to its symplectic form. For simplicity and to conform with the notation in more general situations (see Subsection \ref{s:mixpureini}, last paragraph) we will leave out $J$ from the notation, writing ``holomorphic disk'' instead of ``$J$-holomorphic disk''. The section is organized as follows. In Subsection \ref{s:grad}, the grading is introduced, and in Subsection \ref{s:filt}, the filtration. In Subsections \ref{s:pot} and \ref{s:ham}, distinguished vectors in the vector space of an exact cobordism and in the vector space of a symplectization are defined. (They are called the {\em potential-} and the {\em Hamiltonian} vector, respectively, to conform with notation in \cite{EGH}.) In Subsection \ref{s:diff}, the differential is introduced and in Subsection \ref{s:spsq}, the spectral sequence. Finally, in Subsection \ref{s:cob->chmap}, chain maps induced by joining cobordisms are discussed.

\subsection{Grading}\label{s:grad}
Let $(X,L)$ be an exact cobordism and let $g$ be a formal disk in $(X,L)$. The moduli space $\M(g)$ of holomorphic disks in $(X,L)$ which determines the formal disk $g$ (i.e., the formal disk map associated to the holomorphic disks in $\M(g)$ belong to the homotopy class $g$, see Remark \ref{r:hol->formal}) has  formal dimension $\dim(\M(g))$, which we explain how to compute below.
\begin{dfn}
The $\Z$-grading on the vector space $\V(X,L)$ is induced by the following grading of formal disks: the grading of a formal disk $g\in\V(X,L)$ is
$$
|g|=\dim(\M(g)).
$$
We say that an element $v\in\V(X,L)$ is {\em homogeneous} if it is a formal sum of admissible formal disks which all have the same degree.
\end{dfn}

\begin{rmk}\label{r:invlimgrad}
The grading on $\V(X,L)$ induces gradings of the finite $(+)$-action subspaces $\V_{[\alpha]}(X,L)\subset\V(X,L)$, $\alpha>0$.
\end{rmk}

If $g$ is a formal disk then $\dim(\M(g))$ equals the Fredholm index of the the linearized $\bar\pa_J$-operator at a map $u\colon D_m\to X$ which is asymptotic to Reeb chord strips of the Reeb chords $c_1,\dots,c_m$ at the $m$ boundary punctures of the punctured disk $D_m$, and which induces a formal disk map representing $g$. The source space of this linearized operator splits into a sum of the tangent space of the space of conformal structures on $D_m$, which is $(m-3)$-dimensional, and an infinite dimensional functional analytic space $U$ of vector fields along $u$ which are tangent to $L$ along the boundary. The Fredholm index of the linearized $\bar\pa_J$-operator restricted to $U$ can be computed by relating it to the Riemann-Hilbert problem, see e.g. \cite{EES2, EES3}, as follows. Pick a symplectic trivialization of the tangent bundle of $X$ along $u$ such that the linearized Reeb flow in this trivialization equals the identity along all Reeb chords. The tangent planes of $L$ along $u(\pa D_m)$ give $m$ paths $\gamma_1,\dots,\gamma_m$ of a Lagrangian subspaces of $\C^n$ (where $\C^n$ is determined by the trivialization) where $\gamma_j$ connects the tangent space $T_{e_j^1}L$ of $L$ at the Reeb chord endpoint $e_j^1$ where the formal disk leaves the Reeb chord $c_j$ to the tangent space $T_{e_{j+1}^0}L$ at the Reeb chord endpoint $e_{j+1}^0$ where the formal disk enters the Reeb chord $c_{j+1}$. (Here we use the convention $m+1=1$). The tangent space of $L$ at a Reeb chord endpoint $e^\sigma_j$, $\sigma=0,1$, splits as $T_{e^\sigma_j}\Lambda^\pm\times\R$, where $(Y^\pm,\Lambda^\pm)$ is the $(\pm\infty)$-boundary of $(X,L)$. Our genericity assumptions on the contact forms $\lambda^\pm$ on $Y^\pm$ implies that the linearized Reeb flow $\rho$ along $c_j$ takes the tangent space at the Reeb chord endpoint $e^\sigma_j$, where the Reeb field $R_{\lambda^\pm}$ points into $c_j$, $T_{e^\sigma_j}\Lambda^\pm$ to a Lagrangian subspace $\rho(T_{e^\sigma_j}\Lambda^\pm) \subset\bigl(\krn(\lambda^\pm),d\lambda^\pm\bigr)\subset T_{y^\tau_j} Y$ which is transverse to the tangent space $T_{e^\tau_j}\Lambda^\pm$ at the other endpoint $e^\tau_j$ of $c_j$ ($\sigma\ne \tau\in\{0,1\}$). We close the paths $\gamma_1,\dots,\gamma_m$ to a loop using paths $\hat\gamma_1,\dots,\hat\gamma_m$, where $\hat\gamma_{j-1}$ connects the endpoint $T_{e^0_j}\Lambda^\pm\times\R$ of $\gamma_{j-1}$ to the start-point $T_{e^1_j}\Lambda^\pm\times\R$ of $\gamma_j$, as follows. The path $\hat \gamma_{j-1}$ leaves the $\R$-factor fixed. If $e_j^0$ is the endpoint of a Reeb chord in the $(+\infty)$-boundary, then rotate $\rho(T_{e^0_j}\Lambda^+)$ in the negative direction to $T_{e_{j}^1}\Lambda^+$ in $\krn(\lambda^+)\subset T_{e_j^1}Y^+$. If $e_j^0$ is the endpoint of a Reeb chord in the $(-\infty)$-boundary, then rotate $\rho^{-1}(T_{e^0_j}\Lambda^-)$ in the negative direction to $T_{e_{j}^1}\Lambda^-$ in $\krn(\lambda^-)\subset T_{e_j^1}Y^-$. (Here we use the following notation: if $W$ is a symplectic vector space then a negative rotation of a Lagrangian subspace $V_0\subset W$ to a Lagrangian subspace $V_1\subset W$ transverse to $V_0$, is a path of the form $e^{-s\frac{\pi}{2}I}V_0$, $0\le s\le 1$, where $I$ is a complex structure compatible with the symplectic structure and such that $IV_0=V_1$.) The concatenation $\gamma=\gamma_1\ast\hat \gamma_1\ast\cdots\gamma_m\ast\hat\gamma_m$ is a a loop of Lagrangian subspaces of $\C^n$. The index of the restriction of the linearized $\bar\pa_J$-operator to $U$ is $n+\mu(\gamma)$ and the formal dimension is consequently
$$
|g|=\dim(\M(g))=n-3+\mu(\gamma)+m,
$$
where $2n=\dim(X)$, where $\mu$ denotes Maslov index, and where $m$ is the number of punctures, see \cite{EES1, EES4}.

The dimension is additive in the following sense. With notation as in Section \ref{s:vsp+op}, if $v^a\in\V(X^a,\Lambda^a)$ and $v^b\in\V(X^b,\Lambda^b)$ then the dimension of a formal disk contributing to $\left(v^b\I v^a\right)$ equals the sum of the dimensions of its factors. To see this one uses the following fact: at each Reeb chord where disks are glued, $(n-1)$ negative half-turns are lost (since the closing up rotations at the punctures were erased by the gluing operation) and two punctures are erased.

\subsection{Filtration}\label{s:filt}
Let $(X,L)$ be an exact cobordism and assume that $L=L_1\cup\dots\cup L_k$ is subdivided into $k$ pieces.
\begin{lma}
No admissible formal disk in $(X,L)$ has more than $k$ mixed punctures and consequently no more than $k$ positive punctures.
\end{lma}

\begin{pf}
Assume that $D$ is the source of a formal disk with more than $k$ mixed punctures. Then there exists a pair of boundary components of $D$ which map to the same piece $L_j$ of $L$ and such that both complementary arcs in $\pa D$ of these two boundary components must contain boundary components mapping to some piece of $L$ other than $L_j$. An arc connecting the two boundary components in the pair then contradicts $({\bf a2})$ of Definition \ref{d:adm}. Lemma \ref{l:admpurepos} implies that any admissible formal disk either has one pure positive puncture or all its positive punctures are mixed. The lemma follows.
\end{pf}
If $g$ is a formal disk then let $\pos(g)$ denote the number of positive punctures of $g$.
\begin{dfn}
For $1\le p\le k$, let $F^p\V(X,L)$ denote the subspace of $\V(X,L)$ which contains all elements $v$ of the form
$$
v=\sum_{g\in G,\, \pos(g)\ge p}\delta_g\,g,\quad\quad\delta_g\in\Z_2.
$$
 The filtration of $\V(X,L)$ is
$$
0\subset F^k\V(X,L)\subset F^{k-1}\V(X,L)\subset \dots\subset F^1\V(X,L)=\V(X,L).
$$
\end{dfn}

\begin{rmk}\label{r:invlimfilt}
The number of positive punctures induces filtrations also on the finite $(+)$-action parts of $\V(X,L)$. Using the notation in Remark \ref{r:invlim} we have
$$
0\subset F^k\V_{[\alpha]}(X,L)\subset F^{k-1}\V_{[\alpha]}(X,L)\subset \dots\subset F^1\V_{[\alpha]}(X,L)=\V_{[\alpha]}(X,L),
$$
for all $\alpha> 0$.
\end{rmk}

\subsection{The potential vector}\label{s:pot}
Let $(X,L)$ be an exact cobordism. The moduli space of $0$-dimensional admissible holomorphic disks in $(X,L)$ with action (or which is the same, $\omega$-energy) bounded above is a $0$-dimensional compact manifold, see Lemmas \ref{l:classtv} and \ref{l:nicecomp}. In other words, it is a finite collection of points and to each point is associated a formal disk $g$. In particular, the number of $0$-dimensional admissible holomorphic disks in $(X,L)$ with bounded $(+)$-action is finite. Furthermore, if a formal disk $g$ has a holomorphic representative then $\A(g)\ge 0$, see Lemma \ref{l:action}. Define the {\em potential vector} $f\in\V^+(X,L)$ as
$$
f=\sum_{\dim(\M(g))=0}|\M(g)|\,g,
$$
where $|\M(g)|$ denotes the modulo $2$ number of points in the $0$-dimensional moduli space $\M(g)$. Note that $f\in\V^+(X,L)$ is a homogeneous element of grading $0$.

In the special case when the cobordism is trivial, because of the translational invariance in symplectizations, for generic adjusted almost complex structure, the only holomorphic disks of formal dimension $0$ are Reeb chord strips, see Appendix \ref{A:B}, with one positive and one negative puncture (at the same Reeb chord). Thus for a trivial concordance the potential vector $f$ is extremely simple: if $\Lambda\subset Y$ is a Legendrian submanifold then
$$
f=\sum_{c\in\chord} g_c\in\V(Y\times\R,\Lambda\times\R),
$$
where $g_c$ is the formal disk represented by the Reeb chord strip of the Reeb chord $c$.

\subsection{The Hamiltonian vector of a symplectization}\label{s:ham}
Consider the trivial cobordism $(Y\times\R,\Lambda\times\R)$ associated to a  Legendrian submanifold $\Lambda\subset Y$. Using the fact, pointed out above, that $\R$ acts on the moduli spaces of holomorphic disks and that for generic adjusted almost complex structure, the only $0$-dimensional holomorphic disks are Reeb chord strips, we conclude that that the moduli space $\M$ of holomorphic disks of action bounded above of formal dimension $1$, when divided out by this $\R$-action, forms a compact $0$-dimensional manifold. Write ${\widehat\M}=\M/\R$ and call $\widehat\M$ the {\em reduced} moduli space. Define the {\em Hamiltonian vector} $h\in\V(Y\times\R,\Lambda\times\R)$ of a trivial cobordism as
$$
h=\sum_{\dim({\widehat\M}(g))=0}|{\widehat\M}(g)|\,g,
$$
where the sum ranges over admissible formal disks $g$. Note that $h\in\V(Y\times\R,\Lambda\times\R)$ is a homogeneous element of grading $1$ and that $\A(g)>0$ for any $g$ for which ${\widehat\M}(g)$ is non-empty and hence $h\in\V^+(Y\times\R,\Lambda\times\R)$.

\subsection{The differential}\label{s:diff}
Let $(X,L)$ be an exact cobordism with ends $(Y^\pm\times\R_{\pm},\Lambda^\pm\times\R_{\pm})$. Gluing the trivial cobordisms $(Y^\pm\times\R,\Lambda^\pm\times\R)$ to $(X,L)$ does not change the cobordism. In particular, if $v^\pm\in\V(Y^\pm\times\R,\Lambda^\pm\times\R;\tau)$ and
$w\in\V(X,L;\tau)$ then
$\left(v^-\I w\right)\in\V(X,L;\tau)$ and
$\left(w\I v^+\right)\in\V(X,L;\tau)$, see Subsection \ref{s:glupair} for notation.

Let $f\in\V^+(X,L)$ and $f^\pm\in\V^+(Y^\pm\times\R,\Lambda^\pm\times\R)$ be the potential vectors and let $h^\pm\in\V^+(Y^\pm\times\R,\Lambda^\pm\times\R)$ be the Hamiltonian vectors. We associate operators $h^{\pm}\colon\V(L,X;\epsilon)\to\V(L,X;\epsilon)$ to the Hamiltonian vectors (for simplicity, we denote these operators by the same symbols as the vectors themselves) through the following equations
\begin{align*}
\left(v\I f^+ + \tau h^+\right)= v + \tau\, h^+(v) +\Ordo(\tau^2)\in\V(X,L;\tau,\epsilon),\\
\left(f^-+\tau h^-\I v\right)= v + \tau\, h^-(v) +\Ordo(\tau^2)\in\V(X,L;\tau,\epsilon),
\end{align*}
where $v\in\V(X,L;\epsilon)\subset\V(X;L;\tau,\epsilon)$ and  $f^\pm+\tau h^\pm\in\V^+(Y^\pm\times\R,\Lambda^\pm\times\R;\tau)\subset
\V(Y^\pm\times\R,\Lambda^\pm\times\R;\tau,\epsilon)$.
Define
$$
h=h^++h^-\colon
\V(X,L;\epsilon)\to\V(X,L;\epsilon).
$$

\begin{lma}\label{l:hpot}
If $f\in\V^+(X,L)$ is the potential vector then $h(f)=0$.
\end{lma}

\begin{pf}
To see this we note that any admissible formal disk contributing to $h(f)$ can be viewed as a broken disk which is an endpoint of a $1$-dimensional moduli space of holomorphic disks in $(X,L)$. By Lemma \ref{l:nicecomp}, and uniqueness of gluing, the other endpoint of this $1$-dimensional moduli space correspond to some other broken admissible disk with one factor from $h^+$ or $h^-$ and remaining factors from $f$. Hence also the other endpoint contributes to $h(f)$. The lemma follows.
\end{pf}

We next consider the linearization of the Hamiltonian operator $h$ at the potential vector $f$. This linearization is a map $d^f\colon\V(X,L)\to\V(X,L)$ defined through the following equation for $v\in\V(X,L)$,
$$
h(f+\epsilon v)= \epsilon\,\,d^f(v) +\Ordo(\epsilon^2)\in\V(X,L;\epsilon).
$$

\begin{lma}\label{l:diffconc}
The map $d^f\colon\V(X,L)\to\V(X,L)$ is a filtration preserving differential of degree $1$. That is, $d^f\circ d^f = 0$, $d^f(F^p\V(X,L))\subset F^p\V(X,L)$, and if $v\in\V(X,L)$ is homogeneous of degree $k$ then $d^f(v)$ is homogeneous of degree $k+1$. Furthermore, $d^f$ is $(+)$-action non-decreasing: if $v\in\V(X,L)$ is a formal disk then the $(+)$-action $\A^+(w)$ of any formal disk $w$ contributing to $d^f(v)$ satisfies $\A^+(w)\ge\A^+(v)$.
\end{lma}

\begin{pf}
Let $v\in\V(X,L)$ and let $h^+$ and $h^-$ denote the Hamiltonian vectors at the positive and negative ends of $(X,L)$, respectively. We first show $d^f\circ d^f = 0$. By linearity of $d^f$ and since for each fixed formal disk $g$ of $(X,L)$ there is only a finite number of formal disks in $v$, in $f$, and in $h^\pm$ which can contribute to the coefficient of $g$ in $d^f\circ d^f(v)$, see Lemmas \ref{l:factor1} and \ref{l:factor2}, it suffices to show that $d^f\circ d^f(v)=0$ in the case when $v$ is a single formal disk. In other words, we must show
\begin{equation}\label{e:d^2=0}
h\bigl(f+h(f+\epsilon v)\bigr)=\Ordo(\epsilon^2),
\end{equation}
for formal disks $v$.

The constant term of the left hand side in \eqref{e:d^2=0} vanishes by Lemma \ref{l:hpot}. Consider the linear term. Any formal disk contributing to the linear term in \eqref{e:d^2=0} has one $v$-factor, two $h^\pm$-factors, and all other factors $f$-factors. We distinguish the $h^\pm$-factor corresponding to the {\em outer} Hamiltonian operator $h$ and call it {\em charged}. Note that the non-charged $h^\pm$-factor is connected to $v$. Below we will make use of slightly more general disks. We therefore define an $\lt$-disk ($\lt$ for ``linear term'') as a formal disk with one $v$-factor, two $h^\pm$-factors one which is connected to $v$ and the other one distinguished as charged, and remaining factors $f$-factors.

Consider the positive- or negative end (i.e., all positive- or negative punctures) of the charged $h^\pm$-factor of an $\lt$-disk. We say that this end is {\em un-obstructed} if either there is some $h^\pm$-factor attached at the end, or if there are only $f$-factors attached at the end. An end which is not un-obstructed is called {\em obstructed}. We say that an $\lt$-disk is {\em isolated}, {\em boundary}, or {\em interior} if its charged $h^\pm$-factor has two, one, or zero obstructed ends, respectively.

We first show that the contribution to \eqref{e:d^2=0} from isolated $\lt$-disks vanishes. Both $h^\pm$-factors of an isolated $\lt$-disk are connected to $v$. Changing the factorization by moving the charge from one $h^\pm$-factor to the other we see that the formal disk appears in two ways as an $\lt$-disk in \eqref{e:d^2=0}. Thus the contribution from isolated $\lt$-disks vanishes.

Consider a non-isolated $\lt$-disk with an un-obstructed end of its charged $h^\pm$-factor distinguished and called {\em active}. We define the following propagation law for such an $\lt$-disk.
\begin{itemize}
\item If there are only $f$-factors attached at the active end of the charged $h^\pm$-factor, then view the broken disk consisting of the charged $h^\pm$-factor and all these $f$-factors as the boundary of a $1$-dimensional moduli space of holomorphic disks in $(X,L)$. Moving to the other end of this moduli space we get another broken disk with $f$-factors and one $h^\pm$-factor. Charge the new $h^\pm$-factor and activate the end of it where the new $f$-factors are not attached.
\item If there is an $h^\pm$-factor attached at the active end of the charged $h^\pm$-factor, then one of these two $h^\pm$-factors is attached to $v$. View the broken disk consisting of the two $h^\pm$-factors as the boundary of a $1$-dimensional reduced moduli space in the symplectization of the $(\pm\infty)$-boundary. Moving to the other end of the reduced moduli space we get another broken disk with two $h^\pm$-factors. Charge the one of them which is not attached to $v$ and activate the end of it where the other $h^\pm$-factor is not attached.
\end{itemize}
This propagation law associates to an $\lt$-disk with an un-obstructed end of its charged $h^\pm$-factor activated a unique non-isolated $\lt$-disk with an end of its charged $h^\pm$-factor activated. In particular, using this propagation law repeatedly starting at a boundary $\lt$-disk the process stops at some other uniquely determined boundary $\lt$-disk, see Figure \ref{fig:comb1mfd}.
\begin{figure}[htbp]
\begin{center}
\includegraphics[angle=0, width=12cm]{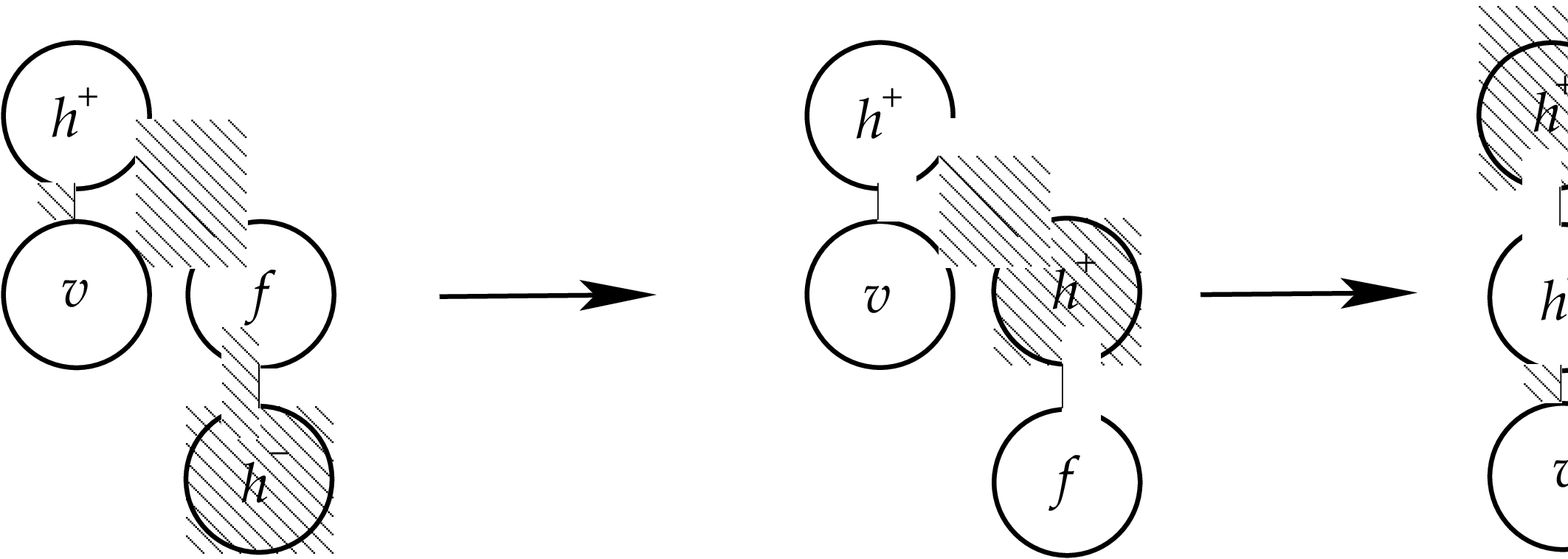}
\end{center}
\caption{The propagation law, the charged $h^\pm$-factor is shaded.}
\label{fig:comb1mfd}
\end{figure}
Noting that every non-isolated $\lt$-disk which contributes to \eqref{e:d^2=0} is a boundary $\lt$-disk and that any boundary $\lt$-disk contributes to \eqref{e:d^2=0} it follows that also the contribution from non-isolated $\lt$-disks to \eqref{e:d^2=0} vanishes. We conclude that $d^f\circ d^f=0$.

The fact that $d^f$ increases grading by $1$ follows from Subsection \ref{s:grad} since $f$ is homogenous of degree $0$ and $h^\pm$ are homogeneous of degree $1$. To see that $d^f(F^p\V(X,L))\subset F^p(\V(X,L))$, note that gluing operations of admissible disks never decreases the number of positive punctures: since each admissible disk has at least one positive puncture, a positive puncture erased by gluing is compensated by at least one new positive puncture.

To see that $d^f$ is $(+)$-action non-decreasing we argue as follows. Formal disks contributing to $d^f(v)$ have one $v$-factor and remaining factors from $h^\pm$ and from $f$. Formal disks in $h^\pm$ and in $f$ have holomorphic representatives and hence have non-negative action. If $g$ is a partial formal disk and some factor is added to it from below giving a new formal disk $g'$ then obviously $\A^+(g')\ge\A^+(g)$. If $g$ is a partial formal disk and a disk $h$ is added to it from above giving a new disk $g'$ then $\A^+(g')\ge \A^+(g)+\A(h)$ so that $\A^+(g')\ge \A^+(g)$ if $\A(h)\ge 0$.  It follows that $d^f$ does not decrease $(+)$-action.
\end{pf}

\begin{rmk}
In the case that $(X,L)=(Y\times\R,\Lambda\times\R)$ is a trivial cobordism, the argument in the proof of Lemma \ref{l:diffconc} shows that the first order terms $d^+$ and $d^-$ in $h^+(f+\epsilon g)$ and $h^-(f+\epsilon g)$, respectively are both differentials, separately. Moreover, it is immediate that $d^+\circ d^-=d^-\circ d^+$. (In this case the {\em vectors} $h^+=h$ and $h^-=h$ are identical. However, the {\em operators} $h^+$ and $h^-$ are not: one corresponds to attaching $h$ from above, the other to attaching $h$ from below.)
\end{rmk}

\begin{rmk}
In this language, the contact homology differential in Legendrian contact homology, see \cite{Ch,EGH,EES1,EES3}, can be described as ${\tilde h}$ acting on the algebra generated by all negative Reeb chords, where $\tilde h$ is the Hamiltonian operator of the trivial cobordism associated to a Legendrian submanifold $\Lambda$ in a contact manifold $Y$  corresponding to formal disks with only one positive puncture.
\end{rmk}

\begin{rmk}\label{r:invlimdiff}
The above construction also induces differentials on the finite action parts of $\V(X,L)$. With notation as in Remark \ref{r:invlim}, for $\alpha>0$, define
$$
d^f_\alpha\colon \V_{[\alpha]}(X,L)\to\V_{[\alpha]}(X,L),\quad d^f_\alpha(v)=\pi_{[\alpha]}\left(d^f(v)\right),
$$
where $\pi_{[\alpha]}\colon \V(X,L)\to\V_{[\alpha]}(X,L)$ is the map which truncates a formal sum in $\V(X,L)$ at action $\alpha$. More precisely, if $w\in\V(X,L)$ then $w$ can be written uniquely as $w=w_{[\alpha]}+w_{[\alpha+]}$ where $w_{[\alpha]}$ is a linear combination of formal disks of $(+)$-action $\A^+\le \alpha$ and where $w_{[\alpha+]}$ is a formal sum of disks of $(+)$-action $\A^+>\alpha$, and $\pi_{[\alpha]}(w)=w_{[\alpha]}$. The fact that $d^f_\alpha$ is a differential (i.e., that $d^f_\alpha\circ d^f_\alpha=0$) is a straightforward consequence of the fact that $d^f$ does not decrease $(+)$-action. This fact also implies that, for $\alpha>\beta>0$, the projections $\pi^\alpha_\beta\colon\V_{[\alpha]}(X,L)\to\V_{[\beta]}(X,L)$, see Remark \ref{r:invlim}, are filtration preserving chain maps. That is,
\begin{equation}\label{e:pichmap}
\pi^\alpha_\beta\circ d^f_\alpha = d^f_\beta\circ\pi^\alpha_\beta, \qquad
\pi^\alpha_\beta\left(F^p\V_{[\alpha]}(X,L)\right)=F^p\V_{[\beta]}(X,L),\qquad
\text{for }\alpha>\beta.
\end{equation}
\end{rmk}

\subsection{The rational admissible SFT spectral sequence}\label{s:spsq}
Let $(X,L)$ be an exact cobordism where $L=L_1\cup\dots\cup L_k$ is subdivided into $k$ pieces and let $\alpha>0$.
\begin{dfn}\label{d:spsqcob}
The {\em $(+)$-action $\alpha$ rational admissible SFT spectral sequence} of the exact cobordism $(X,L)$ is the cohomological spectral sequence
$$
\bigl\{E^{p,q}_{r\,;[\alpha]}(X,L)\bigr\}_{r=1}^k,
$$
induced by the filtration respecting differential $d^f_\alpha\colon\V_{[\alpha]}(X,L)\to\V_{[\alpha]}(X,L)$, where $f$ is the potential vector of $(X,L)$, which has $E_1$-term
$$
E_{1\,;[\alpha]}^{p,q}(X,L)=
H^{p+q}\left(F^p\V_{[\alpha]}(X,L)/F^{p+1}\V_{[\alpha]}(X,L)\right),
$$
(see Remarks \ref{r:invlim}, \ref{r:invlimfilt}, and \ref{r:invlimdiff}).
\end{dfn}

For $\alpha>\beta>0$, \eqref{e:pichmap} implies that the projection map $\pi^\alpha_\beta$ induces a morphism of spectral sequences,
$$
\begin{CD}
\bigl\{E^{p,q}_{r\,;[\beta]}(X,L)\bigr\}
@<{{\pi^\alpha_\beta}_\ast}<<
\bigl\{E^{p,q}_{r\,;[\alpha]}(X,L)\bigr\}
\end{CD}.
$$
The collection $\left\{\bigl\{E^{p,q}_{r\,;[\alpha]}(X,L)\bigr\},{\pi^\alpha_\beta}_\ast\right\}$ is thus an inverse system of spectral sequences. Since the inverse limit functor is not exact, the graded vector space defined by the inverse limit
\begin{equation}\label{e:dfninvlim}
\bigl\{E^{p,q}_{r}(X,L)\bigr\}_{r=1}^k=
\lim_{\longleftarrow\,\alpha} \bigl\{E^{p,q}_{r\,;[\alpha]}(X,L)\bigr\}_{r=1}^k,
\end{equation}
is in general not a spectral sequence. However, if there exists $\alpha_0$ such that for any $\alpha>\beta>\alpha_0$, ${\pi^\alpha_\beta}_\ast$ is an isomorphism of spectral sequences then
\begin{equation}\label{e:bddiso}
\bigl\{E^{p,q}_{r}(X,L)\bigr\}_{r=1}^k=
\bigl\{E^{p,q}_{r\,;[\alpha]}(X,L)\bigr\}_{r=1}^k,
\end{equation}
for any $\alpha>\alpha_0$ and the limit has the structure of a spectral sequence.

If the action of the Reeb chords at the $(+\infty)$-boundary $(Y^+,\Lambda^+)$ of $(X,L)$ is uniformly bounded by $\alpha_0$ then \eqref{e:bddiso} holds. We note that if there is such a bound, then our genericity assumptions for the contact form $\lambda^+$ on $Y^+$ implies that $\Lambda^+$ has only a finite number of Reeb chords as follows. If there are infinitely many then some subsequence of Reeb chord endpoints converges. Since the action is uniformly bounded, it is straightforward to show that some subsequence of the corresponding Reeb chords converges to a Reeb chord. Thus the Reeb chords are non-isolated contrary to our genericity assumption that the linearized Reeb flow takes the tangent space of $\Lambda^+$ at the start-point of any Reeb chord to a Lagrangian subspace of the contact hyperplane at the end-point of the Reeb chord which is transverse to the tangent space of $\Lambda^+$ at this end-point.

Let $\Lambda=\Lambda_1\cup\dots\cup\Lambda_k\subset Y$ be a Legendrian submanifold subdivided into pieces.
\begin{dfn}\label{d:spsqLeg}
The {\em rational admissible SFT invariant} of $(Y,\Lambda)$ is the inverse limit
$$
\bigl\{E^{p,q}_r(Y,\Lambda)\bigr\}_{r=1}^k=
\bigl\{E^{p,q}_r(Y\times\R,\Lambda\times\R)\bigr\}_{r=1}^k,
$$
see \eqref{e:dfninvlim}. If $\Lambda$ has only finitely many Reeb chords, then the rational SFT invariant is a spectral sequence. If the number of Reeb chords of $\Lambda$ is unbounded then the rational SFT invariant is a graded vector space.
\end{dfn}

\begin{rmk}\label{r:invlimspsq}
Let $(X,L)$ be an exact cobordism with infinitely many Reeb chords at its $(+\infty)$-boundary. The differential $d^f\colon\V(X,L)\to\V(X,L)$ respects the filtration and induces a spectral sequence
$$
\bigl\{{\widetilde E}^{p,q}_r(X,L)\bigr\}
$$
with $E_1$-term
$$
{\widetilde E}_{1}^{p,q}(X,L)
=H^{p+q}\left(F^p\V(X,L)/F^{p+1}\V(X,L)\right).
$$
To describe the relation between $\bigl\{{\widetilde E}_r^{p,q}(X,L)\bigr\}$ and the inverse limit  $\bigl\{E_r^{p,q}(X,L)\bigr\}$ of \eqref{e:dfninvlim}, one must take into account the non-exactness of the inverse limit functor measured by the functor $\lim^1$, see \cite{Mi}.
\end{rmk}

\subsection{Joining cobordisms and filtration preserving chain maps}\label{s:cob->chmap}
Let $(X^b,L^b)$ and $(X^a,L^a)$ be exact cobordisms and assume that $(Y^0,\Lambda^0)$ is a $(+\infty)$-boundary of $(X^b,L^b)$ and a $(-\infty)$-boundary of $(X^a,L^a)$. Consider the exact cobordism $(X^{ba},L^{ba})$ obtained by joining the ends corresponding to $(Y^0,\Lambda^0)$. Write $f^a$, $f^b$, and $f^{ba}$ for the potential vectors in $\V^+(X^a,L^a)$, $\V^+(X^b,L^b)$, and $\V^+(X^{ba},L^{ba})$, respectively.
\begin{lma}\label{l:potcob}
The potential vectors satisfy
$$
f^{ba}=\left(f^b\I f^a\right).
$$
\end{lma}
\begin{pf}
Stretching $(X^{ba},L^{ba})$ along $(Y^0,\Lambda^0)$ it follows from Lemma \ref{l:nicecomp}, that any admissible rigid holomorphic disk in $(X^{ba},L^{ba})$ breaks into admissible rigid disks in $(X^a,L^a)$ and $(X^b,L^b)$. Conversely, any such broken configuration can be glued uniquely to an admissible rigid disk in the joined cobordism. Thus the left- and the right hand sides count the same objects.
\end{pf}
Consider the linearizations
\begin{align*}
&\left[\pa_{f^a}\!\left(f^b\I\right.\right]\colon\V(X^a,L^a)\to\V(X^{ba},L^{ba})
\text{ and}\\
&\left[\pa_{f^b}\!\left.\I f^a\right)\right]\colon\V(X^b,L^b)\to\V(X^{ba},L^{ba}),
\end{align*}
introduced in \eqref{e:linglupair1} and \eqref{e:linglupair2}, respectively.
\begin{lma}\label{l:chmap}
The maps $\left[\pa_{f^a}\!\left(f^b\I\right.\right]$ and
$\left[\pa_{f^b}\!\left.\I f^a\right)\right]$ are $(+)$-action non-decreasing, filtration preserving chain maps,
\begin{align}\label{e:cobchmapabove}
\left[\pa_{f^a}\!\left(f^b\I\right.\right]\circ d^{f^a} &=
d^{f^{ba}}\circ \left[\pa_{f^a}\!\left(f^b\I\right.\right],\\\label{e:cobchampbelow}
\left[\pa_{f^b}\!\left.\I f^a\right)\right]\circ d^{f^b} &= d^{f^{ba}}\circ \left[\pa_{f^b}\!\left.\I f^a\right)\right].
\end{align}
\end{lma}

\begin{pf}
We prove \eqref{e:cobchmapabove}. The proof of \eqref{e:cobchampbelow} is similar.
Write $h$ and $H$ for the Hamiltonian operators on $\V(X^a,L^a;\epsilon)$ and $\V(X^{ba},L^{ba};\epsilon)$, respectively. It is enough to show that
\begin{equation}\label{e:cobchmap}
\left(f^b\I f^a + h(f^a+\epsilon v^a)\right) + \left(f^b\I f^a\right) + H\Bigl(\left(f^b\I f^a+\epsilon v^a\right)\Bigr)=\Ordo(\epsilon^2),
\end{equation}
for $v\in\V(X^a,L^a)$, since the coefficient of the linear term of \eqref{e:cobchmap} agrees with the sum of the right- and left hand sides of \eqref{e:cobchmapabove} evaluated on $v$. As in the proof of Lemma \ref{l:diffconc} we find that it suffices to demonstrate \eqref{e:cobchmap} in the case that $v$ is a single formal disk. Let $h^0$ denote the Hamiltonian vector of $(Y^0\times\R,\Lambda^0\times\R)$ and let $h^+$ and $h^-$ denote the Hamiltonian vectors of the positive and negative ends of $(X^{ba},L^{ba})$, respectively. A formal disk from $h^\pm$ or $h^0$ will be called an $h^\ast$-disk.

Any formal disk contributing to the linear term in \eqref{e:cobchmap} is a formal disk with one $v$-factor, one $h^\ast$-factor, and, using Lemma \ref{l:potcob}, with all other factors from the potential vectors $f^a$ and $f^b$. We say that a disk with these properties is an {\em $\lt$-disk}.

Consider an end of the $h^\ast$-factor of an $\lt$-disk. We say that such an end is {\em un-obstructed} if at all of its punctures $f^a$- or $f^b$-factors are attached. An end which is not un-obstructed is {\em obstructed}. An $\lt$-disk is {\em isolated}, {\em boundary}, or {\em interior} if its $h^\ast$-factor has two, one, or zero obstructed ends. The contribution to \eqref{e:cobchmap} from isolated $\lt$-disks vanishes: each  such disk arises once in the first and once in the third term.

Consider a non-isolated $\lt$-disk with an un-obstructed end of its $h^\ast$-factor distinguished and called {\em active}. We define the following propagation law for such an $\lt$-disk.
\begin{itemize}
\item The $h^\ast$-factor together with the $f^a$- or $f^b$-factors attached at its active end is a broken disk which is the boundary of a $1$-dimensional moduli space of holomorphic disks in $(X^a,L^a)$ or in $(X^b,L^b)$. Moving to the other end of the moduli space we find a similar broken disk with one $h^\ast$-factor and several $f^a$- or $f^b$-factors. Activate the end of the new $h^\ast$-factor where the new $f^a$- or $f^b$-factors are not attached.
\end{itemize}
Using this propagation law starting at a boundary $\lt$-disk the process stops at some other uniquely determined boundary $\lt$-disk. Since every non-isolated $\lt$-disk which contributes to \eqref{e:cobchmap} is a boundary $\lt$-disk and since any boundary $\lt$-disk contributes to \eqref{e:cobchmap}, we conclude that \eqref{e:cobchmapabove} holds.

The chain maps respect the filtration since the gluing operation does not decrease the number of positive punctures. They do not decrease $(+)$-action since all formal disks in the potential vectors have holomorphic representatives and hence non-negative action, see Lemma \ref{l:diffconc}.
\end{pf}

\begin{rmk}
Lemma \ref{l:compcoblin} implies that the chain maps corresponding to several successively joined exact cobordisms have natural associativity properties.
\end{rmk}
\section{Deformations and chain homotopies}\label{S:inv}
In this section we show that the rational admissible SFT spectral sequence is invariant under deformations of exact cobordisms and that morphisms induced by joining cobordisms have similar invariance properties. The section is organized a follows. In Subsection \ref{s:incrdisk}, an algebraic description of how the potential vector of an exact cobordism changes as the cobordism is deformed is presented. In Subsection \ref{s:chiso}, it is shown that this change leads to a chain isomorphism which induces isomorphisms between spectral sequences. In Subsection \ref{s:chhom}, it is demonstrated that if a cobordism is joined to a $1$-parameter family of cobordisms then there is an induced chain homotopy connecting the composition of the chain isomorphism of the joined $1$-parameter family and the chain map at the beginning of the original $1$-parameter family to the chain map at the end of the original family.

\subsection{Increment disks}\label{s:incrdisk}
Let $(X^b,L^b)$ and $(X^a,L^a)$ be exact cobordisms. Assume that $(Y^0,\Lambda^0)$ is a $(+\infty)$-boundary of $(X^b,L^b)$ and a $(-\infty)$-boundary of $(X^a,L^a)$. Let $(X^{ba},L^{ba})$ denote the cobordism obtained by joining $(X^b,L^b)$ and $(X^a,L^a)$ at $(Y^0,\Lambda^0)$. Let $\hat\Lambda^0$ be a piece of $\Lambda^0$, let $\hat L^a$ and $\hat L_b$ be corresponding pieces of $L^a$ and $L^b$, respectively, and let $\hat L^{ba}$ be the corresponding piece of $L^{ba}$. Let $h^a\in\V(X^a,L^a)$, let $k^b\in\V(X^b,L^b)$ be such that every formal disk in $k^b$ has some boundary component mapping to $\hat L^b$. We define the {\em $\bigl\{k^b\to h^a\bigr\}$ split gluing operation}
$$
\bigl\{k^b\to h^a\bigr\}\colon
\V(X^b,L^b;\epsilon)\times\V(X^b,L^b;\epsilon)\times\V(X^a,L^a;\epsilon)\to
\V(X^{ba},L^{ba};\epsilon)
$$
in the following way. The vector $\bigl\{k^b\to h^a\bigr\}\!\left(v_0^b,v_1^b,v^a\right)$ is the sum of all admissible disks constructed as follows. First pick a disk in $h^a$ and attach a disk in $k^b$ to it. If the resulting (partial) disk is admissible then the $\hat L^{ba}$-component of its boundary induces an ordering of its punctures, see Remark \ref{r:order}. Second attach $v_0^b$-factors at the negative punctures at Reeb chords in $(Y^0,\Lambda^0)$ of the partial formal disk which lie after the $k^b$-factor, and attach $v_1^b$-factors at such punctures which lie before the $k^b$-factor. Third attach $v^a$-factors at punctures at Reeb chords in $(Y^0,\Lambda^0)$. Continue like this, in each step attaching $v_0^b$ after the $k^b$-factor and $v_1^b$-factors before it. As for the gluing pairing in Subsection \ref{s:vsp+op}, we see that this is a finite process below any given $(+)$-action.

Similarly, if $h^b\in\V(X^b,L^b)$ and if $k^a\in\V(X^a,L^a)$ are such that every formal disk in $k^a$ has some boundary component mapping to $\hat L^a$, then we define the {\em $\bigl\{h^b\leftarrow k^a\bigr\}$ split gluing operation}
$$
\bigl\{h^b\leftarrow k^a\bigr\}\colon
\V(X^b,L^b;\epsilon)\times\V(X^a,L^a;\epsilon)\times\V(X^a,L^a;\epsilon)\to
\V(X^{ba},L^{ba};\epsilon)
$$
as follows. The vector $\bigl\{h^b\leftarrow k^a\bigr\}\!(v^b,v_0^a,v_1^a)$ is the sum of all admissible disks constructed as follows. Pick a disk in $h^b$ and attach a disk in $k^a$ to it. If the resulting (partial) disk is admissible then the $\hat L^{ba}$-component of its boundary induces an ordering of its punctures. Attach $v_0^a$-factors at the positive punctures at Reeb chords in $(Y^0,\Lambda^0)$ of the partial formal disks which lie after the $k^a$-factor, and attach $v_1^a$-factors at such punctures which lie before the $k^a$-factor. Then attach $v^b$-factors at punctures at Reeb chords in $(Y^0,\Lambda^0)$. Continue like this, in each step attaching $v_0^a$ after the $k^b$-factor and $v_1^a$-factors before it. Again this is a finite process below any given $(+)$-action.

Consider an exact cobordism $(X,L)$, let $\hat L$ be one of the pieces of $L$, let $u,v\in\V(X,L)$, and let $k\in\V(X,L)$ be a vector such that every formal disk in $k$ has a boundary component mapping to $\hat L$. Let the $(\pm\infty)$-boundary of $(X,L)$ be $(Y^\pm,\Lambda^\pm)$, with Hamiltonian vectors $h^\pm$ and potential vectors $f^\pm$ in $\V^+(Y^\pm\times\R,\Lambda^\pm\times\R)$. Since the potential vector of a trivial cobordism acts as the identity in gluing pairings we suppress it from the notation for split gluing pairings writing simply
\begin{align*}
&\bigl\{k\to h^+\bigr\}\!(u,v):=\bigl\{k\to h^+\bigr\}\!(u,v,f^+),\\
&\bigl\{h^-\leftarrow k\bigr\}\!(u,v):=\bigl\{k\to h^+\bigr\}\!(f^-,u,v).
\end{align*}

Define the {\em $k$-increment} of $u$ inductively as follows.
First let
$$
\rho^{(1)}(u)=d^u(k)=\bigl\{k\to h^+\bigr\}\!\left(u,u\right)
+\bigl\{h^-\leftarrow k\bigr\}\!\left(u,u\right)
$$
(where, as usual, $d^u(k)$ is defined by $h(u)+h(u+\epsilon k)=\epsilon\,d^u(k) +\Ordo(\epsilon^2)$). For $j>1$, define inductively,
\begin{align}\label{e:defrho(j)}
&\rho^{(j)}(u)=\bigl\{k\to h^+\bigr\}\!\left(u,u+\rho^{(j-1)}(u)\right)
+\bigl\{h^-\leftarrow k\bigr\}\!\left(u,u+\rho^{(j-1)}(u)\right),\\\notag
&\delta^{(j)}(u,k)=\rho^{(j)}(u)+\rho^{(j-1)}(u), \quad j>1,
\end{align}
where we take $\rho^{(0)}(u)=0$. Finally, the $k$-increment of $u$ is
\begin{equation}\label{e:dfnDelta}
\Delta(u,k)=\sum_{j=1}^\infty \delta^{(j)}(u,k).
\end{equation}

\begin{rmk}\label{r:primary}
Consider a formal disk contributing to $\delta^{(j)}(u,k)$. By construction, there is a $h^\pm$-factor to which one $k$-factor, and at least one formal disk from $\delta^{(j-1)}(u,k)$ is attached in the $j^{\rm th}$ step of the construction. We say that this $k$-factor is the {\em primary $k$-factor} of the $h^\pm$-factor. It follows from the inductive construction of $\Delta(u,k)$ that every $h^\pm$-factor in a formal disk contributing to $\Delta(u,k)$ has a unique primary $k$-factor.
\end{rmk}

We next relate the above algebraic construction of increment disks to geometry. As explained in Appendix \ref{A:B}, in a generic $1$-parameter family $(X_s,L_s)$, $0\le s\le 1$, of exact cobordisms there is a finite number of distinct moments $\hat s\in[0,1]$ for which there exists a holomorphic disk $\hat k$ of formal dimension $-1$ in $(X_{\hat s},L_{\hat s})$ (we will call such a disk a $(-1)$-disk below), and the potential vector of $(X_s,L_s)$ changes only when $s$ passes one of these isolated moments. Below we describe algebraically how the passage of one such $(-1)$-disk moment affects the potential vector. The analytical study of moduli spaces behind this algebraic description involves solving transversality problems by introducing, so called, abstract perturbations. In particular, such perturbations give rise to new $(-1)$-disks near the original one. The perturbation scheme is described in Subsections \ref{s:mix} and \ref{s:pure}.

Consider a $1$-parameter family of exact cobordisms $(X_s,L_s)$, $0\le s\le 1$, with a $(-1)$-disk $\hat k$ at $s=\frac12$ as in Subsection \ref{s:mixpureini} which maps some boundary component to the piece $\hat L$ of $L$. Let $k\in\V^+(X_s,L_s)$ be the following vector. If $\hat k$ is mixed, as in Subsection \ref{s:mix}, then $k=\hat k$, if $\hat k$ is pure then $k$ is the vector of $(-1)$ disks which arises from the perturbation described in Subsection \ref{s:pure}. Note that all formal disks in $k$ have some boundary component mapping to $\hat L$, see Remark \ref{r:(-1)hatLcomp}.
\begin{lma}\label{l:Delta}
The potential vectors $f_s\in\V^+(X_s,L_s)$, $s=0,1$, are related as follows
$$
f_1 =f_0+ \Delta(f_0,k):=\phi(f_0).
$$
\end{lma}

\begin{pf}
This is a restatement of  Lemmas \ref{l:dfmix} and \ref{l:dfpure}.
\end{pf}
We next consider $1$-parameter families obtained by joining stationary and moving cobordisms. Let $(X^b_s,L^b_s)$, $0\le s\le 1$, be a $1$-parameter family with a $(-1)$-disk $\hat k$ as in Subsection \ref{s:mixpureini}, with a $(+\infty)$-boundary $(Y^0,\Lambda^0)$ and let $(X^a,L^a)$ be an exact cobordism with a  $(-\infty)$-boundary $(Y^0,\Lambda^0)$. Let $(X^{ba}_s,L^{ba}_s)$ be the $1$-parameter family obtained by joining the two, let $f^b_s$ and $f^a$ denote the potential vectors of $(X^b_s,L^b_s)$ and $(X^a,L^a)$, respectively, and let $k\in\V^+(X^b_s,L^b_s)$ be the vector of $(-1)$-disks as in Subsection \ref{s:mix} if $\hat k$ is mixed or as in Subsection \ref{s:pure} if $\hat k$ is pure. Define the vector $K\in\V^+(X^{ba}_s,L^{ba}_s)$ as
$$
K=\bigl\{k\to f^a\bigr\}\!\left(f^b_0,f^b_1,f^a\right)=
\bigl\{k\to f^a\bigr\}\!\left(f^b_0,f^b_0+\Delta(f^b_0,k),f^a\right).
$$
If $\hat k$ has a boundary component mapping to the piece $\hat L^b_s$ of $L^b_s$ then every disk contributing to $k$ has too, see Remark \ref{r:(-1)hatLcomp}. Consequently, if $\hat L^{ba}_s$ is the piece of $L^{ba}_s$ corresponding to $\hat L^b$ then every disk contributing to $K$ has some boundary component mapping to $\hat L^{ba}_s$.
\begin{lma}\label{l:Deltapureb_s}
The potential vectors $F_s$ of $(X^{ba}_s,L^{ba}_s)$, $s=0,1$, are related by
$$
F_1=F_0+\Delta(F_0,K):=\Phi(F_0).
$$
\end{lma}
\begin{pf}
Consider first the case when $\hat k$ is pure. It is straightforward to extend the ordering perturbation which originates in polyfold charts near holomorphic disks in $(Y^0\times\R,\Lambda^0\times\R)$ and continues to such charts of (broken) disks mapping into $(X^b_s,L^b_s)$ as described in Subsection \ref{s:pure}, to (broken) disks which map into $(X^{ba}_s,L^{ba}_s)$ in such a way that the perturbation-order of the negative $\hat L^a$-punctures agrees with the induced order described in Remark \ref{r:order}. It follows that $K$ counts $(-1)$-disks in the glued cobordism and an argument similar to the proof of Lemma \ref{l:polygen} shows that the glued $1$-parameter family is generic in the sense that all its parameterized moduli spaces of dimension $\le 1$ are transversally cut out. (The polyfold differential is onto also at broken disks, see Subsection \ref{s:pure}.)  In this case, the lemma then follows by the proof of Lemma \ref{l:dfpure}. The case when $\hat k$ is mixed follows from a simpler argument along the same lines.
\end{pf}

Similarly, let $(X^b,L^b)$ be an exact cobordism with a $(+\infty)$-boundary $(Y^0,\Lambda^0)$ and let $(X^a_s,L^a_s)$, $0\le s\le 1$, be a $1$-parameter family of cobordisms with a $(-1)$-disk $\hat k$ as in Subsection \ref{s:mixpureini},  with a $(-\infty)$-boundary $(Y^0,\Lambda^0)$. Let $(X^{ba}_s,L^{ba}_s)$ be the $1$-parameter family obtained by joining the two, let $f^b$ and $f^a_s$ denote the potential vectors of $(X^b,L^b)$ and $(X^a_s,L^a_s)$, respectively, and let $k\in\V^+(X^a_s,L^a_s)$ be the vector of $(-1)$-disks constructed in Subsection \ref{s:mix} or \ref{s:pure} according to whether $\hat k$ is mixed or pure. Define the vector $K\in\V^+(X^{ba}_s,L^{ba}_s)$ as
$$
K=\bigl\{f^b\leftarrow k\bigr\}\!\left(f^b,f^a_0,f^a_1\right)
=\bigl\{f^b\leftarrow k\bigr\}\!\left(f^b,f^a_0,f^a_0+\Delta(f^a_0,k)\right).
$$
If $\hat k$ has a boundary component mapping to the piece $\hat L^a_s$ of $L^a_s$ then every disk contributing to $k$ has too, see Remark \ref{r:(-1)hatLcomp}. Consequently, if $\hat L^{ba}_s$ is the piece of $L^{ba}_s$ corresponding to $\hat L^b$ then every disk contributing to $K$ has some boundary component mapping to $\hat L^{ba}_s$.
\begin{lma}\label{l:Deltapurea_s}
The potential vectors $F_s$ of $(X^{ba}_s,L^{ba}_s)$, $s=0,1$, are related by
$$
F_1=F_0+\Delta(F_0,K):=\Psi(F_0).
$$
\end{lma}
\begin{pf}
Analogous to Lemma \ref{l:Deltapureb_s}
\end{pf}

\subsection{Chain isomorphisms}\label{s:chiso}
Consider a $1$-parameter family $(X^b_s,L^b_s)$, $0\le s\le 1$, with a $(+\infty)$-boundary $(Y^0,\Lambda^0)$ and a stationary cobordism $(X^a,L^a)$ with a $(-\infty)$-boundary $(Y^0,\Lambda^0)$ as in Lemma \ref{l:Deltapureb_s}. We use notation as there. Considering $F_0$ in the definition of $\Phi(F_0)$ as a variable we obtain a function $\Phi\colon\V(X^{ba}_0,L^{ba}_0;\epsilon)\to\V(X^{ba}_1,L^{ba}_1;\epsilon)$. For $F\in\V(X^{ba}_0,L^{ba}_0)$, let $\left[\pa_F\Phi\right]\colon\V(X^{ba}_0,L^{ba}_0)\to\V(X^{ba}_1,L^{ba}_1)$ denote the linearization of this function at $F$ defined by
$$
\Phi(F+\epsilon V)= \Phi(F) +\epsilon\,\left[\pa_F\Phi\right](V)
+\Ordo(\epsilon^2)\in\V(X^{ba}_1,L^{ba}_1;\epsilon),
$$
for $V\in\V(X^{ba}_0,L^{ba}_0)$.

\begin{lma}\label{l:chiso}
The map $\left[\pa_{F_0}\Phi\right]\colon\V(X^{ba}_0,L^{ba}_0)\to \V(X^{ba}_1,L^{ba}_1)$ is a $(+)$-action non-decreasing, filtration preserving chain map,
\begin{equation}\label{e:chiso}
d^{F_1}\circ\left[\pa_{F_0}\Phi\right]=\left[\pa_{F_0}\Phi\right]\circ d^{F_0}.
\end{equation}
Furthermore $\left[\pa_{F_0}\Phi\right]$ induces an isomorphism of spectral sequences,
$$
\begin{CD}
\bigl\{E^{p,q}_{r\,;[\alpha]}(X^{ba}_0,L^{ba}_0)\bigr\}
@>{\left[\pa_{F_0}\Phi\right]_\ast}>{\cong}>
\bigl\{E^{p,q}_{r\,;[\alpha]}(X^{ba}_1,L^{ba}_1)\bigr\}
\end{CD},
$$
for every $\alpha>0$.
\end{lma}
\begin{pf}
For ordering purposes fix one piece $\hat L^{ba}_s\subset L^{ba}_s$ to which the original $(-1)$-disk $\hat k$ maps some boundary component (if $\hat k$ is pure then $\hat L^{ba}_s$ is unique). Let $h^\pm$ denote the Hamiltonian vectors at the ends of $(X^{ba}_s,L^{ba}_s)$ and let $H$ denote the corresponding Hamiltonian operator. Let $V\in\V(X^{ba}_0,L^{ba}_0)$. We rewrite the operators in \eqref{e:chiso}, acting on $V$, as
\begin{equation}\label{e:chiso1}
\epsilon\, d^{F_1}\Bigl(\left[\pa_{F_0}\Phi\right](V)\Bigr)=
H\Bigl((F_0+\epsilon V)+\Delta\bigl((F_0+\epsilon V),K\bigr)\Bigr) +\Ordo(\epsilon^2).
\end{equation}
and
\begin{equation}\label{e:chiso2}
\epsilon \left[\pa_{F_0}\Phi\right]\Bigl(d^{F_0}(V)\Bigr)=
H\bigl(F_0+\epsilon V\bigr)+\Delta\bigl(F_0+H(F_0+\epsilon V), K\bigr)+ \Delta\bigl(F_0,K\bigr) +\Ordo(\epsilon^2).
\end{equation}
As in the proof of Lemma \ref{l:diffconc} we find that \eqref{e:chiso} follows once we show that the order $\epsilon$ terms in the right hand sides of \eqref{e:chiso1} and \eqref{e:chiso2} agree in the special case that $V$ is a single admissible formal disk. The order $\epsilon$ term of formal disks without $K$-factors in \eqref{e:chiso1} cancels with the corresponding term in \eqref{e:chiso2}. Thus, the lemma follows if we show that the formal disks of degree one in $\epsilon$ (i.e., with exactly one $V$-factor) and with at least one $K$-factor in the following expression
\begin{equation*}
H\Bigl((F_0+\epsilon V)+\Delta\bigl((F_0+\epsilon V),K\bigr)\Bigr)+
\Delta\Bigl(F_0+H(F_0+\epsilon V),K\Bigr),
\end{equation*}
cancel out, or equivalently, that the corresponding disks in the expression
\begin{equation}\label{e:chiso12}
H\Bigl((F_0+\epsilon V)+\Delta\bigl(F_0,K\bigr)\Bigr)+
H\Bigl(F_0+\Delta\bigl((F_0+\epsilon V),K\bigr)\Bigr)+
\Delta\Bigl(F_0+H(F_0+\epsilon V),K\Bigr),
\end{equation}
cancel out. In any formal disk with one $V$-factor and at least one $K$-factor from \eqref{e:chiso12}, there is exactly one $h^\pm$-factor which does not have a primary $K$-factor, see Remark \ref{r:primary}. (This is the $h^\pm$-factor corresponding to $H$ in all three summands.) We say that a formal disk is an {\em $\lt$-disk} if it has the following properties.
\begin{itemize}
\item[$(i)$] The disk has one $V$-factor, $n>0$ $K$-factors, $(n+1)$ $h^\pm$-factors, and an arbitrary number of $F_0$-factors.
\item[$(ii)$] All $h^\pm$-factors except one have a primary $K$-factor.
\item[$(iii)$] Any $h^\pm$-factor ${\widetilde h}$ with primary $K$-factor has the following property. Let $q$ be a puncture of ${\widetilde h}$ at the end where the primary $K$-factor is attached. If $q$ lies after the primary $K$-factor in the order described in Remark \ref{r:order}, then either an $F_0$-factor or $V$ is attached at $q$.
\end{itemize}
Note that all disks contributing to \eqref{e:chiso12} are $\lt$-disks. We distinguish the $h^\pm$-factor without primary $K$-factor in an $\lt$-disk and call it {\em charged}. Consider an end (all positive- or all negative punctures) of the charged $h^\pm$-factor of an $\lt$-disk. We say that it is {\em obstructed} if one of the following holds.
\begin{itemize}
\item There are no factors connected to the charged $h^\pm$-factor at any puncture in the end.
\item At one puncture in the end the $V$-factor is attached, at all other punctures in the end $F_0$-factors are attached.
\end{itemize}
We say that an $\lt$-disk is {\em isolated}, {\em boundary}, or {\em interior} if its charged $h^\pm$-factor has two, one, or zero obstructed ends, respectively. We first show that the sum of all isolated $\lt$-disks in \eqref{e:chiso12} vanishes. We distinguish four types of such disks.
\begin{itemize}
\item[$({\rm I})$] An isolated $\lt$-disk of type $({\rm I})$ comes from the first term of \eqref{e:chiso12}.
\item[$({\rm II})$] An isolated $\lt$-disk of type $({\rm II})$ comes from the second term of \eqref{e:chiso12}.
\item[$({\rm III}a)$] An isolated $\lt$-disk of type $({\rm III}a)$ comes from the third term of \eqref{e:chiso12} and its charged $h^\pm$-factor is the only $h^{\pm}$-factor connecting to $V$.
\item[$({\rm III}b)$] An isolated $\lt$-disk of type $({\rm III}b)$ comes from the third term of \eqref{e:chiso12} and its charged $h^\pm$-factor is not the only $h^{\pm}$-factor connected to $V$.
\end{itemize}
It is straightforward to check that the isolated $\lt$-disks cancel according to the following, where $X\leftrightarrow Y$ means that disks of type $X$ and disks of type $Y$ cancel,
$$
\begin{matrix}
({\rm I})\leftrightarrow({\rm III}a), &
({\rm II})\leftrightarrow({\rm III}b)
\end{matrix}.
$$

We consider the following propagation law for a non-isolated disk with one of the un-obstructed ends of its charged $h^\pm$-factor ${\overline h}$ active.
\begin{itemize}
\item If there are only $F_0$-factors attached at the active end of ${\overline h}$, then consider the broken disk consisting of ${\overline h}$ and all the attached $F_0$-factors as one end of a $1$-dimensional moduli space of holomorphic disks in $(X^{ba}_0,L^{ba}_0)$. Move to the other end of the moduli space and find a similar breaking. Charge the new $h^\pm$-factor and activate the end of it where the new $F_0$-factors are not attached.
\item If there are not only $F_0$-factors attached at the end, then let $q$ be the last puncture (in the order described in Remark \ref{r:order}) at the active end of ${\overline h}$ where an $h^\pm$-factor or a $K$-factor is attached.
\begin{itemize}
\item If there is an $h^\pm$-factor ${\widetilde h}$ attached at $q$, then ${\widetilde h}$ has a primary $K$-factor $K'$. View the broken disk consisting of ${\overline h}$ and ${\widetilde h}$ as one end of a $1$-dimensional reduced moduli space of holomorphic disks in $(Y^\pm\times\R,\Lambda^\pm\times\R)$. Move to the other end of the reduced moduli space and find a similar breaking into two $h^\pm$-factors. Charge the new $h^\pm$-factor which is not connected to $K'$ and let $K'$ be the primary $K$-factor of the new non-charged $h^\pm$-factor. Activate the end of the charged $h^\pm$-factor which is not connected to the new non-charged $h^\pm$-factor.
\item If there is a $K$-factor $K'$ attached at $q$, then $K'$ is the primary $K$-factor of some $h^\pm$-factor ${\widetilde h}$. Let the charge flow through $K'$. Charge ${\widetilde h}$ and let $K'$ be the primary $K$-factor of ${\overline h}$. Activate the end of ${\widetilde h}$ where $K'$ is not attached.
\end{itemize}
\end{itemize}
In all cases the result is a new non-isolated $\lt$-disk. Using the propagation law starting at any boundary $\lt$-disk we propagate in a unique way until we meet another boundary $\lt$-disk. It is easy to see that any non-isolated $\lt$-disk which contributes to \eqref{e:chiso12} is a boundary $\lt$-disk. We claim that the converse holds as well: any boundary $\lt$-disk contributes to \eqref{e:chiso12}. This in combination with the propagation law establishes \eqref{e:chiso}.

To see that the claim is true, consider a factor $W$ of an $\lt$-disk and let $p$ be one of its punctures. If the sub-disk cut off by cutting at $p$, which does not contain $W$, contains a $K$-factor then we say that $p$ is a {\em $K$-connection}. (If an $\lt$-disk is admissible and $W$ is one of its factors then $W$ can have at most two mixed $K$-connections: if it has more than two we find a sub-disk of the $\lt$-disk which has at least three mixed $\hat L^{ba}_s$-punctures. This is impossible by Remark \ref{r:no3mixedpunct}.) Define an $\lt'$-disk to have properties as an $\lt$-disk with the following modifications. The number of $h^\pm$-factors in $(i)$ is $n$ instead of $(n+1)$, the number of $V$-factors is one or zero instead of one, and all $h^\pm$-factors in $(ii)$ have a primary $K$-factor instead of all except one. We show that any $\lt'$-disk contributes to the constant- or the linear term in
\begin{equation}\label{e:lt'}
\Delta(F_0+\epsilon V, K),
\end{equation}
as follows. Let $K_1$ be the first $K$-factor in the $\lt'$-disk in the order of Remark \ref{r:order}. Let $h_1$ be the primary $h^\pm$-factor of $K_1$ with $K_1$ attached at the end $E$ of $h_1$. Since $K_1$ is the first $K$-factor it follows that no puncture in $E$ before the puncture where $K_1$ is attached is a $K$-connection. By definition there are only $F_0$- and $V$-factors attached after $K_1$. Thus the sub-disk $D_1$ consisting of $h_1$, $K_1$, and all $F_0$- and $V$-factors at $E$ contributes to $\Delta(F_0+\epsilon V,K)$ and appears in $\rho^{(1)}(F_0+\epsilon V)$, see \eqref{e:defrho(j)}. Let $K_2$ be the $K$-factor following $K_1$ in the order. Assume that it is attached at the end $E$ of the $h^\pm$-factor $h_2$. If some puncture after $K_2$ in $E$ is a $K$-connection then $D_1$ must be attached at this puncture and if that is the case then $h_2$, $K_2$, $D_1$, and all $F_0$- and $V$-factors attached at $E$ constitute a sub-disk $D_2$ which contributes to $\Delta(F_0+\epsilon V,K)$ and appears in $\rho^{(2)}(F_0+\epsilon V)$. If no puncture after $K_2$ in $E$ is a $K$-connection then $h_2$, $K_2$, and all $F_0$- and $V$-factors attached at $E$ constitute a sub-disk $D_2$ which contributes to $\Delta(F_0+\epsilon V,K)$ and appears in $\rho^{(1)}(F_0+\epsilon V)$. Continuing in this manner until we reach the last $K$-factor we find that any $\lt'$-disk comes from \eqref{e:lt'}.

Consider a boundary $\lt$-disk with charged $h^\pm$-factor ${\overline h}$. Assume first that there are no disks attached at the obstructed end of ${\overline h}$. Cutting at all the punctures in the un-obstructed end of ${\overline h}$ we find that the sub-disks which arise are either $F_0$-factors, $V$-factors, or $\lt'$-disks. Such a boundary $\lt$-disk contributes to first terms of \eqref{e:chiso12} if ${\overline h}$ is the only $h^\pm$-factor attached to $V$ and to the second term of \eqref{e:chiso12} otherwise. Assume second that the $V$-factor and only $F_0$-factors are attached at the obstructed end of ${\overline h}$. Let $W$ be the sub-disk consisting of ${\overline h}$ and all the factors attached at its obstructed end. Then $W$ contributes to the linear term of $H(F_0+\epsilon V)$. Considering $W$ as one factor of the formal disk, this formal disk fulfills the requirements of an $\lt'$-disk (if $V$ in the definition is replaced by $W$). It follows that such a boundary $\lt$-disk contributes to the third term in \eqref{e:chiso12}. This completes the proof of \eqref{e:chiso}.

The fact that $\left[\pa_{F_0}\Phi\right]$ respects the filtration follows, as usual, from the fact that the gluing operation does not decrease the number of positive punctures. To see that $\left[\pa_{F_0}\Phi\right]$ does not decrease $(+)$-action we note that all factors in $\Delta(F_0,K)$ have holomorphic representatives and hence have positive action, and argue as in the proof of Lemma \ref{l:diffconc}. It follows that $\left[\pa_{F_0}\Phi\right]$ induces a morphism of spectral sequences
$$
\left[\pa_{F_0}\Phi\right]_\ast\colon \bigl\{E^{p,q}_{r\,;[\alpha]}(X^{ba}_0,L^{ba}_0)\bigr\}\to
\bigl\{E^{p,q}_{r\,;[\alpha]}(X^{ba}_1,L^{ba}_1)\bigr\},
$$
for any $\alpha>0$. It remains to show that these maps are isomorphism. Fix $\alpha$ and let $\left[\pa_{F_0}\Phi\right]_\alpha\colon \V_{[\alpha]}(X^{ba}_0,L^{ba}_0)\to \V_{[\alpha]}(X^{ba}_1,L^{ba}_1)$ be the map induced by $\left[\pa_{F_0}\Phi\right]$. We use a decomposition
\begin{equation}\label{e:actdecomp}
\left[\pa_{F_0}\Phi\right]=\id+B,
\end{equation}
where  $B\colon\V(X^{ba}_0,L^{ba}_0)\to\V(X^{ba}_1,L^{ba}_1)$ strictly increases $(+)$-action. To see that there is such a decomposition consider a formal disk $W$ from $\delta^{(j)}(F_0+\epsilon V)$ contributing to $\Delta(F_0+\epsilon V,K)$. We claim that if $W'$ is sub-disk of $W$ which arises when the last $h^\pm$-factor ${\widetilde h}$ (i.e., the $h^\pm$-factor used in the $j^{\rm th}$ step in the construction of $\rho^{(j)}(F_0+\epsilon V)$) is cut-off, then $\A^+(W)>\A^+(W')$. To prove the claim assume first that ${\widetilde h}$ is of type $h^+$. Then since ${\widetilde h}$ is holomorphic $\A({\widetilde h})>\A^+({\widetilde h})-\A^-({\widetilde h})>0$ which immediately gives $\A^+(W)>\A^+(W')$. Assume second that ${\widetilde h}$ is of type $h^-$. Then the positive punctures of $W'$ is a strict subset of the positive punctures of $W$ (since the positive puncture of the primary $K$-factor of ${\widetilde h}$ is a positive puncture of $W$ but not of $W'$). In either case the claim follows. The claim implies that the $(+)$-action of any disk contributing to the first order term in $\Delta(F_0+\epsilon V, K)$ has strictly larger $(+)$-action than $V$. Defining $B(V)$ as the first order term of $\Delta(F_0+\epsilon V, K)$, \eqref{e:actdecomp} follows.

Our genericity assumptions imply that the actions of the Reeb chords of $L$ form a discrete subset of $\R$. In particular, there exists $\alpha_0>0$ such that if $U$ is any admissible disk with $\A^+(U)\le \alpha$ and if $W$ is any disk contributing to $B(U)$ then $\A^+(W)-\A^+(U)>\alpha_0$. It follows that the map $B_\alpha\colon \V_{[\alpha]}(X^{ba}_0,L^{ba}_0)\to \V_{[\alpha]}(X^{ba}_1,L^{ba}_1)$ induced by $B$ increases $(+)$-action by at least $\alpha_0$. Fix $N=2^M$ such that $N\alpha_0>\alpha$. We have (because of $\Z_2$-coefficients)
$$
(\id+B)\,\circ\stackrel{N}{\cdots}\circ\,(\id+B)=\id+ B^N.
$$
Since $B_\alpha^N$ increases $(+)$-action by at least $N\alpha_0>\alpha$, we find that $B_\alpha^N=0$. We conclude that $\left[\pa_{F_0}\Phi\right]_\alpha=\id+B_\alpha$ satisfies $\left[\pa_{F_0}\Phi\right]^N_\alpha=\id$. In particular, the filtration preserving chain map $\left[\pa_{F_0}\Phi\right]_\alpha\colon \V_{[\alpha]}(X^{ba}_0,L^{ba}_0)\to\V_{[\alpha]}(X^{ba}_1,L^{ba}_1)$ has a filtration preserving inverse $\left[\pa_{F_0}\Phi\right]_\alpha^{N-1}\colon \V_{[\alpha]}(X^{ba}_1,L^{ba}_1)\to\V_{[\alpha]}(X^{ba}_0,L^{ba}_0)$ (we use the fact that $\V_{[\alpha]}(X^{ba}_0,L^{ba}_0)$ and $\V_{[\alpha]}(X^{ba}_1,L^{ba}_1)$ are canonically isomorphic as vector spaces) which is easily seen to be a chain map. It follows that $\left[\pa_{F_0}\Phi\right]_\alpha$ induces an isomorphism of spectral sequences.
\end{pf}

Consider a $1$-parameter family $(X^a_s,L^a_s)$, $0\le s\le 1$, with a negative end at $(Y^0,\Lambda^0)$ and a stationary cobordism $(X^b,L^b)$ with positive end at $(Y^0,\Lambda^0)$. Considering $F_0$ in the definition of $\Psi(F_0)$, see Lemma \ref{l:Deltapurea_s}, as a variable we obtain a function $\Psi\colon\V(X^{ba}_0,L^{ba}_0;\epsilon)\to\V(X^{ba}_1,L^{ba}_1;\epsilon)$. For $F\in\V(X^{ba}_0,L^{ba}_0)$, let $\left[\pa_F\Psi\right]\colon\V(X^{ba}_0,L^{ba}_0)\to\V(X^{ba}_1,L^{ba}_1)$ denote the linearization of this function at $F$ defined by
$$
\Psi(F+\epsilon V)= \Psi(F) +\epsilon\,\left[\pa_F\Psi\right](V)+\Ordo(\epsilon^2)
\in\V(X^{ba}_1,L^{ba}_1;\epsilon).
$$

\begin{lma}\label{l:chiso'}
The map $\left[\pa_{F_0}\Psi\right]$ is a $(+)$-action non-decreasing, filtration preserving chain map,
\begin{equation}\label{e:chiso'}
d^{F_1}\circ\left[\pa_{F_0}\Psi\right]=\left[\pa_{F_0}\Psi\right]\circ d^{F_0}.
\end{equation}
Furthermore $\left[\pa_{F_0}\Psi\right]$ induces an isomorphism of spectral sequences,
$$
\begin{CD}
\bigl\{E^{p,q}_{r\,;[\alpha]}(X^{ba}_0,L^{ba}_0)\bigr\}
@>{\left[\pa_{F_0}\Psi\right]_\ast}>{\cong}>
\bigl\{E^{p,q}_{r\,;[\alpha]}(X^{ba}_1,L^{ba}_1)\bigr\}
\end{CD},
$$
for every $\alpha>0$.
\end{lma}
\begin{pf}
Analogous to Lemma \ref{l:chiso}.
\end{pf}

\subsection{Chain homotopy}\label{s:chhom}
We consider the same situations as in Subsection \ref{s:chiso}, now focusing on how the deformation affects chain maps. We will first treat the case when a stationary cobordism is joined to a moving one from above and second the case when the stationary cobordism is joined to the moving one from below.

Let $(X^a,L^a)$ be a stationary cobordism joined to a moving one $(X^b_s,L^b_s)$, $0\le s\le 1$, from above along $(Y^0,\Lambda^0)$ forming a new $1$-parameter family $(X^{ba}_s,L^{ba}_s)$. With notation as in Subsection \ref{s:chiso}, we have
\begin{align}\label{e:F_1exp}
F_1&=\Phi(F_0)=F_0+\Delta(F_0,K)\\\label{e:F_1exp'}
&=\left(f^b_0+\Delta(f^b_0,k)\I f^a \right),
\end{align}
where $K\in\V^+(X^{ba}_s,L^{ba}_s)$ is given by
\begin{equation}\label{e:Kexp}
K=\bigl\{k\to f^a\bigr\}\!\left(f^b_0, f^b_1,f^a\right)=
\bigl\{k\to f^a\bigr\}\!\left(f^b_0, f^b_0+\Delta(f^b_0,k),f^a\right).
\end{equation}
Subsections  \ref{s:cob->chmap} and \ref{s:chiso} give the following diagram of chain maps
$$
\begin{CD}
\V(X^a,L^a) @>{\left[\pa_{f^a}\!\left(f^b_0\right.\I\right]}>> \V(X^{ba}_0,L^{ba}_0)\\
@V{\id}VV  @VV{\left[\pa_{F_0}\Phi\right]}V\\
\V(X^a,L^a) @>{\left[\pa_{f^a}\!\left(f^b_1\right.\I\right]}>> \V(X^{ba}_1,L^{ba}_1)
\end{CD}.
$$
Starting at the upper left corner of this diagram, going right and then down corresponds to linearizing \eqref{e:F_1exp} with respect to $f^a$-factors in $F_0$. Starting at the upper left corner, going down and then right corresponds to linearizing \eqref{e:F_1exp'} with respect to $f^a$. We show below that these two filtration preserving chain maps induce the same morphism of spectral sequences. In order to do so we will make use of chain homotopies which are built from the following maps constructed using the split gluing pairing, see Subsection \ref{s:incrdisk}. For $v\in\V(X^a,L^a)$, define
$$
\theta(v)=\bigl\{k\to v\bigr\}\!\left(f^b_0,f^b_0+\Delta(f^b_0,k),0\right)\in\V(X^{ba}_s,L^{ba}_s).
$$

Let $W\in\V(X^{ba}_s,L^{ba}_s;\epsilon)$ and let $K$ be as above. Recall the inductive construction of the $K$-increment disks of $F\in\V(X^{ba}_0,L^{ba}_0)$. We consider a deformed version of that construction as follows. With $h^+$ and $h^-$ denoting the Hamiltonian vectors at the positive and negative ends of $(X^{ba}_s,L^{ba}_s)$, define
$$
\rho_W^{(1)}(F)= d^F(K) + W=\bigl\{K\to h^+\bigr\}\!\left(F,F\right)+
\bigl\{h^-\leftarrow K\bigr\}\!\left(F,F\right)+W,
$$
and inductively
$$
\rho^{(j)}_W(F)=\bigl\{K\to h^+\bigr\}\!\left(F,F+\rho_W^{(j-1)}(F)\right)+
\bigl\{h^-\leftarrow K\bigr\}\!\left(F,F+\rho_W^{(j-1)}(F)\right)+W.
$$
Let
$$
\delta_W^{(j)}(F)=\rho^{(j)}_W(F)+\rho^{(j-1)}_W(F),\quad j>0,
$$
where we take $\rho^{(0)}_W(F)=0$, and finally
\begin{equation}\label{e:defDelta}
\Delta_W(F,K)=\sum_{j=1}^\infty \delta^{(j)}_W(F).
\end{equation}

The map we will use to construct chain homotopies is $\Theta\colon\V(X^a,L^a)\to\V(X^{ba}_1,L^{ba}_1)$ defined through the following equation
\begin{equation}\label{e:Theta}
\Delta_{\epsilon\,\theta(v)}(F_0,K)+ \Delta(F_0,K) + \epsilon\,\Theta(v)= \Ordo(\epsilon^2)\in\V(X^{ba}_1,L^{ba}_1;\epsilon).
\end{equation}
Recall the notation $\hat L^{ba}_s$ for the piece of $L^{ba}_s$ where all $K$-disks have a boundary component.

\begin{lma}\label{l:chhom}
The $(+)$-action non-decreasing, filtration preserving chain maps $\left[\pa_{F_0}\Phi\right]\circ\left[\pa_{f^a}\!\left(f^b_0\right.\I\right]$ and $\left[\pa_{f^a}\!\left(f^b_1\right.\I\right]$ induce the same morphism of spectral sequences. In other words the following diagram commutes
\begin{equation}\label{e:chhom}
\begin{CD}
\bigl\{E_{r\,;[\alpha]}^{p,q}(X^a,L^a)\bigr\} @>{\left[\pa_{f^a}\!\left(f^b_0\right.\I\right]_\ast}>> \bigl\{E_{r\,;[\alpha]}^{p,q}(X^{ba}_0,L^{ba}_0)\bigr\}\\
@V{\id}VV  @VV{\left[\pa_{F_0}\Phi\right]_\ast}V\\
\bigl\{E_{r\,;[\alpha]}^{p,q}(X^a,L^a)\bigr\} @>{\left[\pa_{f^a}\!\left(f^b_1\right.\I\right]_\ast}>> \bigl\{E_{r\,;[\alpha]}^{p,q}(X^{ba}_1,L^{ba}_1)\bigr\}
\end{CD},
\end{equation}
for every $\alpha>0$.
\end{lma}

\begin{pf}
Consider the graded vector spaces
\begin{align*}
\CC^p(X^a,L^a)&=F^p\V(X^a,L^a)/F^{p+1}\V(X^a,L^a),\\
\CC^p_{[\alpha]}(X^a,L^a)&=F^p\V_{[\alpha]}(X^a,L^a)/F^{p+1}\V_{[\alpha]}(X^a,L^a),\\
\CC^p(X^{ba}_s,L^{ba}_s)&=F^p\V(X^{ba}_s,L^{ba}_s)/F^{p+1}\V(X^{ba}_s,L^{ba}_s),\\
\CC^p_{[\alpha]}(X^{ba}_s,L^{ba}_s)&=
F^p\V_{[\alpha]}(X^{ba}_s,L^{ba}_s)/F^{p+1}\V_{[\alpha]}(X^{ba}_s,L^{ba}_s),
\end{align*}
where $\alpha>0$ and $s=0,1$. Since the differentials $d^{f_a}$, $d^{f_a}_\alpha$, $d^{F_s}$, and $d^{F_s}_\alpha$ respect the filtration, they induce differentials on $\CC^p(X^a,L^a)$, $\CC^p_{[\alpha]}(X^a,L^a)$,   $\CC^p(X^{ba}_s,L^{ba}_s)$, and $\CC^p_{[\alpha]}(X^{ba}_s,L^{ba}_s)$, respectively, for every $p$, making them graded chain complexes. We denote these induced differentials by the same symbols as the original ones. The homology of the graded complexes $\CC^p_{[\alpha]}(X^a,L^a)$ and $\CC^p_{[\alpha]}(X^{ba}_s,L^{ba}_s)$ are the $E_1$-terms of the spectral sequences of $(X^a,L^a)$ and of $(X^{ba}_s,L^{ba}_s)$, respectively. Furthermore, as the chain maps
$\left[\pa_{F_0}\Phi\right]$, $\left[\pa_{f^a}\!\left(f^b_0\right.\I\right]$, and $\left[\pa_{f^a}\!\left(f^b_1\right.\I\right]$ are all defined using the gluing pairing, they respect filtration and hence induce chain maps
\begin{equation}\label{e:comp1}
\begin{CD}
\CC^p(X^a,L^a) @>{\left[\pa_{f^a}\!\left(f^b_0\right.\I\right]^p}>>
\CC^p(X^{ba}_0,L^{ba}_0)
@>{\left[\pa_{F_0}\Phi\right]^p}>>
\CC^p(X^{ba}_1,L^{ba}_1)
\end{CD}
\end{equation}
and
\begin{equation}\label{e:comp2}
\begin{CD}
\CC^p(X^a,L^a) @>{\left[\pa_{f^a}\!\left(f^b_1\right.\I\right]^p}>>
\CC^p(X^{ba}_1,L^{ba}_1),
\end{CD}
\end{equation}
which are all $(+)$-action non-decreasing. Thus we have induced chain maps
\begin{equation}\label{e:comp'1}
\begin{CD}
\CC^p_{[\alpha]}(X^a,L^a) @>{\left[\pa_{f^a}\!\left(f^b_0\right.\I\right]^p}>>
\CC^p_{[\alpha]}(X^{ba}_0,L^{ba}_0)
@>{\left[\pa_{F_0}\Phi\right]^p}>>
\CC^p_{[\alpha]}(X^{ba}_1,L^{ba}_1)
\end{CD}
\end{equation}
and
\begin{equation}\label{e:comp'2}
\begin{CD}
\CC^p_{[\alpha]}(X^a,L^a) @>{\left[\pa_{f^a}\!\left(f^b_1\right.\I\right]^p}>>
\CC^p_{[\alpha]}(X^{ba}_1,L^{ba}_1),
\end{CD}
\end{equation}
for all $\alpha>0$. We show below that \eqref{e:comp1} and \eqref{e:comp2} are chain homotopic for every $p$, via a $(+)$-action non-decreasing chain homotopy. This means that the chain maps \eqref{e:comp'1} and \eqref{e:comp'2} induce identical maps on the $E_1$-term of the spectral sequences. Since morphisms of spectral sequences which agree on the $E_1$-term agree everywhere this will complete the proof of the lemma. More specifically, $\Theta\colon\V(X^a,L^a)\to\V(X^{ba}_1,L^{ba}_1)$, see \eqref{e:Theta}, does not decrease $(+)$-action (since it is defined by gluing of holomorphic curves, see the proof of Lemma \ref{l:diffconc}) and respects the filtration. Therefore it induces $(+)$-action non-decreasing maps $\Theta^p\colon\CC^p(X^a,L^a)\to\CC^p(X^{ba}_1,L^{ba}_1)$ and we will show that
\begin{equation}\label{e:chhomp}
\left[\pa_{F_0}\Phi\right]^p\circ \left[\pa_{f^a}\!\left(f^b_0\right.\I\right]^p
+\left[\pa_{F_0}\Phi\right]^p=d^{F_1}\circ\Theta^p +\Theta^p\circ d^{f^a},
\end{equation}
for all $p$.  As in the proof of Lemma \ref{l:diffconc}, we note that it is enough to show that the left- and right hand sides of \eqref{e:chhomp} agree in $\CC^p(X^{ba}_1,L^{ba}_1)=F^p\V(X^{ba}_1,L^{ba}_1)/F^{p+1}\V(X^{ba}_1,L^{ba}_1)$, when evaluated on an element $v\in F^p\V(X^a,L^a)$ which is represented by a single formal disk with $p$ positive punctures.

Let $v$ be a formal disk with $p$ positive punctures. We compute
\begin{align}\label{e:chhomT1}
&\epsilon\,\left[\pa_{F_0}\Phi\right]
\Bigl(\left[\pa_{f^a}\!\left(f^b_0\right.\I\right](v)\Bigr)\\\notag
&=\Phi(F_0)+\Phi\Bigl(\left(f_0^b\I f^a+\epsilon v\right)\Bigr)
+ \Ordo(\epsilon^2)\\\notag
&= F_0 +\Delta(F_0,K)+\left(f_0^b\I f_0^a+\epsilon v\right)+
\Delta\Bigl(\left(f_0^b\I f^a_0+\epsilon v\right),K\Bigr)
+ \Ordo(\epsilon^2)\\\notag
&=\left(f_0^b\I f^a\right) +
\Delta\Bigl(\left(f_0^b\I f^a\right),\bigl\{k\to f^a\bigr\}\!\left(f_0^b,f_0^b+\Delta(f_0^b,k),f^a\right)\Bigr)+\left(f_0^b\I f^a+\epsilon v\right)\\\notag
&+
\Delta\Bigl(\left(f_0^b\I f^a+\epsilon v\right),\bigl\{k\to f^a\bigr\}\!\left(f_0^b,f_0^b+\Delta(f_0^b,k),f^a\right)\Bigr)
+ \Ordo(\epsilon^2),
\end{align}
and
\begin{align}\label{e:chhomT2}
\epsilon\,\left[\pa_{f^a}\!\left(f_1^b\right.\I\right](v) &=
\left(f^b_1\I f^a +\epsilon v\right) +
\left(f^b_1\I f^a\right) + \Ordo(\epsilon^2)\\\notag
&=\left(f^b_0+\Delta\left(f^b_0,k\right)\I f^a+\epsilon v\right)+
\left(f^b_0+\Delta\left(f^b_0,k\right)\I f^a\right)+\Ordo(\epsilon^2).
\end{align}
The sum of the formal disks linear in $\epsilon$ in \eqref{e:chhomT1} and \eqref{e:chhomT2} which contain no $k$-factors vanish. Thus, the contributions to the left hand side of \eqref{e:chhomp} all come from formal disks in \eqref{e:chhomT1} and \eqref{e:chhomT2} which have one $v$-factor and at least one $k$-factor. We study such disks more closely taking into account the fact that formal disks with more than $p$ positive punctures do not contribute. We first note that in any contributing disk, each factor except for the $v$-factor has only one positive puncture. This follows from the observation that if $u$ is a factor of a formal disk $w$ and if some factor of $w$ other than $u$ has at least two positive punctures then the number of positive punctures of $w$ is larger than the number of positive punctures of $u$. The contributing disks from \eqref{e:chhomT1} thus all come from
\begin{equation}\label{e:chhomT1'}
\Delta\Bigl(\left(f_0^b\I f^a+\epsilon v\right),\bigl\{k\to f^a\bigr\}\!\left(f_0^b,f_0^b+\Delta(f_0^b,k),0\right)\Bigr),
\end{equation}
where every factor except the $v$-factor has only one positive puncture. (The last $0$ follows since a non-trivial $f^a$-factor implies that the corresponding $K$-factor has at least two positive punctures.) The contributing disks from \eqref{e:chhomT2} come from
\begin{equation}\label{e:chhomT2'}
\left(f^b_0+\Delta\left(f^b_0,k\right)\I \epsilon v\right)
\end{equation}
where every factor except the $v$-factor has only one positive puncture and the absence of the $f^a$-term follows since any formal disk \eqref{e:chhomT2} with a $v$- and an $f^a$-factor has more than $p$ positive punctures.

Let $H$ denote the Hamiltonian operator on $\V(X^{ba}_1,L^{ba}_1)$ and let $h$ denote the Hamiltonian operator on $\V(X^a,L^a)$. The two terms in the right hand side of \eqref{e:chhomp} have the form
\begin{equation}\label{e:chhomT3}
H\Bigl(F_1+\epsilon\,\Theta(v)\Bigr)+\Ordo(\epsilon^2)
\end{equation}
and
\begin{equation}\label{e:chhomT4}
\Theta\bigl(h(f^a+\epsilon v)\bigr)+\Ordo(\epsilon^2).
\end{equation}
As above we find that in any formal disk from \eqref{e:chhomT3} or \eqref{e:chhomT4} which contributes to \eqref{e:chhomp}, all factors except for the $v$-factor has only one positive puncture.

In order to demonstrate how the disks in \eqref{e:chhomT1'} -- \eqref{e:chhomT4} cancel, we first introduce some notation. Let $h^0$ denote the Hamiltonian vector of the trivial cobordism $(Y^0\times\R,\Lambda^0\times\R)$ where the cobordisms $(X^b_s,L^b_s)$ and $(X^a,L^a)$ are joined. Let $h^{a+}$ and $h^{b-}$ denote the Hamiltonian vectors at the positive end of $(X^a,L^a)$ and at the negative end of $(X^b_s,L^b_s)$, respectively. Let $h^{a-}+h^0$ and $h^{b+}+h^0$ be the Hamiltonian vectors at the negative end of $(X^a,L^a)$ and at the positive end of $(X^b_s,L^b_s)$, respectively. Note that the Hamiltonian vectors at the positive- and negative ends of $(X^{ba}_s,L^{ba}_s)$ are $h^+=h^{a+}+h^{b+}$ and $h^-=h^{a-}+h^{b-}$, respectively. We use the collective name {\em $h^\ast$-factor} for a factor from $h^{\pm}$ or $h^0$.

We define an {\em $\lt$-disk} to be a formal disk with factors as follows.
\begin{itemize}
\item One $v$-factor, $n>0$ $h^\ast$-factors, $n>0$ $k$-factors, and arbitrary numbers of $f^a$- and $f^b_0$-factors.
\item All factors except the $v$-factor has only one positive puncture.
\item All of the $n+1$ factors of type $h^\ast$ and $v$, except one, has a primary $k$-factor, if it is of type $h^{b\pm}$, $h^0$, or $v$, or a primary $K$-factor, if it is of type $h^{a\pm}$.
\item Sub-disks obtained by cutting at a negative puncture of a $h^\ast$- or $v$-factor with primary $k$- or $K$-factor which lie {\em after} this primary factor, in the ordering induced as described in Remark \ref{r:order}, contain no $k$-factors.
\end{itemize}
Note that the first three requirements imply that no $\lt$-disk has a $h^{a-}$ factor with a primary $K$-factor, since such a disk must have at least one factor, which is not the $v$-factor, with more than one positive puncture. Similarly, no $\lt$-disk has a $h^{b+}$-factor since if such a factor is to be connected to $v$ there is at least one factor, which is not the $v$-factor, with more than one positive puncture.
Let $q$ be a  puncture of some factor of an $\lt$-disk. If the sub-disk cut-off by cutting at $q$ contains some $k$-factor then we say that $q$ is an {\em $k$-connection}. If one of the negative punctures of a factor of an $\lt$-disk is a $k$-connection then its negative punctures are ordered, see Remark \ref{r:order}.

We will distinguish one $h^\ast$- or $v$-factor in each $\lt$-disk as charged and, in the case it is an $h^\ast$-factor, define two activity modes for the charged disk. The following combinations are possible.
\begin{itemize}
\item[$({\rm a})$] A $v$-factor may be charged.
\item[$({\rm b})$] An $h^\ast$-factor without primary $k$- or $K$-factor may be charged and either its positive- or its negative end may be active.
\item[$({\rm c})$] An $h^0$-factor with primary $k$-factor may be charged if its positive puncture is attached either to the $v$-factor, or to an $f^a$-factor, the positive puncture of which in turn is attached to the last negative $k$-connection of an $h^+$-factor without primary $K$-factor. The positive end of such a charged $h^0$-factor may be either flow-active or glue-active.
\end{itemize}
We define the {\em dimension} of a sub-disk of an $\lt$-disk to be the number of $h^\ast$- and $v$-factors it contains minus the number of $k$-factors it contains. Since all $h^\ast$- and $v$-factors of an $\lt$-disk except one have a primary $K$- or $k$-factor it follows that any sub-disk of an $\lt$-disk obtained by cutting at the negative puncture of some $h^*$-factor has dimension $0$ or $-1$.

Consider the charged factor ${\overline h}$ of an $\lt$-disk of type $({\rm b})$. This is an $h^\ast$-factor. Let $p$ denote the positive puncture of ${\overline h}$ and let $q$ denote the last negative puncture of ${\overline h}$ which is a $k$-connection. We use the following terminology.
\begin{itemize}
\item The positive end of ${\overline h}$ is {\em obstructed} if one of the following holds.
\begin{itemize}
\item There is no disk attached to ${\overline h}$ at $p$,
\item ${\overline h}$ is attached to $v$ at $p$, or
\item ${\overline h}$ is attached to a $k$-factor at $p$ and this $k$-factor is attached to $v$ at its positive puncture.
\end{itemize}
\item The negative end of ${\overline h}$ is {\em obstructed} if one of the following holds.
\begin{itemize}
\item There is no disk attached to ${\overline h}$ at any of its negative punctures,
\item $v$ is attached to ${\overline h}$ at $q$, or
\item the sub-disk cut off by cutting at $q$ has dimension $-1$.
\end{itemize}
\end{itemize}
An end which is not obstructed is called {\em un-obstructed}. An $\lt$-disk with charged factor of type $({\rm b})$ with two, one, or zero ends obstructed is called {\em isolated}, {\em boundary}, or {\em interior}, respectively. An $\lt$-disk with charged factor of type $({\rm a})$ is {\em isolated} and an $\lt$-disk with charged factor of type $({\rm c})$ is {\em interior}.

We view the disks in \eqref{e:chhomT1'} and \eqref{e:chhomT2'} as $\lt$-disks by letting their $v$-factors be charged. They are all isolated. Vi view the disks in \eqref{e:chhomT3} and \eqref{e:chhomT4} as $\lt$-disks by charging the $h^\ast$-factor which corresponds to the operator $H$ and the operator $h$, respectively. We first show that isolated $\lt$-disks cancel.

Isolated $\lt$-disks from \eqref{e:chhomT1'} are of two types
\begin{itemize}
\item[$({\rm I}a)$] The $v$-factor is connected to one $h^+$-factors.
\item[$({\rm I}b)$] The $v$-factor is connected to two $h^+$-factors.
\end{itemize}
It is easy to see that the $v$-factor has at least one $h^+$-factor attached. Any $h^+$-factor has a primary $K$-factor. Thus, a $v$-factor with more than two $h^+$-factors attached would be non-admissible since it has at least three positive $k$-connections, see Remark \ref{r:no3mixedpunct}. Consider a disk $D$ of type $({\rm I}b)$. Since both $h^+$-factors attached to $v$ have primary $K$-factors it follows by admissibility that both the positive punctures of $v$ where $h^+$-factors are attached are mixed $\hat L^{ba}_s$-punctures. Hence one is the first $\hat L^{ba}_s$-puncture of $v$ and the other one the last. Let $h^+_{\rm \,fi}$ and $h^+_{\rm \,la}$ denote the $h^+$-factors attached at the first and last puncture, respectively. We claim that the positive puncture of $h^+_{\rm \,fi}$ is a positive puncture of $D$ as well. To see this we argue as follows. Cutting $D$ at the last $\hat L^{ba}_s$-puncture of $v$ we obtain two sub-disks $D_{\rm \,fi}$ and $D_{\rm \,la}$ where $D_{\rm \,fi}$ contains $v$ and $h^+_{\rm \,fi}$ and where $D_{\rm \,la}$ contains $h^+_{\rm \,la}$. Since $D$ is admissible it follows that the last $\hat L^{ba}_s$-puncture of $D_{\rm \,fi}$ is the last $\hat L^{ba}_s$-puncture of $v$. Assume that the positive puncture of $h^+_{\rm \,fi}$ is attached to some other $h^+$-factor $h^+_{\rm \,ot}$. Since $h^+_{\rm \,ot}$ has a primary $K$-factor and since $D_{\rm \,fi}$ contains a $k$-factor, the positive puncture of $h^+_{\rm \,fi}$ must be attached at a negative puncture of $h^+_{\rm \,ot}$ before the $K$-factor. Furthermore, since $D_{\rm \,fi}$ contains mixed punctures the positive puncture of $h^+_{\rm \,fi}$ is mixed and $h^+_{\rm \,fi}$ must be attached at the first $\hat L^{ba}_s$-puncture in $h^+_{\rm \,ot}$. This, however, is impossible  for orientation reasons: the first $\hat L^{ba}_s$-puncture of $D_{\rm \,fi}$ should be attached to the last $\hat L^{ba}_s$-puncture of $h^+_{\rm \,ot}$-factor. We conclude that the positive puncture of $h^+_{\rm \,fi}$ is also a positive puncture of $D$.

Isolated $\lt$-disks from \eqref{e:chhomT2'} are of two types as well. Let $q$ be the last $k$-connection of the charged $v$-factor.
\begin{itemize}
\item[$({\rm II}a)$] An $h^0$-factor is attached at $q$.
\item[$({\rm II}b)$] A $k$-factor is attached at $q$.
\end{itemize}

Isolated $\lt$-disks from \eqref{e:chhomT3} are of three types. Let ${\overline h}$ denote the charged $h^\ast$-factor, let $p$ be its positive puncture, and let $q$ be its last negative $k$-connection.
\begin{itemize}
\item[$({\rm III}a)$] ${\overline h}$ is of type $h^+$, $v$ is attached at $q$, and $v$ is not attached to any non-charged $h^+$-factor.
\item[$({\rm III}b)$]  ${\overline h}$ is of type $h^+$, $v$ is attached at $q$, and $v$ is attached to one non-charged $h^+$-factor.
\item[$({\rm III}c)$] ${\overline h}$ is of type $h^-$, a $k$-factor $k'$ is attached at $p$, the positive puncture of $k'$ is attached to $v$.
\end{itemize}

Isolated $\lt$-disks from \eqref{e:chhomT4} are of five types. Let ${\overline h}$ denote the charged $h^\ast$-factor, let $p$ be its positive puncture, and let $q$ be its last negative $k$-connection.
\begin{itemize}
\item[$({\rm IV}a)$] ${\overline h}$ is of type $h^+$, $v$ is attached at $q$, and $v$ is not attached to any non-charged $h^+$-factor.
\item[$({\rm IV}b)$] ${\overline h}$ is of type $h^+$, $v$ is attached at $q$, and $v$ is attached to one non-charged $h^+$-factor.
\item[$({\rm IV}c)$] ${\overline h}$ is of type $h^+$, an $f^a$-factor is attached at $q$, the disk cut off by cutting at $q$ has dimension $-1$, and no non-charged $h^+$-factor is attached to $v$.
\item[$({\rm IV}d)$] ${\overline h}$ is of type $h^+$, an $f^a$-factor is attached at $q$, the disk cut off by cutting at $q$ has dimension $-1$, and one non-charged $h^+$-factor is attached to $v$.
\item[$({\rm IV}e)$] ${\overline h}$ is of type $h^0$, $v$ is attached at $p$, and a $k$-factor is attached at $q$.
\end{itemize}

Reinterpreting the factorizations of the isolated $\lt$-disks it is straightforward to check that they cancel according to the following.
$$
\begin{matrix}
({\rm I}a)\leftrightarrow ({\rm IV}c), &
({\rm I}b)\leftrightarrow ({\rm IV}d), &
({\rm II}a)\leftrightarrow ({\rm IV}e), \\
({\rm II}b)\leftrightarrow ({\rm III}c), &
({\rm III}a)\leftrightarrow ({\rm IV}a), &
({\rm III}b)\leftrightarrow ({\rm IV}b).
\end{matrix}
$$

\begin{figure}[htbp]
\begin{center}
\includegraphics[angle=0, width=12cm]{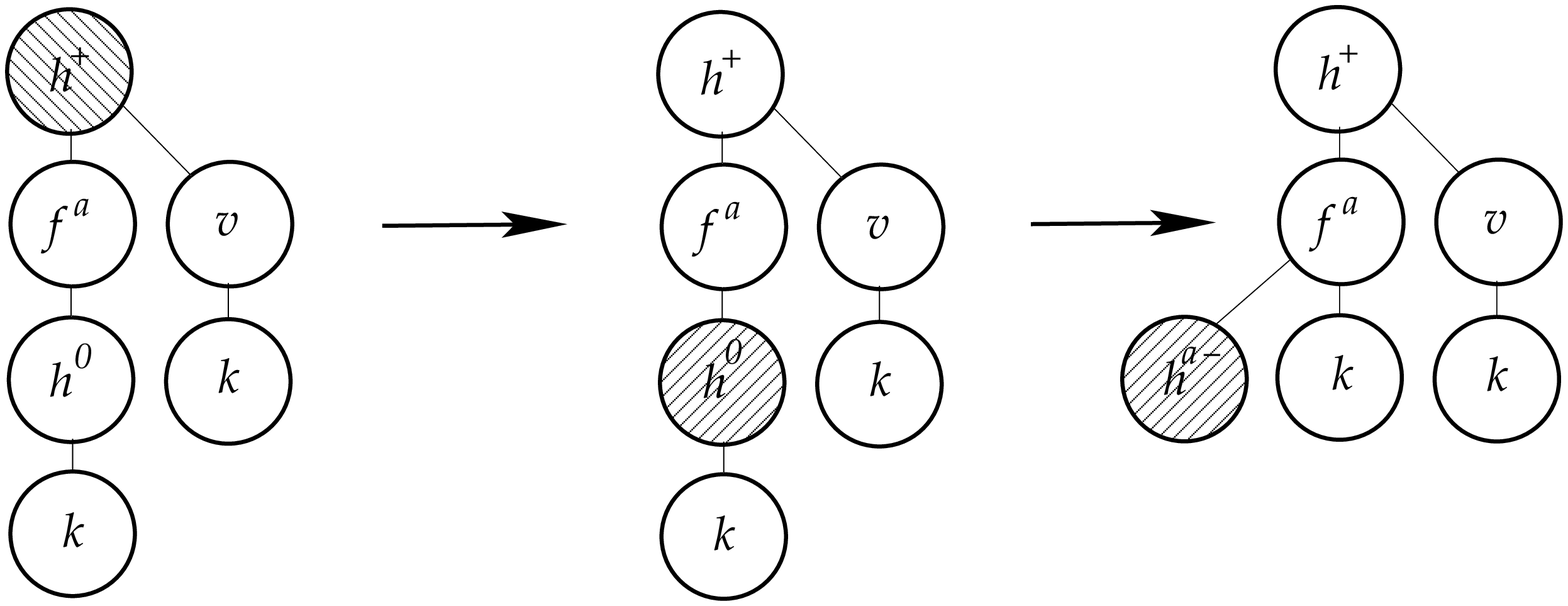}
\end{center}
\caption{The propagation law: start at a boundary $\lt$-disk with $h^+$ charged and with negative end active. Flow the charge to $h^0$ with primary $k$-factor. Its positive end becomes glue-active. Glue $h^0$ to $f^a$. At the other end of the moduli space the disk splits into $f^a$ and $h^{a-}$. This gives the end boundary $\lt$-disk with $h^{a-}$ charged.}
\label{fig:chflow}
\end{figure}

We next define propagation laws for non-isolated $\lt$-disks, with a specified activity mode of its charged factor. Let ${\overline h}$ denote the charged $h^\ast$-factor.
\begin{itemize}
\item
${\overline h}$ is of type $h^+$ and its positive end is active. Then the positive puncture of ${\overline h}$ is attached to some other $h^+$-factor, ${\widetilde h}$ which has a primary $K$-factor $K'$. The broken disk obtained by gluing ${\overline h}$ to ${\widetilde h}$ is the boundary of a $1$-dimensional reduced moduli space of holomorphic disks in $(Y^+\times\R,\Lambda^+\times\R)$. Moving to the other end of this reduced moduli space, we find a similar breaking into two $h^+$-disks. Charge the one of the new $h^+$-factors which is not connected to $K'$ and activate the end of it where it is not connected to the new non-charged $h^+$-factor.
\item ${\overline h}$ is of type $h^+$ and its negative end is active. Consider the last negative $k$-connection $q$ of ${\overline h}$.
\begin{itemize}
\item If there is another $h^+$-factor ${\widetilde h}$ attached at $q$ then it has a primary $K$-factor $K'$. In this case, we glue the two $h^+$-factors, move to the other end of the moduli space and interpret the new factorization as an $\lt$-disk exactly as described above.
\item If there is an $f^a$-factor attached at $q$, with last $k$-connection $q'$ such that a $k$-factor $k'$ is attached at $q'$, then $k'$ must be the primary $k$-factor of some $h^-$-factor ${\widetilde h}$. Let the charge flow through the $f^a$-factor and $k'$. Charge ${\widetilde h}$ and activate its negative end.
\item If there is an $f^a$-factor attached at $q$, with last $k$-connection $q'$ such that a $h^0$-factor ${\widetilde h}$ is attached at $q'$, then ${\widetilde h}$ must have a primary $k$-factor. Let the charge flow through the $f^a$-factor. Charge ${\widetilde h}$ and let its positive end be glue active.
\end{itemize}
\item
${\overline h}$ is of type $h^0$ with a primary $k$-factor $k'$ and its positive end is glue-active. The positive puncture of ${\overline h}$ is attached to an $f^a$-factor. Gluing ${\overline h}$ to this $f^a$-factor gives a broken disk which is the boundary of a $1$-dimensional moduli space of holomorphic disks in $(X^a,L^a)$. Moving to the other end of this moduli space we find a similar breaking into one new $h^{a\pm}$- or $h^0$-factor ${\widetilde h}$ and new $f^a$-factors. If ${\widetilde h}$ is of type $h^{a+}$ then there may be several new $f^a$-factors but if it is of type $h^{a-}$ or $h^0$ then there is only one new $f^a$-factor. Charge ${\widetilde h}$.
\begin{itemize}
\item If ${\widetilde h}$ is of type $h^{a+}$ ($h^{a-}$) then activate its positive (negative) end.
\item If ${\widetilde h}$ is of type $h^0$ and is attached to $k'$ then regard $k'$ as its primary $k$-factor and let the positive end of ${\widetilde h}$ be flow active.
\item If ${\widetilde h}$ is of type $h^0$ and is not attached to $k'$, then let the sub-disk obtained by cutting at the positive puncture of the new $f^a$-factor which contains the new $f^a$-factor and $k'$ be the $K$-factor of the $h^+$-factor to which the positive puncture of $f^a$ was attached (see $(\rm c)$ above), and let the negative end of ${\widetilde h}$ be active.
\end{itemize}
\item
${\overline h}$ is of type $h^0$ and its positive end is flow-active. Let the charge flow through the $f^a$-factor attached at the positive puncture of ${\overline h}$, into the $h^+$-factor without primary $K$-factor where the positive puncture of the $f^a$-factor is attached. Charge this $h^+$-factor and activate its positive end.
\item
${\overline h}$ is of type $h^-$ and its positive end is active.
\begin{itemize}
\item If ${\overline h}$ is of type $h^{a-}$ then its positive puncture is attached to an $f^a$-factor. Gluing ${\overline h}$ to this $f^a$-factor we obtain a broken disk which is the boundary of a $1$-dimensional moduli space in $(X^a,L^a)$. Moving to the other end we find a similar broken disk with one new $h^{a\pm}$- or $h^0$-factor. Charge this new $h^\ast$-factor and let its positive end be active if it is of type $h^{a+}$ and let its negative end be active if it is of type $h^{a-}$ or $h^0$.
\item If ${\overline h}$ is of type $h^{b-}$ and its positive puncture is attached to an $f^b_0$-factor, then gluing ${\overline h}$ to this $f^b_0$-factor we obtain a broken disk which is the boundary of a $1$-dimensional moduli space in $(X^b_0,L^b_0)$. Moving to the other end we find a similar broken disk with one new $h^{b-}$- or $h^0$-factor ${\widetilde h}$. Charge ${\widetilde h}$ and let its positive end be active if it is of type $h^{0}$ and let its negative end be active if it is of type $h^{b-}$.
\item If ${\overline h}$ is of type $h^{b-}$ and its positive puncture is attached to a $k$-factor $k'$ which is the primary $k$-factor of some $h^0$- or $h^{b-}$-factor ${\widetilde h}$, then let $k'$ be the new primary $k$-factor of ${\overline h}$. Charge ${\widetilde h}$ and let the end of it where $k'$ is not attached be active.
\item If ${\overline h}$ is of type $h^{b-}$ and its positive puncture is attached to a $k$-factor $k'$ with positive puncture connected to an $f^a$-factor which in turn is connected to an $h^+$-factor ${\widetilde h}$ then let $D$ denote the sub-disk obtained by removing ${\overline h}$ from the $\lt$-disk and let $D'$ be the sub-disk of $D$ obtained by cutting at the positive puncture of $f^a$. We claim that $D'$ is the primary $K$-factor of ${\widetilde h}$. To see this we argue as follows. The disk $k'$ is not the primary $k$-factor of any $h^\ast$-factor and every $h^\ast$-factor of $D$ has a primary $k$- or $K$-factor. This implies on the one hand that the dimension of $D'$ is $-1$ and on the other that every sub-disk of $D$ cut off by cutting at a negative puncture of ${\widetilde h}$, except for the primary $K$-factor, has dimension $0$. The claim follows. Let $k'$ be the new primary $k$-factor of ${\overline h}$, let the charge flow through $k'$ and the $f^a$-factor. Charge ${\widetilde h}$ and let its positive end be active.
\item If ${\overline h}$ is of type $h^{b-}$ and is attached at another disk of type $h^{b-}$ with a primary $k$-factor $k'$, then gluing these two $h^{b-}$-disks, we obtain a broken disk which is the boundary of an $1$-dimensional reduced moduli space of holomorphic disks in $(Y^-\times\R,\Lambda^-\times\R)$. Moving to the other end of the reduced moduli space we obtain a broken disk with two new $h^{b-}$-factors. Charge the one of them which is not connected to $k'$ and let its negative end be active.
\end{itemize}
\end{itemize}
In all cases the disk resulting from the propagation law is a new non-isolated $\lt$-disk. Using the propagation law starting from any boundary $\lt$-disk, we can continue in a unique way until we reach another boundary $\lt$-disk in a finite number of steps. (For an example, see Figure \ref{fig:chflow}.) Any non-isolated $\lt$-disk contributing to \eqref{e:chhomT1'} -- \eqref{e:chhomT4} is a boundary $\lt$-disk and an argument similar to that in the proof of Lemma \ref{l:chiso} shows that any boundary $\lt$-disk contributes to \eqref{e:chhomT1'} -- \eqref{e:chhomT4}. Equation \eqref{e:chhomp} follows.
\end{pf}

Let $(X^b,L^b)$ be a stationary cobordism joined to a moving one $(X^a_s,L^a_s)$, $0\le s\le 1$, from below along $(Y^0,\Lambda^0)$ forming a new $1$-parameter family $(X^{ba}_s,L^{ba}_s)$. With notation as in Subsection \ref{s:chiso}, we have
\begin{align}\label{e:F_1expB}
F_1&=\Psi(F_0)=F_0+\Delta(F_0,K)\\\label{e:F_1expB'}
&=\left(f^b\I f^a_0+\Delta(f^a_0,k) \right),
\end{align}
where $K\in\V^+(X^{ba}_s,L^{ba}_s)$ is given by
\begin{equation}\label{e:KexpB}
K=\bigl\{f^b\leftarrow k\bigr\}\!\left(f^b, f^a_0,f^a_1\right)=
\bigl\{f^b\leftarrow k\bigr\}\!\left(f^b, f^a_0, f^a_0+\Delta(f^a_0,k)\right).
\end{equation}
Subsections  \ref{s:cob->chmap} and \ref{s:chiso} give the following diagram of chain maps
$$
\begin{CD}
\V(X^b,L^b) @>{\left[\pa_{f^b}\!\left.\I f^a_0\right)\right]}>> \V(X^{ba}_0,L^{ba}_0)\\
@V{\id}VV  @VV{\left[\pa_{F_0}\Psi\right]}V\\
\V(X^a,L^a) @>{\left[\pa_{f^b}\!\left.\I f^a_1\right)\right]}>> \V(X^{ba}_1,L^{ba}_1)
\end{CD}.
$$
We show below that the right-down and down-left chain maps in this diagram induce the same morphism of spectral sequences. We build chain homotopies using a map $\Omega\colon\V(X^b,L^b)\to\V(X^{ba}_1,L^{ba}_1)$ defined as follows. For  $v\in\V(X^b,L^b)$ let
$$
\Omega(v)=
\bigl\{v\leftarrow k\bigr\}\!
\bigl(f^b,f^a_0,f^a_0+\Delta(f^a_0,k)\bigr)\in\V(X^{ba}_1,L^{ba}_1).
$$

\begin{lma}\label{l:chhom'}
The $(+)$-action non-decreasing, filtration preserving chain maps $\left[\pa_{F_0}\Psi\right]\circ\left[\pa_{f^b}\!\left.\I f^a_0\right)\right]$ and $\left[\pa_{f^b}\!\left.\I f^a_1\right)\right]$ induce the same morphism of spectral sequences. In other words the following diagram commutes
\begin{equation}\label{e:chhom'}
\begin{CD}
\bigl\{E_{r\,;[\alpha]}^{p,q}(X^b,L^b)\bigr\}
@>{\left[\pa_{f^b}\!\left.\I f^a_0\right)\right]_\ast}>> \bigl\{E_{r\,;[\alpha]}^{p,q}(X^{ba}_0,L^{ba}_0)\bigr\}\\
@V{\id}VV  @VV{\left[\pa_{F_0}\Phi\right]_\ast}V\\
\bigl\{E_{r\,;[\alpha]}^{p,q}(X^b,L^b)\bigr\}
@>{\left[\pa_{f^b}\!\left.\I f^a_1\right)\right]_\ast}>> \bigl\{E_{r\,;[\alpha]}^{p,q}(X^{ba}_1,L^{ba}_1)\bigr\}
\end{CD},
\end{equation}
for every $\alpha>0$.
\end{lma}
\begin{pf}
The proof is analogous to the proof of Lemma \ref{l:chhom}. We will not repeat all arguments given there and we will use the same notation for Hamiltonian vectors ($h^{a\pm}$, $h^{b\pm}$, $h^0$, and $h^{\pm}=h^{a\pm}+h^{b\pm}$) and for Hamiltonian operators ($H$ for the one on $\V(X^{ba}_1,L^{ba}_1)$ and $h$ for the one on $\V(X^b,L^b)$). The map $\Omega$ respects the filtration and does not decrease $(+)$-action. Thus, with notation as in the proof of Lemma \ref{l:chhom}, $\Omega$ induces maps
$$
\Omega^p\colon \CC^p(X^b,L^b)\to\CC^p(X^{ba}_1,L^{ba}_1).
$$
for all $p$ and it suffices to show that
\begin{equation}\label{e:chhomp'}
\left[\pa_{F_0}\Psi\right]^p\circ\left[\pa_{f^b}\!\left.\I f^a_0\right)\right]^p(v)+ \left[\pa_{f^b}\!\left.\I f^a_1\right)\right]^p(v)=
d^{F_1}\circ \Omega^p(v) + \Omega^p\circ d^{f^b}(v),
\end{equation}
in $\CC^p(X^{ba}_1,L^{ba}_1)$, where $v$ is a formal disk with $p$ positive punctures.

As in the proof of Lemma \ref{l:chhom} we find that the terms contributing to the left hand side of \eqref{e:chhomp'} come from formal disks with one $v$-factor and at least one $k$-factor in the expressions
\begin{equation}\label{e:chhom'T1}
\Delta\Bigl(
\left(f^b+\epsilon v\I f^a_0\right),
\bigl\{f^b\leftarrow k\bigr\}\!(f^b,f^a_0, f^a_1)
\Bigr)
\end{equation}
and
\begin{equation}\label{e:chhom'T2}
\left(f^b+\epsilon v\I f^a_0+\Delta(f^a_0,k)\right).
\end{equation}
Since we work in $\CC^p(X^{ba}_1,L^{ba}_1)$ disks with more than $p$ positive punctures can be disregarded. It follows in particular that in any contributing disk, all factors except the $v$-factor has only one positive puncture. Therefore, \eqref{e:chhom'T1} reduces to
\begin{equation}\label{e:chhom'T1'}
\Delta\Bigl(
\left(f^b+\epsilon v\I f^a_0\right),
\bigl\{f^b\leftarrow k\bigr\}\!(f^b,0,0)\Bigr).
\end{equation}
The terms contributing to the right hand side of \eqref{e:chhomp'} come from formal disks with one $v$-factor, with at least one $k$-factor, and with all factors except the $v$-factor having only one positive puncture, in the expressions
\begin{equation}\label{e:chhom'T3}
H\bigl(F_1+\epsilon\,\Omega(v)\bigr)
\end{equation}
and
\begin{equation}\label{e:chhom'T4}
\Omega\bigl(h(f^b+\epsilon v)\bigr).
\end{equation}

We define an {\em $\lt$-disk} to be a formal disk with factors as follows.
\begin{itemize}
\item One $v$-factor, $n>0$ $h^\ast$-factors, $n>0$ $k$-factors, and arbitrary numbers of $f^a_0$- and $f^b$-factors.
\item All factors except the $v$-factor has only one positive puncture.
\item At most one of the $h^\ast$-factors are of type $h^{b+}$ or $h^{b-}$.
\item All of the $n+1$ factors of type $h^\ast$ and $v$, except one, has a primary $k$-factor.
\item A sub-disk obtained by cutting at a negative puncture $q$ of a $h^{a+}$-factor with primary $k$-factor, or at a positive puncture $p$ of a $v$-factor with primary $k$-factor, where $q$ or $p$ lies {\em after} the primary $k$-factor, does not contain any $k$-factors.
\end{itemize}
Note that an $h^\ast$-factor of type $h^{b\pm}$ cannot have a primary $k$-factor. If $q$ is a puncture of some factor of an $\lt$-disk and if the sub-disk cut-off by cutting at $q$ contains some $k$-factor then we say that $q$ is a {\em $k$-connection}.

We will distinguish one $h^\ast$- or $v$-factor in each $\lt$-disk as charged according to the following.
\begin{itemize}
\item[$({\rm a})$] A $v$-factor with all positive punctures connected to $f^a_0$-factors may be charged.
\item[$({\rm b})$] An $h^\ast$-factor without primary $k$-factor may be charged.
\item[$({\rm c})$] A $h^0$-factor with primary $k$-factor may be charged if one of its negative punctures is connected to $v$.
\end{itemize}

Consider the charged factor ${\overline h}$ of an $\lt$-disk of type $({\rm b})$ or $({\rm c})$. Let $p$ denote the positive puncture of ${\overline h}$ and let $q$ denote the last negative puncture of ${\overline h}$ which is a $k$-connection. We use the following terminology.
\begin{itemize}
\item The positive end of ${\overline h}$ is {\em obstructed} if one of the following holds.
\begin{itemize}
\item There is no disk attached to ${\overline h}$ at $p$,
\item ${\overline h}$ is attached to $v$ at $p$,
\item ${\overline h}$ is attached to a $k$-factor at $p$ and the positive puncture of this $k$-factor is not attached to any disk, or
\item ${\overline h}$ is attached to the primary $k$-factor of the $v$-factor at $p$.
\end{itemize}
\item The negative end of ${\overline h}$ is {\em obstructed} if one of the following holds.
\begin{itemize}
\item There is no disk attached to ${\overline h}$ at any of its negative punctures,
\item $v$ is attached to ${\overline h}$ at one of its negative punctures, or
\item the $k$-factor attached to ${\overline h}$ at $q$ is the primary $k$-factor of the $v$-factor.
\end{itemize}
\end{itemize}
An end which is not obstructed is called {\em un-obstructed}. An $\lt$-disk with charged factor of type $({\rm b})$ or $({\rm c})$ with two, one, or zero ends obstructed is called {\em isolated}, {\em boundary}, or {\em interior}, respectively. An $\lt$-disk with charged factor of type $({\rm a})$ is {\em isolated}.

We view the disks in \eqref{e:chhom'T1'} as $\lt$-disks by letting their $v$-factors be charged. They are all isolated. The disks in \eqref{e:chhom'T2} are of two kinds. In the first kind, all positive punctures of the $v$-factor are connected to $f^a_0$-factors. We view such disks as $\lt$-disks by letting their $v$-factor be charged. In the second kind some positive puncture of the $v$-factor is not connected to an $f^a_0$-factor and this puncture must be a $k$-connection. Let $q$ be the last positive puncture of the $v$-factor which is a $k$-connection. If there is an $h^0$-factor attached at $q$ then let that be the charged $h^\ast$-factor. If there is a $k$-factor $k'$ attached to $v$ at $q$ then we view $k'$ as the primary $k$-factor of $v$ and charge the $h^\ast$-factor in $\Delta(f_0^a,k)$ which had $k'$ as primary $k$-factor. We view the disks in \eqref{e:chhom'T3} and \eqref{e:chhom'T4} as $\lt$-disks by charging the $h^\ast$-factor which corresponds to the operator $H$ and to the operator $h$, respectively. We first show that isolated $\lt$-disks cancel. They are of the following types.
\begin{itemize}
\item[$({\rm I})$] Disks from \eqref{e:chhom'T1'} are of type $({\rm I})$.
\end{itemize}

Isolated disks from \eqref{e:chhom'T2} are of four types. Let $p$ denote the first positive $k$-connection of $v$.
\begin{itemize}
\item[$({\rm II}a)$] $f^a$-factors are attached at all positive punctures of the $v$-factor.
\item[$({\rm II}b)$] An $h^0$-factor is attached at $p$.
\item[$({\rm II}c)$] A $k$-factor, with positive puncture attached to a $h^{a+}$-factor with no disk attached at its positive puncture, is attached at $p$.
\item[$({\rm II}d)$] A $k$-factor, with a $h^{a-}$-factor attached at some negative puncture, is attached at $p$.
\end{itemize}

Isolated disks from \eqref{e:chhom'T3} are of three types. Let ${\overline h}$ denote the charged $h^\ast$-factor.
\begin{itemize}
\item[$({\rm III}a)$] ${\overline h}$ is of type $h^{a+}$ and is attached to the primary $k$-factor of $v$ in $\Omega(v)$.
\item[$({\rm III}b)$] ${\overline h}$ is of type $h^{a-}$ and is attached to the primary $k$-factor of $v$ in $\Omega(v)$.
\item[$({\rm III}c)$] ${\overline h}$ is of type $h^{b\pm}$ and is attached to the $v$-factor in $\Omega(v)$.
\end{itemize}

Isolated disks from \eqref{e:chhom'T4} are of two types. Let ${\overline h}$ denote the charged $h^\ast$-factor.
\begin{itemize}
\item[$({\rm IV}a)$] ${\overline h}$ is of type $h^0$ and the positive puncture of ${\overline h}$ is attached to the primary $k$-factor of the $h(f^b+\epsilon v)$-factor.
\item[$({\rm IV}b)$] ${\overline h}$ is of type $h^{b\pm}$.
\end{itemize}
It is straightforward to check that the isolated disks cancel according to the following.
$$
\begin{matrix}
({\rm I})\leftrightarrow({\rm II}a), &
({\rm II}b)\leftrightarrow({\rm IV}a), &
({\rm II}c)\leftrightarrow({\rm III}a), \\
({\rm II}d)\leftrightarrow({\rm III}b), &
({\rm III}c)\leftrightarrow({\rm IV}b).
\end{matrix}
$$

To show that non-isolated $\lt$-disks cancel we use the following propagation law for such disks. Let ${\overline h}$ denote the charged $h^\ast$-factor.
\begin{itemize}
\item The positive end of ${\overline h}$ with positive puncture $p$ is active.
\begin{itemize}
\item If $p$ is attached to an $f^a_0$-factor (an $f^b$-factor) then glue ${\overline h}$ and the $f^a_0$-disk ($f^b$-disk) to get a broken disk which is the boundary of a $1$-dimensional moduli space of holomorphic disks in $(X^a_0,L^a_0)$ (in $(X^b,L^b)$). Moving to the other end of the moduli space we find a similar breaking into one $h^\ast$-factor and $f^a_0$-factors ($f^b$-factors). Charge the new $h^\ast$-factor and activate the end of it where no new $f^a_0$-factor ($f^b$-factor) is attached.
\item If $p$ is attached to another $h^\ast$-factor then this other $h^\ast$-factor has a primary $k$-factor
$k'$. Glue the two $h^\ast$-factors to obtain a broken disk which is one end of a $1$-dimensional reduced moduli space in the symplectization of a $(\pm\infty)$-boundary of $(X^a,L^a)$ or $(X^b_s,L^b_s)$. Moving to the other end of the reduced moduli space we find a similar splitting. Charge the new $h^\ast$-factor not connected to $k'$ and let the end of it where the new uncharged $h^\ast$-factor is not attached be active.
\item If $p$ is attached to a $k$-factor $k'$ then let the charge flow through $k'$ into the $h^\ast$-factor ${\widetilde h}$ which has $k'$ as primary $k$-factor. Let $k'$ be the primary $k$-factor of ${\overline h}$ and let ${\widetilde h}$ be the new charged factor with the end not connected to $k'$ active.
\end{itemize}
\item The negative end of ${\overline h}$, with last $k$-connection $q$, is active.
\begin{itemize}
\item
If all negative punctures of ${\overline h}$ are attached to $f^a_0$-factors ($f^b$-factors) then glue ${\overline h}$ and the $f^a_0$-disks ($f^b$-disks) to get a broken disk which is the boundary of a $1$-dimensional moduli space of holomorphic disks in $(X^a_0,L^a_0)$ (in $(X^b,L^b)$). Moving to the other end of the moduli space we find a similar breaking into one $h^\ast$-factor and $f^a_0$-factors ($f^b$-factors). Charge the new $h^\ast$-factor and activate the end of it where no new $f^a_0$-factor ($f^b$-factor) is attached.
\item
If some negative puncture of ${\overline h}$ is a $k$-connection and if there is an $h^\ast$-factor ${\widetilde h}$ attached at $q$ then ${\widetilde h}$ has a primary $k$-factor $k'$. Glue ${\overline h}$ and ${\widetilde h}$ to get a broken disk which is the boundary of a $1$-dimensional reduced moduli space of holomorphic disks in the symplectization of a $(\pm\infty)$-boundary of $(X^a,L^a)$ or $(X^b_s,L^b_s)$. Moving to the other end of the reduced moduli space we find a similar splitting. Charge the new $h^\ast$-factor not connected to $k'$ and let the end of it where the new uncharged $h^\ast$-factor is not attached be active.
\item
If some negative puncture of ${\overline h}$ is a $k$-connection and if there is a $k$-factor $k'$ attached at $q$ then $k'$ is the primary $k$-factor of some $h^\ast$-factor ${\widetilde h}$. Let the charge flow through $k'$. Charge ${\widetilde h}$ and let $k'$ be the primary $k$-factor of ${\overline h}$.
\end{itemize}
\end{itemize}
In all cases the disk resulting from the propagation law is a new non-isolated $\lt$-disk. Using the propagation law starting from any boundary $\lt$-disk, we can continue in a unique way until we reach another boundary $\lt$-disk in a finite number of steps. Any non-isolated $\lt$-disk contributing to \eqref{e:chhom'T2} -- \eqref{e:chhom'T4} is a boundary $\lt$-disk and an argument similar to that in the proof of Lemma \ref{l:chiso} shows that any boundary $\lt$-disk contributes \eqref{e:chhom'T2} -- \eqref{e:chhom'T4}. Equation \eqref{e:chhomp'} follows.
\end{pf}

\section{Proofs}\label{S:proofs}

\begin{pf}[Proof of Theorem \ref{t:main1}]
If $(X_s,L_s)$, $0\le s\le 1$, is a $1$-parameter family of exact cobordisms fixed outside a compact set then it follows from Subsection \ref{s:mixpureini} that after small perturbation the differential is independent of $s$ except when values of $s$ for which there are $(-1)$-disks are passed. It follows from Lemmas \ref{l:chiso} and \ref{l:chiso'} that the spectral sequence remains unchanged at such instances.
\end{pf}

\begin{pf}[Proof of Theorem \ref{t:main2}]
Let $Y$ be an odd-dimensional manifold. Let $Y_0$ and $Y_1$ denote $Y$ equipped with two contact forms which determine isotopic contact structures. Let furthermore  $\Lambda_0\subset Y_0$ and $\Lambda_1\subset Y_1$ be Legendrian submanifolds and assume that the pairs $(Y_0,\Lambda_0)$ and $(Y_1,\Lambda_1)$ are contact isotopic through $(Y_s,\Lambda_s)$, $0\le s\le 1$.

Using Lemmas \ref{l:isotopy} and \ref{l:formch}, we see that such an isotopy determines  an exact cobordism $(X_{01},L_{01})$ with $(+\infty)$-boundary $(Y_1,\Lambda_1)$ and $(-\infty)$-boundary $(Y_0,\Lambda_0)$. Using the inverse isotopy we get a cobordism $(X_{10},L_{10})$ with $(+\infty)$-boundary $(Y_0,\Lambda_0)$ and $(-\infty)$-boundary $(Y_1,\Lambda_1)$. Joining the negative end of $(X_{01},L_{01})$ to the positive end of $(X_{10},L_{10})$ we obtain a cobordism $(X_{11},L_{11})$ which admits a compactly supported isotopy deforming it to $(Y_1\times\R,\Lambda_1\times\R)$. Let $f^j$, $j=0,1$, denote the potential vectors in $\V^+(Y_j\times\R,\Lambda_j\times\R)$ and let $F^{ij}$, $i,j\in\{0,1\}$ denote the potential vectors in $\V^+(X_{ij},L_{ij})$. Then there are grading respecting chain maps as follows
$$
\begin{CD}
\V_{[\alpha]}(Y_1\times\R,\Lambda_1\times\R)
@>{\left[\pa_{f^1}\!\left(F^{01}\I\right.\right]}>>
\V_{[\alpha]}(X_{01},L_{01})\\
{} @V{\id}VV  @VV{\left[\pa_{F^{01}}\!\left(F^{11}\I\right.\right]}V\\
\V_{[\alpha]}(Y_1\times\R,\Lambda_1\times\R)
@<{\Phi}<<
\V_{[\alpha]}(X_{11},L_{11}),
\end{CD}
$$
for any $\alpha>0$, where $\Phi$ is the chain map from Theorem \ref{t:main1}, see also Lemma \ref{l:chiso}, induced by the deformation of $(X_{11},L_{11})$. Lemma \ref{l:compcoblin} implies that the composition of the first two maps equals ${\left[\pa_{f^1}\!\left(F^{11}\I\right.\right]}$ which, after composition with $\Phi$, by Lemma \ref{l:chhom}, is chain homotopic to $\left[\pa_{f^1}\!\left(f^1\I\right.\right]=\id$. It follows that the composition
$$
\Phi\circ\left[\pa_{F^{01}}\!\left(F^{11}\I\right.\right]\circ
\left[\pa_{f^1}\!\left(F^{01}\I\right.\right]
$$
induces an isomorphism of spectral sequences. Theorem \ref{t:main1} implies that $\Phi$ induces an isomorphism. Thus, the map induced by $\left[\pa_{f^1}\!\left(F^{01}\I\right.\right]$ is a monomorphism and that induced by
$\left[\pa_{F^{01}}\!\left(F^{11}\I\right.\right]$ is an epimorphism.

Gluing $(X_{11},L_{11})$ at its negative end to the positive end of $X_{01}$ we get a cobordism $(X_{011},L_{011})$, with potential vector $F^{011}\in\V^+(X_{011},L_{011})$, which can be deformed to $(X_{01},L_{01})$. Arguing as above using the diagram
$$
\begin{CD}
\V_{[\alpha]}(X_{01},L_{01})
@>{\left[\pa_{F^{01}}\!\left(F^{11}\I\right.\right]}>>
\V_{[\alpha]}(X_{11},L_{11})\\
{} @V{\id}VV  @VV{\left[\pa_{F^{11}}\!\left(F^{011}\I\right.\right]}V\\
\V_{[\alpha]}(X_{01},L_{01})
@<{\Phi}<<
\V_{[\alpha]}(X_{011},L_{011}),
\end{CD}
$$
we find that the map induced by $\left[\pa_{F^{01}}\!\left(F^{11}\I\right.\right]$ is also a monomorphism and it follows that
\begin{equation}\label{e:pf1}
\left[\pa_{f^1}\!\left(F^{01}\I\right.\right]\colon
\V_{[\alpha]}(Y_1\times\R,\Lambda_1\times\R)\to
\V_{[\alpha]}(X_{01},\Lambda_{01})
\end{equation}
induces an isomorphism of spectral sequences. Taking inverse limits, we find
$$
\bigl\{E^{p,q}_r(Y_1,\Lambda_1)\bigr\}=
\lim_{\longleftarrow\,\alpha}\bigl\{E^{p,q}_{r\,;[\alpha]}(Y_1\times\R,\Lambda_1\times\R)\bigr\}
=\lim_{\longleftarrow\,\alpha}\bigl\{E^{p,q}_{r\,;[\alpha]}(X_{01},L_{01})\bigr\}.
$$

A similar argument, which uses Lemmas \ref{l:chiso'} and \ref{l:chhom'} instead of Lemmas \ref{l:chiso} and \ref{l:chhom}, shows that
\begin{equation}\label{e:pf2}
\left[\pa_{f^0}\!\left.\I F^{01}\right)\right]\colon
\V_{[\beta]}(Y_0\times\R,\Lambda_0\times\R)\to
\V_{[\beta]}(X_{01},\Lambda_{01}).
\end{equation}
induces an isomorphism of spectral sequences. Taking inverse limits, we find that
$$
\bigl\{E^{p,q}_r(Y_0,\Lambda_0)\bigr\}=
\lim_{\longleftarrow\,\beta}\bigl\{E^{p,q}_{r\,;[\beta]}(Y_0\times\R,\Lambda_0\times\R)\bigr\}
=\lim_{\longleftarrow\,\beta}\bigl\{E^{p,q}_{r\,;[\beta]}(X_{01},L_{01})\bigr\}.
$$
as well. In the case that $\Lambda_0$ and $\Lambda_1$ has only finitely many Reeb chords, all inverse systems above stabilizes and it suffices to take $\alpha$ in \eqref{e:pf1} and $\beta$ in \eqref{e:pf2} sufficiently large to establish the required isomorphism of spectral sequences. This proves the theorem.
\end{pf}

\appendix

\section{}\label{A:A}
In this section we describe two standard contact geometric constructions. Let $Y$ be a contact manifold with contact form $\lambda$ and let $\Lambda\subset Y$ be a Legendrian submanifold. Let $\Lambda_t\subset Y$, $0\le t\le M$, $\Lambda_0=\Lambda$, be a Legendrian isotopy. Using a standard argument, see e.g. \cite{G} Theorem  2.41, we find a time dependent function (a contact Hamiltonian) $H_t\colon Y\to\R$ such that the flow $\phi_t\colon Y\to Y$ of the time dependent vector field $X_t$, uniquely specified by the properties $\lambda(X_t)=H_t$ and $\iota_{X_t} d\lambda = dH_t(R_\lambda)\lambda-dH_t$, is a $1$-parameter family of contact-diffeomorphisms such that $\phi_t(\Lambda)=\Lambda_t$.
\begin{lma}\label{l:isotopy}
Let $\tilde H_t\colon V\to \R$ be a time dependent contact Hamiltonian with corresponding flow $\tilde\phi_t\colon Y\to Y$, $0\le t\le 1$. If $M>0$ and if $H_t=M^{-1}\tilde H_t$ is the re-scaled contact Hamiltonian with corresponding flow $\phi_t$, $0\le t\le M$ then $\phi_{M^{-1}t}(x)=\tilde\phi_t(x)$. Moreover, if $M>0$ is sufficiently large then the form
\begin{equation}\label{e:sfisotopy}
d\left(e^t(\lambda - dH_t)\right)
\end{equation}
is a symplectic form on $V\times\R$ for which the map $\Lambda\times\R\to Y\times\R$
\begin{equation}\label{e:sfisotopy'}
(q,t)\mapsto(\phi_t(q),t)
\end{equation}
is an exact Lagrangian embedding which agrees with $\Lambda_0\times \R$ and $\Lambda_1\times\R$ in neighborhoods of $t=-\infty$ and $t=+\infty$, respectively. \end{lma}

\begin{pf}
We may assume that the isotopy is constant near the end points of the interval. In particular the contact Hamiltonian vanishes there. The statement on the properties of the flows follows by a change of variables in the corresponding differential equations. By taking $M>0$ sufficiently large, we may make the perturbation $dH_t$ of the symplectic form $d(e^t\lambda)$ in \eqref{e:sfisotopy} arbitrarily small. Since the perturbed form is obviously closed and since the symplectic non-degeneracy condition is open, it follows that also the perturbed form is symplectic. The pull-back of $e^t(\lambda+dH_t)$ under \eqref{e:sfisotopy'} is
\begin{align*}
e^t\left(\left(\lambda(X_t)+\frac{\pa H_t}{\pa t}\right)dt + \frac{\pa H_t}{\pa q}dq\right)=
e^t\left(\left(H_t+\frac{\pa H_t}{\pa t}\right)dt + \frac{\pa H_t}{\pa q}dq\right)
=d (e^t H_t).
\end{align*}
This shows exactness and finishes the proof.
\end{pf}

We consider secondly a situation where we change the contact form. We have $\tilde \lambda_t = e^{\tilde g_t(y)}\lambda$, $0\le t\le 1$, for some $1$-parameter family of functions $\tilde g_t\colon Y\to Y$.

\begin{lma}\label{l:formch}
If $M>0$ is sufficiently large and $\lambda_t=e^{g_t(y)}\lambda$, where $g_t(y)=M^{-1}\tilde g_t(y)$ then the form
\begin{equation}\label{e:form1}
d\left(e^t(\lambda_t)\right)
\end{equation}
is a symplectic form on $Y\times\R$ for which the map $\Lambda\times\R\to Y\times\R$
\begin{equation}\label{e:form2}
(q,t)\mapsto(q,t-g_t(q))
\end{equation}
is an exact Lagrangian embedding which agrees with $\Lambda_0\times \R$ and $\Lambda_1\times\R$ in neighborhoods of $t=-\infty$ and $t=+\infty$, respectively.
\end{lma}

\begin{pf}
Take $M>0$ sufficiently large to ensure that the perturbation of the symplectic form in \eqref{e:form1} is sufficiently small. The pull-back of the form $e^t\lambda_t$ under \eqref{e:form2} equals $e^t\lambda$, which is identically equal to $0$ on $\Lambda\times\R$.
\end{pf}

\section{}\label{A:B}
In this appendix we give a brief description of the properties of holomorphic disks that are being used in the paper. Details on the subject can be found in \cite{BEHWZ}, \cite{H,HWZ}, and \cite{EES1,EES4}.

\subsection{Compactness properties of moduli spaces of holomorphic disks}\label{s:compactness}
Let $(X,L)$ be an exact symplectic cobordism with symplectic form $\omega=d\beta$, with $(\pm\infty)$-boundary $(Y^\pm,\Lambda^\pm)$ where the contact form is $\lambda^\pm$ and where $\beta=e^t\lambda^\pm$, $t\in\R_\pm$, in the corresponding end $(Y^\pm\times\R_\pm,\Lambda^\pm\times\R_\pm)$. Fix an adjusted almost complex structure $J$ in $(X,L)$. That is, $J$ is required to be compatible with the symplectic form, and in the ends, $J$ is required to be invariant under $\R$-translations and to pair the $\R$-direction with the direction of the Reeb field. A map $u\colon S\to X$, where $S$ is a Riemann surface with complex structure $j$ is called {\em $J$-holomorphic} if it satisfies the differential equation $\bar\pa_J u = 0$, where
$$
\bar\pa_J u = du+J\circ du\circ j.
$$
We say that a Riemann surface $S$ is {\em simple} if it is either a sphere or a disk. We will consider $J$-holomorphic maps of Riemann surfaces into $X$ with boundary on $L$. Let $S$ denote a Riemann surface with some number of punctures on the boundary and some number of interior punctures. We write $(X,L)=(\bar X,\bar L)\cup (Y^+\times\R_+,\Lambda\times\R_+)\cup(Y^-\times\R_-,\Lambda^-\times\R_-)$, where $(\bar X,\bar L)$ is the compact part of the cobordism. Define the {\em $\omega$-energy} of a map $u\colon (S,\pa S)\to(X,L)$ as
$$
E_\omega(u)=\int_{u^{-1}(\bar X)} u^\ast\omega +
\int_{u^{-1}(Y^-\times\R_-)}u^\ast d\lambda^- +
\int_{u^{-1}(Y^+\times\R_+)}u^\ast d\lambda^+,
$$
and the {\em $\lambda$-energy} as, writing $u=(v,t)$, where $t\colon S\to\R_\pm$ and $v\colon S\to Y^\pm$ in the ends,
$$
E_\lambda(u)=\sup_{\phi_-}\left(\int_{u^{-1}(Y^-\times\R_-)}\phi_-(t)\, dt\wedge (v^\ast\lambda^-)\right)+
\sup_{\phi_+}\left(\int_{u^{-1}(Y^+\times\R_+)}\phi_+(t)\, dt\wedge (v^\ast\lambda^+)\right),
$$
where $\phi_+$ and $\phi_-$ range over all smooth functions with
$$
\int_{-\infty}^{0}\phi_-(t)\,dt=1\quad\text{and}\quad
\int_0^{\infty}\phi_+(t)\,dt=1,
$$
respectively.

Since the almost complex structure $J$ is adjusted, it follows that if $c$ is a Reeb orbit or a Reeb chord then $c\times\R_\pm\subset Y^\pm\times\R_\pm$ is a $J$-holomorphic surface in the end. We call these a {\em Reeb orbit cylinder} and a {\em Reeb chord strip} according to whether $c$ is an orbit or a chord. As in \cite{BEHWZ}, Proposition 6.2, one shows the following
\begin{lma}\label{l:asympt}
Let $u\colon(S,\pa S)\to (X,L)$ be a $J$-holomorphic surface of finite energy and without removable singularities. Assume that the contact forms at the $(\pm\infty)$-boundary of $(X,L)$ are generic so that Reeb chords and Reeb orbits are isolated. Then, near each boundary puncture, $u$ is (exponentially) asymptotic to a Reeb chord strip and, near each interior puncture, $u$ is (exponentially) asymptotic to a Reeb orbit cylinder.\qed
\end{lma}

\begin{rmk}\label{r:hol->formal}
It follows from Lemma \ref{l:asympt} that any $J$-holomorphic disk in an exact cobordism without interior punctures determines a formal disk.
\end{rmk}

A puncture of $u\colon (S,\pa S)\to(X,L)$ is positive (negative) if $u$ is asymptotic to a Reeb chord or Reeb orbit in the $(+\infty)$-boundary (in the $(-\infty)$-boundary) of $(X,L)$ near this puncture. The various energy concepts are easily controlled in the case of an exact cobordism. Recall the notion of the action of a curve $\gamma$ in $Y^\pm$, $\A(\gamma)=\int_\gamma\lambda^\pm$.
\begin{lma}\label{l:action}
For any finite energy $J$-holomorphic surface $u\colon(S,\pa S)\to (X,L)$ with positive punctures at Reeb chords or orbits $c_1^+,\dots,c_p^+$, and with negative punctures at Reeb chords or orbits $c_1^-,\dots,c_q^-$, the following holds.
\begin{align*}
0\le E_\omega(u)&=\sum_{j=1}^p\A(c^+_j)-\sum_{k=1}^q\A(c^-_k),\\
E_\lambda(u)&\le 2\sum_{j=1}^p\A(c^+_j).
\end{align*}
\end{lma}
\begin{pf}
That $E_\omega(u)\ge 0$ follows since $\omega$ is positive on $J$-complex lines. Assume that the levels $0$ in the ends are regular levels of the function $t\circ u$. (If they are not, the argument below together with a straightforward approximation argument using nearby regular levels give the desired result.) Since $\omega=d\beta$ is exact and since $\beta|L=df$ for a function $f\colon L\to\R$ which is constant in the ends we find
\begin{equation}\label{e:in}
\int_{u^{-1}(\bar X)} u^\ast\omega=\int_{u^{-1}(\pa\bar X\cup \bar L)}u^\ast \beta=\pm\int_{u^{-1}(\pa\bar X)}u^\ast\lambda^\pm.
\end{equation}
Also,
\begin{equation}\label{e:out-}
\int_{u^{-1}(Y^-\times\R_-)}u^\ast d\lambda^-=\int_{u^{-1}(\pa\bar X)}u^\ast\lambda^- -\sum_{k=1}^q\A(c^-_k)
\end{equation}
and
\begin{equation}\label{e:out+}
\int_{u^{-1}(Y^+\times\R_+)}u^\ast d\lambda^+=-\int_{u^{-1}(\pa\bar X)}u^\ast\lambda^+ +\sum_{j=1}^p\A(c^+_j).
\end{equation}
Adding \eqref{e:in}, \eqref{e:out-}, and \eqref{e:out+} gives the result on $E_\omega(u)$.

Consider the $\lambda$-energy. We have
$$
\int_{u^{-1}(Y^+\times\R_+)}\phi(t)dt\wedge \left(v^\ast\lambda^+\right)
=\int_{u^{-1}(Y^+\times\R_+)} v^\ast d(\psi\lambda^+)
-\int_{u^{-1}(Y^+\times\R_+)}\psi \left(v^\ast d\lambda^+\right),
$$
where $\psi(s)=\int_0^s\phi(t)\,dt$. The last integral in the right hand side is positive and the first term equals $\sum_{j=1}^p\A(c^+_j)$. We also have
$$
\int_{u^{-1}(Y^-\times\R_-)}\phi(t)dt\wedge \left(v^\ast\lambda^-\right)=
\int_{u^{-1}(Y^-\times\R_-)} v^\ast d(\psi\lambda^-)
-\int_{u^{-1}(Y^-\times\R_-)}\psi\left(v^\ast d\lambda^-\right),
$$
where $\psi(s) = \int_{-\infty}^{s}\phi(t)\,dt$. Again the last integral is positive and the first term equals $\int_{u^{-1}(Y^-\times\{0\})}u^\ast \beta$. Using the $\omega$-positivity of $J$-holomorphic surfaces, we can estimate this first term, first by $\int_{u^{-1}(Y^+\times\{0\})}u^\ast \beta$ and then by $\sum_{j=1}^p\A(c^+_j)$. The lemma follows.
\end{pf}

Since the sum of $\omega$- and $\lambda$-energy of any non-constant $J$-holomorphic surface is positive, any such surface has at least one positive puncture. Furthermore, a bound on the action of the positive punctures of a $J$-holomorphic surface automatically gives an energy bound. Applying Subsection 11.3 of \cite{BEHWZ} we get the following result.
\begin{lma}\label{l:gencomp}
The moduli space $\M$ of $J$-holomorphic disks in an exact cobordism $(X,L)$ with positive punctures at Reeb chords $c_1,\dots,c_l$ has a natural compactification consisting of several level simple $J$-holomorphic surfaces joined at Reeb chords or at Reeb orbits.\qed
\end{lma}

\begin{rmk}\label{r:1-jetcmp}
The results in this section can be generalized to some cases when $X$ does not have compact ends. We discuss one such case important for applications to Legendrian isotopy problems. Let $M$ be a smooth manifold and let $J^1(M)$ be its $1$-jet space with a contact form $\lambda$ which agrees with its standard contact form outside a compact set, see e.g. \cite{E1}. Let $\Lambda^+$ and $\Lambda^-$ be closed Legendrian submanifolds of $J^1(M)$. Consider the symplectization $J^1(M)\times\R$ and its exact Lagrangian submanifolds $L^+=\Lambda^+\times[1,\infty)$ and $L^-=\Lambda^-\times(-\infty,-1]$. Let $K'$ be a compact subset of $J^1(M)$ containing $\Lambda^+\cup \Lambda^-$ and the region where $\lambda$ is not equal to the standard contact form. Let $K$ be a compact subset of $J^1(M)\times\R$ of the form $K=K'\times[-T,T]$, $T>1$. We say that $(J^1(M)\times\R, L^+\cup L^-) - K$ is a {\em closed end}.

We say that an exact cobordisms $(X,L)$ has {\em tame ends} if there is a compact subsect $C$ of $X$ such that the connected components of $(X,L)-C$ are either standard ends $(Y\times\R_{\pm},\Lambda\times\R_{\pm})$ for a compact contact manifold $Y$ with a Legendrian submanifold $\Lambda$, or closed ends. All the results in this section hold for exact cobordisms with tame ends. To see that note that by monotonicity of the area of $J$-holomorphic disks, an energy bound on the disk implies that the parts of the disk which map into a closed end $(J^1(M)\times\R,L^+\cup L^-)-K$ must have projection into $J^1(M)$ which stays in some compact set $K''\subset J^1(M)$ such that $\Lambda^+\cup\Lambda^-\subset K''$ and such that the region where the contact form is non-standard is contained in $K''$. (The size of $K''$ depends on the energy bound.)

We define a $(\pm\infty)$-boundary of an exact cobordism $(X,L)$ with tame ends to be either $(Y,\Lambda)$ for one of the usual ends $(Y\times\R_{\pm},\Lambda\times\R_{\pm})$ of $(X,L)$, or $(J^1(M),\Lambda^{\pm})$, for closed ends $\bigr(J^1(M)\times[T,\infty),\Lambda^+\times[T,\infty)\bigl)$ or $\bigl(J^1(M)\times(-\infty,T],\Lambda^-\times(-\infty,T]\bigr)$ of $(X,L)$.
\end{rmk}

\subsection{Good ends}\label{s:goodends}
With the results in Subsection \ref{s:compactness} established, we are ready to state the technical condition we impose on the ends of our exact cobordisms. The two main consequences of this technical condition are that it allows us to work with $\Z_2$-coefficients and that it guarantees that there is no input from the contact homology of the $(\pm\infty)$-boundary, see Remark \ref{r:CH}.

Let $(X,L)$ be an exact cobordism with tame ends. A contact form $\lambda^+$ on a contact manifold $Y^+$ in the $(+\infty)$-boundary of $X$ is {\em good} if the following condition holds. For any Reeb orbit $\gamma^+$ in $Y^+$ which is contractible in $X$, the formal dimension $\dim(\M(\gamma^+))$ of any moduli space $\M(\gamma^+)$ of one-punctured $J$-holomorphic spheres in $X$ with positive puncture at $\gamma^+$ satisfies $\dim(\M(\gamma^+))\ge 2$. A contact form $\lambda^-$ on a contact manifold $Y^-$ in the $(-\infty)$-boundary of $X$ is {\em good} if the following condition hold. For any Reeb orbit $\gamma^-$ in $Y^-$ which is contractible in $Y^-$, the formal dimension $\dim(\M(\gamma^-))$ of any moduli space $\M(\gamma^-)$ of one-punctured $J$-holomorphic spheres with positive puncture at $\gamma^-$ in the symplectization $Y^-\times\R$ satisfies $\dim(\M(\gamma^-))\ge 2$. If $Y$ is a contact manifold with a contact form $\lambda$ we say that $\lambda$ is {\em good} if the symplectization $Y\times\R$ has good ends.

We note that an end of an exact cobordism is good if the corresponding $(\pm\infty)$-boundary is a sphere of any dimension endowed with a contact form which is sufficiently close to the standard contact form, see \cite{B, EGH}. Also, an end is good if the corresponding $(\pm\infty)$-boundary is a $1$-jet space with its standard contact form: in this case there are no closed Reeb orbits at all. For general closed contact manifolds, monotonicity conditions which guarantee that the index contribution of a Reeb chord is (roughly) proportional to its action can be used to show that a contact form is good.

\begin{lma}\label{l:nicecomp}
Let $\M$ be a moduli space of admissible $J$-holomorphic disks of bounded energy and without interior punctures in an exact cobordism $(X,L)$, where $J$ is a generic almost complex structure. Assume that the formal dimension of $\M$ is $\le 1$ ($\le 2$ in trivial cobordisms) and that $(X,L)$ has good ends. Then the compactification of $\M$ consists of many level curves, where each level has only admissible disks without interior punctures. It follows in particular that the boundary of a moduli space of admissible disks of dimension $1$ consists of broken curves of two levels: admissible disks of dimension $0$ in $(X,L)$ and one admissible disk of dimension $1$ in the symplectization of a $(\pm\infty)$-boundary.
\end{lma}

\begin{pf}
Lemma \ref{l:gencomp} implies that $\M$ has a compactification with boundary consisting of broken curves. Define an {\em admissible disk with interior punctures} exactly as an admissible disk but allowing its domain to have interior punctures.
Note that no two positive boundary punctures of an admissible disk (with interior punctures) can be asymptotic to the same Reeb chord, see the proof of Lemma \ref{l:no2mixed}. Transversality for such disks can thus be achieved by perturbing near a positive boundary puncture as in \cite{EES1,EES3} and it follows that for generic $J$, any admissible $J$-holomorphic disk (with interior punctures) in $(X,L)$ (in $(Y^\pm\times\R,\Lambda^\pm\times\R)$) has dimension $\ge 0$ (has dimension $\ge 1$ unless it is a Reeb chord strip).

Consider a broken curve in the boundary of $\M$. Our assumption on good ends rules out the possibility that one piece of this broken curve is a sphere with one puncture: any such curve has dimension at least $2$ and other pieces have dimensions at least $0$ (at least $1$ in a symplectization).

Using the dimension bound, we find that remaining possibilities for breaking is either as a $1$-dimensional disk in some end and $0$-dimensional components in the cobordism joined at punctures,
or as a boundary bubbling, breaking into two disk families of dimensions $d_1$ and $d_2$, where $d_1+d_2+2=n$, $n=\dim(L)$. However, as the disk is admissible, interior bubbling is ruled out: one of the disks which results from the boundary bubbling would not have a positive puncture by the definition of admissibility. Thus, the later case does not appear. In the former case, it follows from Lemma \ref{l:nonadmsub} that all factors in the broken disk are admissible.
\end{pf}

\begin{rmk}\label{r:CH}
In order to extend the theory described in the current paper to exact cobordisms $(X,L)$ with tame ends which are not good, also disks with interior punctures must be taken into account. This requires multi-valued perturbations and $\Q$-coefficients. In contact homology, similar mixed theories are important in the description of the effect of surgery, see \cite{BEE}.
\end{rmk}

\subsection{Initial transversality for $1$-parameter families}\label{s:mixpureini}
Consider a $1$-parameter family of exact cobordisms $(X_s,L_s)$, $0\le s\le 1$, with $(\pm\infty)$-boundary $(Y^\pm,\Lambda^\pm)$ and with good ends, which is constant near its ends.
\begin{lma}\label{l:classtv}
For generic $1$-parameter families $J_s$ of adjusted almost complex structures the following holds. For all $s$, the moduli space of admissible $J_s$-holomorphic disks (with interior punctures) in $(X_s,L_s)$ of formal dimension $<-1$ ($<0$ if the cobordism is trivial)  is empty. There is a finite number of isolated instances $0<s_1<s_2<\dots<s_n<1$ where the moduli spaces of admissible $J_s$-holomorphic disks (with interior punctures) in $(X_s,L_s)$ of formal dimension $-1$ ($0$ if the cobordism is trivial) are non-empty. At each such instance the moduli space contains exactly one disk which is a transversely cut out $0$-dimensional component of the parameterized moduli space.

If $\bigl\{u_{s(j)}^j\bigr\}_{j=1}^\infty$ is a sequence of admissible $J_{s(j)}$-holomorphic disks in $(X_{s(j)},L_{s(j)})$ without interior punctures such that each $u_{s(j)}^j$ has formal dimension $0$ and if $u^j_{s(j)}$ converges to a broken curve $u_s$ in $(X_s,L_s)$, then every factor of $u_s$ is an admissible disk without interior punctures.
\end{lma}
\begin{pf}
As in Lemma \ref{l:nicecomp}, we see that, for generic $J_s$, (parameterized) moduli spaces of admissible disks with interior punctures are transversely cut out. This implies the first statement of the lemma. To see that the second statement holds, note that for each factor of $u_s$ which is a disk $v$, in $(X_s,L_s)$ or in $(Y^\pm\times\R,\Lambda^\pm\times\R)$, with $r$ interior punctures, there are $r$ $1$-punctured spheres in $X$ or in $Y^-\times\R$ which are attached to $v$ at its interior punctures. The good end assumption implies that such a configuration contributes at least $\dim(v)+2r\ge 2r-1$ to the formal dimension of $u_s$. It follows that $r=0$.
\end{pf}

Remaining parts of this section are devoted to explaining how the moduli space of admissible $J_s$-holomorphic disks in $(X_s,L_s)$ varies with $s$. More precisely, we will describe how the count of such disks, or, in the terminology of Subsection \ref{s:pot}, the potential vector $f_s$ of $(X_s,L_s)$ depends on $s$. First observe that ordinary Morse modifications of the $0$-dimensional moduli space (when $s$ passes a maximum or minimum of the $1$-dimensional parameterized moduli space) do not alter the modulo $2$ count of rigid $J_s$-holomorphic disks in the cobordism. Thus, to understand how the count changes in a $1$-parameter family, it suffices to understand how it changes as a moment where an admissible $(-1)$-disk without interior punctures and with properties as in Lemma \ref{l:classtv} is passed (such an moment will be called {\em $(-1)$-disk moment}). We therefore  specialize further and consider a $1$-parameter family $(X_s,L_s)$ with only one $(-1)$-disk moment, at $s=\frac12$. We denote the $(-1)$-disk by $\hat k$. Moreover, for convenience we will move only the complex structure and keep the Lagrangian submanifold $L_s$ fixed: if $\Phi_s\colon X\to X$ is the symplectomorphism moving $L$ and if $J_s$ is the almost complex structure then, defining $\hat J_s = d \Phi_s^{-1}\circ J_s\circ d \Phi_s$, it is easy to see that composition with $\Phi_s^{-1}$ takes $J_s$-holomorphic curves with boundary on $L_s$ to $\hat J_s$-holomorphic curves with boundary on $L_0$.

We will describe how the moduli space of rigid {\em admissible} disks at $s=0$ differs from the corresponding space at $s=1$. If the $(-1)$-disk $\hat k$ is mixed then, see Lemma \ref{l:no2mixed}, any  broken admissible disk has at most one $\hat k$-component. Because of that, it is straightforward to describe the change in the moduli space in the case when $\hat k$ is mixed, see Subsection \ref{s:mix}. In contrast to the mixed case, when the $(-1)$-disk $\hat k$ is pure, it may be used several times to build new rigid admissible disks. This creates certain transversality problems. To control what happens algebraically one must use more elaborate perturbations. Perturbation schemes to handle similar situations have been introduced by Fukaya, Oh, Ohta, and Ono, \cite{FOOO} and by  Hofer, Wysocki, and Zehnder \cite{HWZ}. We follow the latter set of authors and use their polyfold language. The perturbation scheme involves perturbations of the complex structure near the original moduli space of $J_s$-holomorphic disks. We will use the term  ``holomorphic disk'' to denote a disk which is a solution of the possibly perturbed equation.

\subsection{$M$-polyfold structure and $M$-polyfold bundles}\label{s:polyfold}
Let $(X,L)$ be an exact cobordism with good ends. In this section we will briefly discuss how to incorporate the analysis of admissible holomorphic disks in $(X,L)$ into the language of $M$-polyfolds, see \cite{H} Section 1.2. Intuitively, the $M$-polyfold ${\mathcal X}$ we will make use of is a mapping space which is a ``complete'' configuration space for admissible holomorphic disks in the following sense. The elements in this space are several level maps of punctured disks $\Delta$ into $(X,L)$ and into its ends $(Y^\pm\times\R,\Lambda^\pm\times\R)$ which match at Reeb chords. That is, each level is a collection of maps $u\colon(\Delta,\pa\Delta)\to(X,L)$ or $v\colon(\Delta,\pa\Delta)\to(Y^\pm\times\R,\Lambda^\pm\times\R)$ which are asymptotic to Reeb chord disks at their punctures and where punctures at one level matches those at the next so that by joining the disks at these punctures we obtain an admissible formal disk. In particular, a neighborhood of a broken map of several levels contains the unbroken map obtained by gluing the broken map at Reeb chords in some natural way. Over the $M$-polyfold ${\mathcal X}$ there is an $M$-polyfold bundle, see \cite{H} Section 1.2.4, ${\mathcal Y}$ and the $\bar\pa_J$-operator gives a Fredholm section, see \cite{H} Section 2.3, $\bar\pa_J\colon {\mathcal X}\to {\mathcal Y}$ such that the moduli space of admissible holomorphic disks is $\bar\pa_J^{-1}(0)$, where $0$ denotes the zero-section. In particular, this framework includes a notion of transversality between the section $\bar\pa_J({\mathcal X})$ and the $0$-section also at broken solutions, see \cite{H} Section 2.4.

We sketch how to make these notions precise in the case under study. For the sources of the maps one may use models $\Delta$ for $m$-punctured disks ($m\ge 3$) which are strips in the plane of width $m-1$ and with $m-2$ slits of small width removed, see \cite{E1} for a discussion of this way of parameterizing the space of conformal structures on the punctured disk. The $M$-polyfold chart, see \cite{H} Definition 1.10, around a smooth map $u\colon(\Delta,\pa\Delta)\to (X,L)$, which satisfies $\bar\pa_J u=0$ on $\pa\Delta$, is modeled on the space of vector fields along $u$, tangent to $L$ along the boundary, and which satisfies $\bar \nabla_J v=0$ along the boundary, where $\nabla$ is the Levi-Civita connection of a metric in which $L$ is totally geodesic  and for which $J$ takes Jacobi fields along geodesics in $L$ to Jacobi fields, see \cite{EES2,EES4}. Here we use the nested sequence of Sobolev spaces with small positive exponential weight at the punctures which increases with the number of derivatives starting with the maps with two derivatives in $L^2$. (The weights must be small compared to the complex angle of the tangent spaces to $\Lambda^\pm$ in the contact hyperplane at the Reeb chord endpoints, see Subsection \ref{s:grad} and \cite{EES1,EES3}.) Furthermore, these Sobolev spaces of vector fields should be augmented by a $1$-dimensional space of cut-off solutions supported near each puncture and corresponding to shifts in the $\R$-direction of the Reeb chord disk, much like in \cite{BM} or \cite{E1}. The $M$-polyfold bundle ${\mathcal Y}\to{\mathcal X}$ consists of the corresponding Sobolev spaces of complex anti-linear maps $T\Delta\to w^\ast(TX)$, $w\in {\mathcal X}$. At broken disks, the splicing core, see \cite{H} Section 1.2.1, is modeled on the natural gluing and anti-gluing maps, much like in Morse theory, see \cite{H} Theorem 1.21, but with the following additional properties: the subspace of the sum of cut-off solutions which contains the solutions with opposite values at matching infinities are mapped by the gluing map to a new-born conformal variation, and the subspace of cut-off solutions which vanishes at the end of the upper disk is mapped to a cut-off constant section by the anti-gluing map. The filler, see \cite{H} Section 2.2, can now be taken as the linearized $\bar\pa_J$-operator on the strip $\R\times[0,1]$ with boundary conditions at the two tangent spaces of $L$ at the ends of the Reeb chord and with positive exponential weights at the infinities acting on the usual space augmented by the cut-off solution from the anti-gluing. (This operator is an isomorphism as it should).

In order to demonstrate that the construction sketched above gives a Fredholm section of an $M$-polyfold bundle, one may build on the gluing analysis in \cite{EES1,EES3}, and \cite{E1}. In particular, the estimates derived there for the $\bar\pa_J$-operator on functional analytic spaces, with a weighted Sobolev $2$-norm, augmented by cut-off solutions in combination with the Seeley extension theorem, see \cite{S}, and ordinary boot-strap arguments for elliptic operators give the required estimates on higher derivatives.

\subsection{Polyfold-transversality for $1$-parameter families -- mixed case}\label{s:mix}
Consider a $1$-parameter family $(X_s,L_s)$, $0\le s\le 1$, of exact cobordisms as in Subsection \ref{s:mixpureini} such that there is a mixed $(-1)$-disk $\hat k$ at $s=\frac12$. Recall the increment function $\Delta\colon\V(X_0,L_0)\to\V(X_1,L_1)$ introduced in \eqref{e:dfnDelta}. If $f\in\V(X,L)$ and if $\hat k$ is mixed then Lemma \ref{l:no2mixed} implies that no admissible formal disk can have more than one $\hat k$-factor and consequently if $k\in\V(X,L)$ is the formal disk represented by $\hat k$ then $\Delta(f,k)$ simplifies and in this case
$$
\Delta(f,k)=d^f(k)=\left\{k\to h^+\right\}\!(f,f)+\left\{h^-\leftarrow k\right\}\!(f,f),
$$
where $h^\pm$ are the potential vectors at the ends of $(X_s,L_s)$.
\begin{lma}\label{l:dfmix}
Let $(X_s,L_s)$, $0\le s\le 1$, be a $1$-parameter family of cobordisms as above with a {\em mixed} $(-1)$-disk $k=\hat k$ at $s=\frac12$. Let $f_s\in\V^+(X_s,L_s)$, $s=0,1$, be the potential vectors. Then
$$
f_1=f_0+\Delta(f_0,k).
$$
\end{lma}

\begin{pf}
Consider the limit of an admissible rigid holomorphic disk in $(X_s,L_s)$ as $s\to\frac12$. If the limit does not contain a $k$-factor then the disk contributes to both $f_0$ and $f_1$. Since the disk is admissible it follows from Lemma \ref{l:dfmix} that the limit disk may contain at most one $k$-factor. Counting dimensions we find that the limit has one $h^\pm$-factor, one $k$-factor, and  remaining factors rigid factors which do not break (i.e., the remaining factors contribute to both $f_0$ and $f_1$). It follows that any disk contributing to $f_0$ and not to $f_1$, or vice versa, must have the form of a disk in $\Delta(f_0,k)$.  On the other hand, a standard gluing argument shows that any broken disk contributing to $\Delta(f_0,k)$ indeed does contribute also to $f_1+f_0$. The lemma follows.
\end{pf}

\subsection{Polyfold-transversality for $1$-parameter families -- pure case}\label{s:pure}
Consider a $1$-parameter family $(X_s,L_s)$, $0\le s\le 1$, as discussed in Subsection \ref{s:mixpureini} which has a pure $(-1)$-disk $\hat k$ at $s=\frac12$. Consider a parameterized moduli space $\M^s$ of admissible holomorphic disks with boundary on $(X_s,L_s)$. A point in this space is a pair $(u,s)$ where $s\in[0,1]$ and where $u$ is a holomorphic disk in $(X_s,L_s)$. The non-transversality at $s=\frac12$ arises as follows. Consider a broken disk with at least one factor in $(Y^\pm\times\R,\Lambda^\pm\times\R)$ and at least two $\hat k$-factors. Consider the differential of the (parameterized) $\bar\pa_{J_s}$-operator at this broken disk. At each of the factors not equal to $\hat k$, the linearized $\bar\pa_{J_s}$-operator is surjective already when restricted to a subspace which is a complement of the tangent vector $\pa_s$ of the deformation. At a $\hat k$-factor, the operator is surjective on the full source space (i.e., $\pa_s$ is included in the source space). However, independently of how many $\hat k$-factors there are, there is only one tangent vector $\pa_s$ of the deformation. Hence, at a broken disk of this type, the linearized $\bar\pa_{J_s}$-operator has a cokernel of dimension equal to the number of $\hat k$-factors minus one.

General arguments using the polyfold machinery, see \cite{H} Section 2.4 and \cite{HWZ} Theorem 6.14, show the existence of perturbations which take us out of this non-transverse situation. In fact, the argument just given shows that after such perturbations a broken holomorphic disk in a parameterized moduli space can contain at most one $(-1)$-disk factor. For our purposes, the mere existence of a perturbation is not quite enough. In order to obtain an algebraic formula for the change in potential vector, the perturbation must be carefully designed. In particular, as we shall see, while perturbing out of the non-transverse situation we create new $(-1)$-disks and the main difficulty is controlling these new $(-1)$-disks. Intuitively, the perturbation consists of separating the $s$-coordinates of the complex structures at the negative punctures of disks in $(Y^+\times\R,\Lambda^+\times\R)$.

Let the piece of the exact Lagrangian submanifold $L\subset X$ to which the pure $(-1)$-disk $\hat k$ maps the boundary be denoted $\hat L$ and let the corresponding pieces of $\Lambda^\pm$ be denoted $\hat\Lambda^\pm$. We say that a Reeb chord with two endpoints on $\hat L$ is a {\em pure $\hat L$-chord} and that a Reeb chord with only one endpoint on $\hat L$ is a {\em mixed $\hat L$-chord}.

\begin{rmk}\label{r:order}
The punctures of an admissible formal disk $w$ which maps to a pure or mixed $\hat L$-chord are ordered in a natural way as follows.
\begin{itemize}
\item[$({\bf op})$] If $w$ is a pure admissible disk then the orientation of the boundary of the disk induces an ordering of the punctures of $w$ as follows. The first puncture in the order is the first puncture in the direction of the orientation of the boundary as seen from the positive puncture, continuing in the direction of this orientation we meet all the other negative punctures in a certain order, and finally we meet the positive puncture as the last one in the order.
\item[$({\bf om})$] If $w$ is a mixed admissible disk with some boundary component mapping to $\hat L$ then the orientation of the boundary of the disk induces an ordering as follows. At exactly one puncture mapping to a mixed $\hat L$-chord the boundary orientation points into the boundary component mapping to $\hat L$, this is the first puncture. Then follows all negative punctures mapping to pure $\hat L$-chords in the order of the boundary orientation, and finally the second puncture mapping to a mixed $\hat L$-chord where the boundary orientation points out of the boundary component mapping to $\hat L$.
\end{itemize}
Furthermore, once the $\hat L$-punctures of the boundary of an admissible formal disk is ordered in this way, we order the remaining punctures in the order they appear as the boundary of the disk is traversed in the positive direction starting at the last $\hat L$-puncture.
\end{rmk}

\begin{rmk}\label{r:no3mixedpunct}
Note that the condition that the disk $w$ is admissible is used in the above construction. For example if there were a third mixed $\hat L$-puncture in $({\bf om})$ then $w$ would not be admissible. To see this let $p_1$, $p_2$, and $p_3$ be three mixed $\hat L$-punctures of a disk $D$. Choose notation so that, following the $\hat L$-boundary component near $p_1$ away from $p_1$, we first meet $p_2$ and then $p_3$. An arc connecting the $\hat L$-boundary component near $p_1$ to the $\hat L$-boundary component near $p_3$ then either separates $p_1$ from $p_2$ and $p_3$, or separates $p_1$ and $p_3$ from $p_2$. In either case the disk is non-admissible.
\end{rmk}

\begin{lma}\label{l:mixstr}
Let $w$ be a formal disk in $(X,L)$ with two admissible factors $v$ and $v'$ which are mixed disks with boundary components mapping to $\hat L$.  If $w$ is admissible then there is a unique chain of mixed disks $v_0,v_1,\dots, v_k$ in $w$, all with boundary components mapping to $\hat L$, which connects $v$ and $v'$ (i.e., $v_0=v,\,v_k=v'$ or $v_0=v',\,v_k=v$), and which is such that the first mixed $\hat L$ puncture of $v_j$ is connected to the last mixed $\hat L$-puncture of $v_{j+1}$.
\end{lma}

\begin{pf}
Straightforward from the admissibility of $w$, see Remark \ref{r:no3mixedpunct}.
\end{pf}

We next describe how to associate a vector $\delta^e=(\delta^e_1,\dots,\delta^e_q)\in[0,1]^q$ to each $u\in \M^+_{\rm pure}(e)$ where $\M^+_{\rm pure}(e)$ denotes the moduli space of pure admissible holomorphic disks in $(Y^+\times\R,\Lambda^+\times\R)$ with the property that the action $\A(c^+)$ of the Reeb chord $c^+$ at the positive puncture satisfies $\A(c^+)\le e$, and where $q$ is the number of negative punctures of $u$. (This vector will be key in the definition of our perturbation.) The space $\M^+_{\rm pure}(e)$ is a compact $M$-polyfold. The perturbation constructed using $\delta^e$ suffices to compute the change in the moduli space of disks of energy below $e$. As we shall see the change in a low energy moduli space is independent of the cut-off energy $e$ used to compute it.

Since each Reeb chord $c$ of $\hat L$ satisfies $\A(c)>\A_0>0$ and since the $\omega$-energy of any holomorphic disk is non-negative, we find that there is $N>0$ such that no disk in $\M^+_{\rm pure}(e)$ has more than $N$ negative punctures.
We let $\sigma^e\colon\M^+_{\rm pure}(e)\to[0,1]$ be a function such that $\sigma^e(u)$ depends only on the $(+)$-action $\A^+(u)$ of $u$, and such that if $u,v\in\M^+_{\rm pure}(e)$ and $\A^+(v)<\A^+(u)$ then
\begin{equation}\label{e:growth}
10N \sigma^e(v)<\frac{1}{10N}\sigma^e(u).
\end{equation}
Using a function $\R\to\R$ of sufficiently rapid growth and scaling, it is easy to construct $\sigma^e$ with this property which takes values in an arbitrarily small half-interval around $0$. We will think of the values of $\sigma^e$ as very small.

To describe the vector $\delta^e$ we first fix a small neighborhood of all broken disks in $\M^+_{\rm pure}(e)$. For $u$ outside this neighborhood we let
$$
\delta^e(u)=\sigma^e(u)\bigl((q-1),(q-2),\dots,1,0\bigr),
$$
where we think of the $j^{\rm th}$ component of the vector as belonging to the $j^{\rm th}$ negative puncture of $u$ (in the order described in $({\bf op})$). Consider next a pure broken admissible disk $v=(u_1,\dots,u_k)$ of several levels such that the vector $\delta^e(u_j)$ is defined for each $j$. Note that $v$ has the structure of a tree, where vertices correspond to its factors and edges to glued punctures. Let $U_1=u_1$ be the factor of the positive puncture of $v$. At a negative puncture $p$ of $U_1$ which is the $r^{\rm th}$ negative puncture of $u_1$, we define $\tilde\delta^e_p(U_1)=\delta^e_r(u_1)$. We proceed by induction: assume that $\tilde\delta^e_p(U_j)$ has been defined at all negative punctures $p$ of the sub-disk $U_j$ of $v$. Let $U_{j+1}$ be the sub-disk of $v$ consisting of $U_j$ and all factors $u_k$ attached to $U_j$. At any negative puncture $p$ of $U_{j+1}$ which is also a negative puncture of $U_j$ let $\tilde\delta^e_p(U_{j+1})=\tilde\delta^e_p(U_j)$. At a negative puncture $p$ of $U_{j+1}$ which is the $r^{\rm th}$ negative puncture of some factor $u_k$ attached to $U_j$ at the puncture $P$ let $\tilde\delta^e_p(U_{j+1})=\tilde\delta^e_P(U_j)+\delta^e_r(u_k)$.

Fix $k$ such that $U_k=v$. If $p$ is the $r^{\rm th}$ puncture of $v$, define $\delta^e_r(v)=\tilde\delta^e_p(U_k)$. Thus, this inductive construction gives a vector
$$
\delta^e(v)=\bigl(\delta^e_1(v),\dots,\delta^e_q(v)\bigr),
$$
which, because of \eqref{e:growth}, satisfies
\begin{equation}\label{e:order}
\delta^e_j(v)>\delta^e_{j+1}(v),\quad\text{for all $j$.}
\end{equation}
Consequently, we can extend the vector $\delta^e$ smoothly over all of $\M^+_{\rm pure}(e)$ so that \eqref{e:order} holds for all $v\in\M^+_{\rm pure}(e)$.

With the vector $\delta^e$ defined we are in position to describe our initial perturbation for {\em pure} disks in $(X,L)$. This perturbation is defined in polyfold charts near the moduli space $\M_{\rm pure}(X,L;e)$ of pure admissible holomorphic disks in $(X,L)$ of $\omega$-energy at most $e$ as follows. Consider a broken disk with a factor $u$ near $\M^+_{\rm pure}(e)$ and with several factors $v$ of broken disks mapping into $(X,L)$ and $(Y^-\times\R,\Lambda^-\times\R)$. We perturb the $\bar\pa_{J_s}$-equation so that the $s$-coordinate of $J_s$, at a disk $v$ attached to the $j^{\rm th}$ negative puncture of $u$, in the operator  $\bar\pa_{J_{s}}$, equals $s'=s_0+\delta^e_j(u)$, where $s_0$ is the  $s$-coordinate of the positive puncture of $u$, and where $\delta_j^e(u)$ is the component of the vector $\delta^e$ corresponding to the $j^{\rm th}$ negative puncture of $u$. The new perturbed operator is then extended over a neighborhood of the whole parameterized moduli space. In order to make the new operator surjective we must use an additional perturbation which can be chosen arbitrarily small. More importantly, the main perturbation we use may create new $(-1)$-disks. We study this more closely.

We start with an intuitive explanation of this phenomenon. Consider a disk $h^+\in \M^+_{\rm pure}(e)$ of dimension $1$, with two negative punctures $q_1$ and $q_2$ where the complex structures are separated in $s$, $s(q_1)<s(q_2)$. If we change the perturbation smoothly so that at the end of the change $s(q_1)>s(q_2)$ then there is some moment where $s(q_1)=s(q_2)$. At this moment we can form a new $(-1)$-disk by gluing two copies of $\hat k$ to $h^+$, at $q_1$ and $q_2$ when the $s$-parameter passes $0$. Before we introduce the perturbation described above, all negative punctures of all disks in $\M^+_{\rm pure}(e)$ have the same $s$-coordinate. Hence, we perturb from a very degenerate position and must expect to find many new $(-1)$-disks.

\begin{figure}[htbp]
\begin{center}
\includegraphics[angle=0, width=8cm]{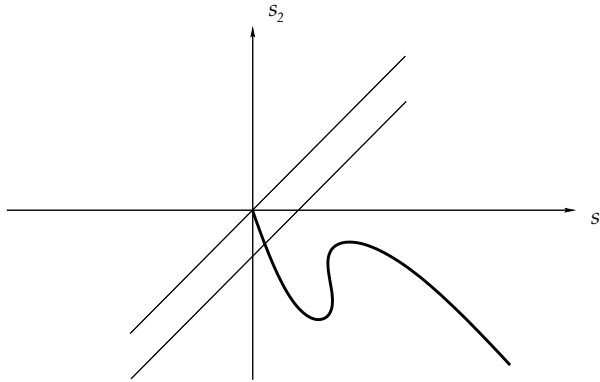}
\end{center}
\caption{Intersection between the hypersurface of $(-1)$-disks and the path of perturbed complex structures}
\label{fig:new(-1)}
\end{figure}

In order to control the new $(-1)$-disks which our perturbation gives rise to we consider first the disk $u\in\M^+_{\rm pure}(e)$ with the smallest action Reeb chord at its positive puncture. Since no other disk in $\M^+_{\rm pure}(e)$ can be glued to this disk, the only way that it can give rise to new $(-1)$-disks is when $\hat k$-factors and disks in the cobordism are attached to it. Letting the complex structures of the various $\hat k$-factors be independent we can perform this gluing. (Formally, we multiply the standard polyfold with one copy of $\R$ for each $\hat k$-factor. The differential at the broken curve is now surjective and we may use standard arguments to glue.) The dimension of the glued moduli space equals the sum of the dimensions of its factors and the number of $\R$-factors. In particular, if the dimension of the bare disk (without $\R$-factors counted) should equal $-1$ we find that the number $d$ of $\hat k$-factors is one larger than the sum of the dimensions of non-$\hat k$-factors and the standard gluing argument shows that the glued $(-1)$-disks form a manifold with boundary with corners which project to a hyper-surface in the $d$-dimensional space of complex structures near the ends, see Figure \ref{fig:new(-1)}.
We can now read off the disks that we meet: the path of complex structures the original $1$-parameter family follows is $s\mapsto (s,\dots,s)$, the perturbed $1$-parameter family corresponds to $s\mapsto\bigl(s+\delta^e_{j_1}(u), s+\delta^e_{j_2}(u),\dots, s+\delta^e_{j_d}(u)\bigr)$, where $j_1<j_2<\dots<j_d$. Perturbing it slightly it intersects the hyper-surface of $(-1)$-disks in a finite number of points at distance $c\sigma^e(u)\le {\rm distance}\le C\sigma^e(u)$ from $0$ for some constants $c\le C$. We then bring the $1$-parameter family back to the initial $1$-parameter family. The result is thus a finite number of new $(-1)$-disks with complex structures at their positive punctures distinct and in a $C \sigma^e(u)$ neighborhood of $0$. We continue in this way with the disk $u_2$ with the Reeb chord of second smallest action, gluing the original $(-1)$-disk and all $(-1)$-disks which was constructed in the previous step. The result is a finite number of new $(-1)$-disks with positive punctures separated by $C \sigma^e(u_2)$, where $\sigma^e(u_2)$ is the separation of the negative punctures for the second smallest Reeb chord. We continue in this way until we have passed through all Reeb chords of action $\le e$.

A similar argument shows that after small perturbation there are no disks of formal dimension $\le -2$ along our path: such disks correspond to a submanifold of $\R^d$ of codimension at least $2$ and after arbitrarily small perturbation, our path does not intersect a codimension two submanifold.

\begin{figure}[htbp]
\begin{center}
\includegraphics[angle=0, width=8cm]{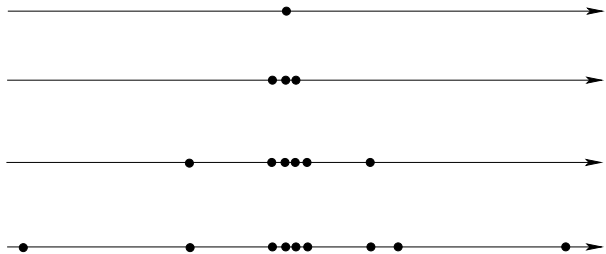}
\end{center}
\caption{Locations of positive punctures of new $(-1)$-disks after the first steps of the construction.}
\label{fig:posp(-1)}
\end{figure}

Note the following important feature of the perturbation: no $(-1)$-disk constructed in step $r$ can be glued to any negative puncture of a $\M^+_{\rm pure}(e)$-disk used in step $s$, if $s<r$. This follows since the Reeb chord at the positive puncture of the step $r$ $(-1)$-disk has larger action than the Reeb chord at the positive puncture of the $\M^+_{\rm pure}(e)$-disk in step $s$ which in turn has larger action than any of its negative puncture Reeb chords. As mentioned above the perturbation size grows rapidly with Reeb chord length and we can depict the location of the positive punctures of the new-born Reeb chord disks as in Figure \ref{fig:posp(-1)}.

\begin{figure}[htbp]
\begin{center}
\includegraphics[angle=0, width=8cm]{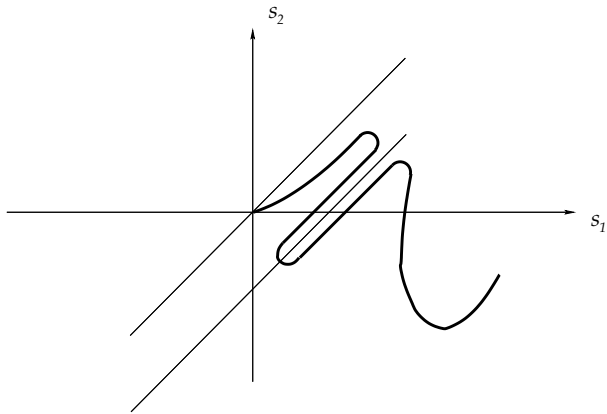}
\end{center}
\caption{The result of pushing the hyper-surface of $(-1)$-disks.}
\label{fig:finger}
\end{figure}

In order to finish the construction of the perturbation we must consider also mixed disks with some boundary component mapping to $\hat L$. The construction is similar to the one just discussed. Note that any mixed admissible disk in $(Y^+\times\R,\Lambda^+\times\R)$ has exactly two mixed $\hat \Lambda^+$-punctures and that all pure $\hat \Lambda^+$-punctures lie between these. We use a perturbation near mixed disks with some boundary component mapping to $\hat \Lambda^+\times\R$ much like the one described above for pure disks but which separates the punctures in mixed disks by an amount much larger than the separation between the punctures of pure disks. The first $\hat \Lambda^+$-puncture will have the largest $s$-coordinate, the last and all non $\hat \Lambda^+$-punctures the smallest. It follows from Lemma \ref{l:mixstr} that such a perturbation can be extended to a neighborhood of the moduli space. We must control the new $(-1)$-disks that the perturbation gives rise to. Arguing as above we find that each perturbed moduli space component together with $(-1)$-disks from earlier steps gives a hyper-surface in the space of complex structures and that the intersection points of the shifted diagonal line with this hyper-surface determines the new $(-1)$-disks in the perturbed $1$-parameter family. At this point we use an auxiliary perturbation which moves the intersections into the region where the complex structure of the first puncture lies well inside the region $\{s>0\}$ and the complex structure of the last puncture lies well inside $\{s<0\}$. (Here ``well inside'' refers to a big amount compared to earlier perturbations.) This can be obtained by pushing the hypersurface of $(-1)$-disks locally along the chosen path (and changing the differential equation accordingly), see Figure \ref{fig:finger}. Thus, all newly created $(-1)$-disks have their first puncture to the far right of $0$ and their last to the far left of $0$. As above the construction is inductive and we consider next the moduli space of mixed disks of the next action. Note that only one other mixed $(-1)$-disk can be glued to such a disks if it is to be admissible and that the structure of first/last punctures is preserved: first punctures stay to the right of $0$ and last punctures to the left.

To conclude, we have found a perturbation of the $\bar\pa_{J_s}$-equation by ordering the punctures of disks in $(Y^+\times\R,\Lambda^+\times\R)$ in the deformation parameter $s$ according to the order induced by the orientation of the boundary of the disk as explained in $({\bf op})$ and $({\bf om})$. In this way we created new $(-1)$-disks: all $(-1)$-disks in previous steps give rise to new $(-1)$-disks in later steps. However, the complex structures at the positive- and at the mixed punctures of the new $(-1)$-disks are separated in $s$.

\begin{rmk}\label{r:(-1)hatLcomp}
All new admissible $(-1)$-disks which arise in the above perturbation scheme has some boundary component mapping to $\hat L$. To see this we note that this holds true for the disk $\hat k$ and assume inductively that it holds true for all $(-1)$-disks arising after $r$ steps in the perturbation scheme. Since the $(-1)$-disks which arise in the $(r+1)^{\rm th}$ step contains at least one $(-1)$-disk that arouse in some previous step it follows that all $(-1)$-disk which arise have some boundary component mapping to $\hat L$.
\end{rmk}

\begin{lma}\label{l:polygen}
The parameterized moduli spaces of admissible holomorphic disks in $(X_s,L_s)$ of dimensions $\le 1$ for the perturbation described above are being transversally cut out.
\end{lma}

\begin{pf}
The fact that no broken admissible disk can have two copies of the same mixed $(-1)$-disk follows as before. A mixed disk in $(Y^-\times\R,\Lambda^-\times\R)$ cannot glue to two mixed $(-1)$-disks since its first and last punctures have complex structure at the same $s$-coordinate. Any disk in $(Y^-\times\R,\Lambda^-\times\R)$ which is attached to a pure $(-1)$-disk cannot be attached to any other disk since its positive puncture is pure. We show that no disk in $(Y^+\times\R,\Lambda^+\times\R)$ has two negative punctures attached to $(-1)$-disks. If the disk $u$ is pure then the $(-1)$-disks which could be glued to it have Reeb chords of smaller action than $\A^+(u)$. Hence on the scale of the $s$-separation of the negative punctures of $u$ the complex structures at the positive punctures of the $(-1)$-disks which could be glued are at approximately the same $s$-value. The same holds for pure punctures of a mixed disks $u$. Hence, the $s$-value of at most one pure negative puncture of $v$ at a time matches some $(-1)$-disk. If one of the two negative $\hat L$-punctures of $u$ is mixed (first or last), the fact that it matches the mixed puncture (last or first) of a mixed $(-1)$-disk implies that all its pure negative $\hat L$-punctures and the second mixed $\hat L$-puncture are outside the region where the pure $(-1)$-disks have their positive punctures. We conclude that any broken disk has at most one $(-1)$-disk factor and it is easy to see that this implies that the linearized operator is surjective.
\end{pf}

Recall the increment function $\Delta\colon\V(X_0,L_0)\to\V(X_1,L_1)$ introduced in \eqref{e:dfnDelta}. Let $k\in\V^+(X,L)$ denote the vector obtained by counting all $(-1)$-disks that appear in the construction of the perturbation. (Looking at intersections between the perturbed path and the hyper-surface of $(-1)$-disks, see Figure \ref{fig:new(-1)}, it is straightforward to check that the part of the vector $k$ of $\omega$-energy $\le e_0$ is independent of the cut-off $\omega$-energy level $e>e_0$ used to construct the perturbation for all sufficiently small perturbations.)

\begin{lma}\label{l:dfpure}
The potential vectors $f_s\in\V^+(X_s,L_s)$, $s=0,1$, are related as follows
$$
f_1 =f_0+ \Delta(f_0,k).
$$
\end{lma}

\begin{pf}
Consider the contribution to $f_0$ and $f_1$ from disks below some fixed $\omega$-energy level. After Lemma \ref{l:polygen} we know that the parameterized moduli space of holomorphic disks of dimension $0$ (parameterized dimension $1$) is a compact $1$-manifold. We check that $f_1+f_0+\Delta(f_0,k)$ counts the points in its boundary.

Consider first gluing of a disk $h^+$ of dimension $1$ in $(Y^+\times\R,\Lambda^+\times\R)$ to a $(-1)$-disk. If the $h^+$-factor is pure then on the scale of the $s$-separation of its negative punctures, the $k$-factors which can be glued to $h^+$ have complex structures at their positive punctures which all lie very close to each other and very close to $0$. Thus, we can glue only one such $k$-factor at a time. Furthermore, the negative punctures of the $h^+$-factor before an attached $k$-factor lie in the region where moduli spaces of rigid disk with positive punctures at Reeb chords of action less than the action of the positive Reeb chord of $h^+$ have already changed, whereas the punctures after the attached $k$-factor lie in the region where this has not yet happened. A mixed $h^+$-factor which is glued to a $(-1)$-disk at a pure puncture has the properties just mentioned. If its first puncture hits a $(-1)$-disk then only rigid disks which have not yet changed can be glued to its other negative punctures. If its last puncture hits a $(-1)$-disk then only rigid disks which have already changed can be glued to other punctures.

Finally, we consider a $1$-dimensional disk $h^-$ in $(Y^-\times\R,\Lambda^-\times\R)$. For pure $h^-$ there is nothing to prove, they can join to only one $k$-factor. For mixed $h^-$ glued at its first (last) puncture to a $k$-factor only rigid disks which have not (have) changed can be glued at other punctures.

It is straightforward to see that the description given here of the part of the boundary of the parameterized moduli space which does not come from $f_0+f_1$ equals the vector $\Delta(f_0,k)$ constructed in \eqref{e:dfnDelta}.
\end{pf}


\begin{thebibliography}{999}

\bibitem{B}
F. Bourgeois,
{\em A Morse-Bott approach to Contact Homology},
PhD thesis, Stanford University (2002)

\bibitem{BEE}
F. Bourgeois, T. Ekholm, Y. Eliashberg,
{\em A surgery exact sequence for linearized contact homology},
in preparation.

\bibitem{BEHWZ}
F. Bourgeois, Y. Eliashberg, H. Hofer, K. Wysocki, E. Zehnder,
{\em Compactness results in symplectic field theory},
Geom. Topol. {\bf 7} (2003), 799--888

\bibitem{BM}
F. Bourgeois, K. Mohnke,
{\em Coherent orientations in symplectic field theory},
Math. Z. {\bf 248} (2004), no. 1, 123--146


\bibitem{Ch}
Y. Chekanov,
{\em Differential algebra of Legendrian links},
Invent. Math. {\bf 150} (2002), no. 3, 441--483

\bibitem{CL}
O. Cornea, F. Lalonde,
{\em Cluster homology: an overview of the construction and results},
Electron. Res. Announc. Amer. Math. Soc. {\bf 12} (2006), 1--12


\bibitem{E1}
T. Ekholm,
{\em Morse flow trees and Legendrian contact homology in 1-jet spaces}, math.SG/0509386

\bibitem{E2}
T. Ekholm,
{\em Legendrian $2n$-spheres distinguished by rational SFT},
in preparation

\bibitem{EENS}
T. Ekholm, J. Etnyre, L. Ng, M. Sullivan,
{\em Knot contact homology},
in preparation

\bibitem{EES1}
T. Ekholm, J. Etnyre, M. Sullivan,
{\em Non-isotopic Legendrian submanifolds in ${\mathbb R}\sp {2n+1}$},
J. Differential Geom. {\bf 71} (2005), no. 1, 85--128

\bibitem{EES2}
T. Ekholm, J. Etnyre, M. Sullivan
{\em The contact homology of Legendrian submanifolds in ${\mathbb R}\sp {2n+1}$}, J. Differential Geom. {\bf 71} (2005), no. 2, 177--305

\bibitem{EES3}
T. Ekholm, J. Etnyre and M. Sullivan
{\em Orientations in Legendrian contact homology and exact Lagrangian immersions}, Internat. J. Math. {\bf 16} (2005), no. 5, 453--532

\bibitem{EES4}
T. Ekholm, J. Etnyre and M. Sullivan,
{\em Legendrian Contact Homology in $P\times\R$},
to appear in Trans. Amer. Math. Soc, math.SG/0505451


\bibitem{El}
Y. Eliashberg,
{\em Invariants in contact topology},
Proceedings of the International Congress of Mathematicians, Vol. II (Berlin, 1998). Doc. Math. 1998, Extra Vol. II, 327--338


\bibitem{EGH}
Y. Eliashberg, A. Givental, H. Hofer,
{\em Introduction to symplectic field theory},
GAFA 2000 (Tel Aviv, 1999). Geom. Funct. Anal. 2000, Special Volume, Part II, 560--673

\bibitem{FO}
K. Fukaya, Y.-G. Oh,
{\em Zero-loop open strings in the cotangent bundle and Morse homotopy},
Asian J. Math. {\bf 1} (1997), no. 1, 96--180.


\bibitem{FOOO}
K. Fukaya, Y.-G. Oh, H. Ohta, K. Ono,
{\em Lagrangian intersection Floer theory - anomaly and obstructon},
preprint.

\bibitem{G}
H. Geiges,
{\em Contact geometry},
Handbook of differential geometry. Vol. II, 315--382, Elsevier/North-Holland, Amsterdam, 2006.


\bibitem{H}
H. Hofer,
{\em A General Fredholm Theory and Applications},
math.SG/0509366

\bibitem{HWZ}
H. Hofer, K. Wysocki, E. Zehnder,
{\em Polyfolds and Fredholm Theory I},
in preparation.

\bibitem{K}
T. Kalman, {\em Contact homology and one parameter families of Legendrian knots}, Geom. Topol. {\bf 9} (2005), 2013--2078

\bibitem{Mi}
J. Milnor,
{\em On axiomatic homology theory},
Pacific J. Math. {\bf 12} 1962 337--341.

\bibitem{Ng1}
L. Ng, {\em Knot and braid invariants from contact homology. I.},
Geom. Topol. {\bf 9} (2005), 247--297

\bibitem{Ng2}
L. Ng,
{\em Knot and braid invariants from contact homology. II. With an appendix by the author and Siddhartha Gadgil},
Geom. Topol. 9 (2005), 1603--1637



\bibitem{S}
R. Seeley, {\em Extension of $C\sp{\infty }$ functions defined in a half space}, Proc. Amer. Math. Soc. {\bf 15} (1964) 625--626
\end{thebibliography}
\end{document}